\documentclass[onefignum,onetabnum]{siamonline220329}

\usepackage{lipsum}
\usepackage{amsfonts}
\usepackage{graphicx}
\usepackage{epstopdf}
\usepackage{algorithmic}
\ifpdf
  \DeclareGraphicsExtensions{.eps,.pdf,.png,.jpg}
\else
  \DeclareGraphicsExtensions{.eps}
\fi

\usepackage{enumitem}
\setlist[enumerate]{leftmargin=.5in}
\setlist[itemize]{leftmargin=.5in}

\newsiamremark{remark}{Remark}
\newsiamremark{hypothesis}{Hypothesis}
\crefname{hypothesis}{Hypothesis}{Hypotheses}
\newsiamthm{claim}{Claim}

\headers{Kernel EDMD For Sparse Reconstruction of Closable Koopman Operators}{Himanshu S. \& Nishant P. \& J. Nathan K.}
\title{
Kernel Dynamic Mode Decomposition For Sparse Reconstruction of Closable Koopman Operators
\thanks{Submitted to the editors DATE.
\funding{This work (LA-UR-25-24301) is funded by Los Alamos National Laboratory (LANL)-LDRD.}}}

\author{Himanshu Singh\thanks{ Information Sciences, CCS-3, LANL, Los Alamos, NM 87545, USA
  (\email{singh\_h@lanl.gov}).}
\and Nishant Panda\thanks{(Corresponding Author) Information Sciences, CCS-3, LANL, Los Alamos, NM 87545, USA 
  (\email{nishpan@lanl.gov}).}
\and J. Nathan Kutz\thanks{Applied Mathematics \& Electrical and Computer Engineering, University of Washington, Seattle, WA 98195, USA
  (\email{kutz@uw.edu}).}
}

\usepackage{amsopn}

\ifpdf
\hypersetup{
  pdftitle={Kernel Dynamic Mode Decomposition For Sparse Reconstruction of Closable Koopman Operators}
  pdfauthor={Himanshu Singh and Nishant Panda and J. Nathan Kutz}
}
\fi

\usepackage{cite,float,nicefrac}
\newtheorem{assumption}[theorem]{Assumption}
\newtheorem{example}[theorem]{Example}
\newtheorem{aim}[theorem]{Experimental Aim}

\newtheorem{details}[theorem]{Experimental Details}
\newtheorem{datasetinfo}[theorem]{Data set Information}
\newtheorem{note}[theorem]{Note}
\newtheorem{datamatrix}[theorem]{Data Matrix}

\usepackage{bm,mathtools,dsfont,nicefrac}
\newcommand\numberthis{\addtocounter{equation}{1}\tag{\theequation}}
\newcommand{\vx}{\bm{x}}
\newcommand{\vy}{\bm{y}}
\newcommand{\vb}{\bm{b}}

\newcommand{\Cn}{\mathbb{C}^D}
\newcommand{\mI}{\mathbf{I}}
\newcommand{\mA}{\mathbf{A}}
\newcommand{\bmz}{\bm{z}}

\newcommand{\Hl}{H_{\mathds{L}}}

\newcommand{\koop}{\mathcal{K}_{\varphi}}
\newcommand{\Lap}{K_{{exp}}^{1,1}}
\newcommand{\Lapsigma}{K_{{exp}}^{1,\sigma}}
\newcommand{\GRBF}{K_{{exp}}^{2,1}}
\newcommand{\GRBFsigma}{K_{{exp}}^{2,\sigma}}
\newcommand{\LOmega}{L^2(\Omega)}

\usepackage{amssymb}
\usepackage{mathrsfs}

\usepackage{tabularx,booktabs}
\newcolumntype{C}{>{\centering\arraybackslash}X} 
\usepackage{longtable}
\usepackage{booktabs, siunitx}
\newcolumntype{M}[1]{>{\centering\arraybackslash}m{#1}}
\usepackage{subcaption,float}
\usepackage[font=footnotesize,labelfont=bf]{caption}
\usepackage{pgfplots}
\usepgfplotslibrary{patchplots}
\usetikzlibrary{patterns, positioning, arrows}
\pgfplotsset{compat=1.15}
\usepackage[compact]{titlesec} 
\begin{document}

\maketitle

\begin{abstract}
Spatial temporal reconstruction of dynamical system is indeed a crucial problem with diverse applications ranging from climate modeling to numerous chaotic and physical processes. These reconstructions are based on the harmonious relationship between the Koopman operators and the choice of dictionary, determined implicitly by a kernel function. This leads to the approximation of the Koopman operators in a reproducing kernel Hilbert space (RKHS) associated with that kernel function. Data-driven analysis of Koopman operators demands that Koopman operators be closable over the underlying RKHS, which still remains an unsettled, unexplored, and critical operator-theoretic challenge. We aim to address this challenge by investigating the embedding of the Laplacian kernel in the measure-theoretic sense, giving rise to a rich enough RKHS to settle the closability of the Koopman operators. We leverage Kernel Extended Dynamic Mode Decomposition with the Laplacian kernel to reconstruct the dominant spatial temporal modes of various diverse dynamical systems. After empirical demonstration, we concrete such results by providing the theoretical justification leveraging the closability of the Koopman operators on the RKHS generated by the Laplacian kernel on the avenues of Koopman mode decomposition and the Koopman spectral measure. Such results were explored from both grounds of operator theory and data-driven science, thus making the Laplacian kernel a robust choice for spatial-temporal reconstruction. 
\end{abstract}

\begin{keywords}
Earth Science, Dynamic Mode Decomposition, Sparse Reconstruction,
Closable Koopman Operators, $0$\textsuperscript{th} Snapshot $m$\textsuperscript{th} Koopman Mode Difference
\end{keywords}

\begin{MSCcodes}
86A08, 47N60, 47N70, 46E22, 47B32, 46E20, 46E22, 70G60, 76F20
\end{MSCcodes}

\section{Introduction}
Spatial temporal (ST) data have become ubiquitous in the current age where spatial-temporal sensors are storing large volumes of data in almost every geophysical avenue. An important aspect of dealing with ST data is the task of prediction based on the data collection with sparsity and non-uniformly (irregular). Numerous applications in traffic prediction, video action recognition, and anomaly detection have led to advancement in deep learning techniques such as spatial-temporal convolution, spatial-temporal transformers, graph-based spatial-temporal transformers, to name a few~(see, e.g.~\cite{STCNN,STCONV,STForecast,STGraph,STtransformer,STuniformer}). On the other hand, there is an increased focus on understanding ST data from a dynamical systems perspective (see~\cite{budivsic2012applied}) where the task is not ST prediction, but rather \emph{spatial-temporal reconstruction}. 
Currently, data-driven methods provide a promising route in ST reconstruction by learning \emph{spatial-temporal modes} of the dynamical system. One of the key and useful data-driven algorithms for providing ST modal information is \emph{Dynamic Mode Decomposition} developed by \cite{schmid2010dynamic}, which in principle is an unsupervised machine learning (ML) algorithm \cite{fujii2019supervised} based on the Krylov subspaces \cite{saad2003iterative} and the Arnoldi algorithm \cite{arnoldi1951principle}.

\begin{figure}[H]
    \centering
    \frame{\includegraphics[scale=.45]{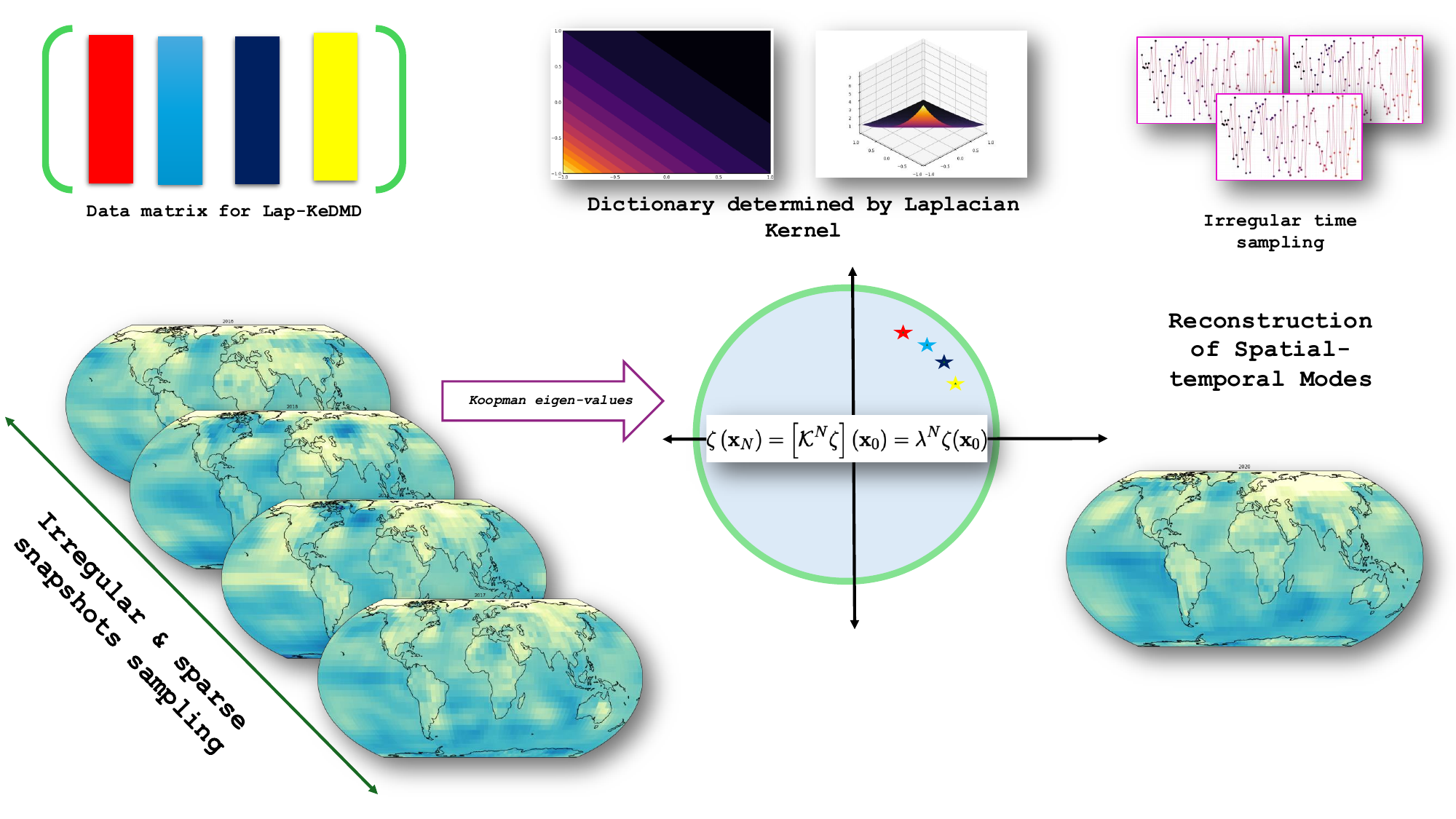}}
    \frame{\includegraphics[scale=.7]{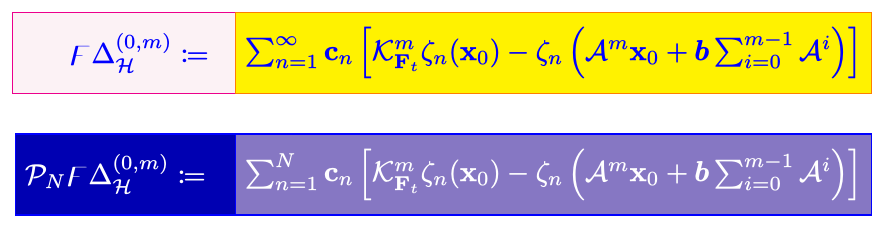}}
    \footnotesize
    \caption{The paradigm of interpret-ability of a dynamical system through Kernel Extended Dynamic Mode Decompositions. 
    We will incorporate the Laplacian kernel to overcome our \Cref{aim} whose Koopman operator theoretic inference is given in the last couple of equations capturing Koopman mode decomposition difference.}
    \label{fig:spectral_0thmthsnapshot}
\end{figure}

\begin{aim}\label{aim}
    We are interested in learning to recover the spatial-temporal reconstruction for a given data set when only a limited number of data snapshots are provided, which are also collected in an arbitrary manner with respect to time; illustrated in \Cref{fig:spectral_0thmthsnapshot}. 
\end{aim}

As far as ST reconstruction is concerned, ML-centered architectures \cite{rubanova2019latent} (ODE-RNNs), \cite{ye2020flow,murata2020nonlinear,hasegawa2020cnn} and \cite[Chapter 6 Section 8, Page 236]{brunton2022data} based on Convolutional Neural Networks and Autoencoders offer a competitive platform to simulate, model, and forecast complex, chaotic nonlinear dynamical systems, but face serious issues in analyzing the underlying dynamics. This is because they lack a tendency to operate or evolve with respect to \emph{time} (cf. \cite{mezic2021koopman,sarker2021deep,haggerty2023control}) followed by random initialization in their down-streaming (fine-tuning) task as well. In particular, such a pitfall is overcome by data-driven techniques powered by the operator theoretic and spectral analysis of infinite-dimensional linear Koopman operators (cf. \cite{koopman1931hamiltonian,koopman1932dynamical,mezic2005spectral}) over the underlying reproducing kernel Hilbert space (RKHS) \cite{williams2015data,li2017extended}, and such a DMD framework is referred to as \emph{Extended DMD} (eDMD). Approximation of Kooopman operators also crop up in applications like data assimilation~\cite{panda2021data} and uncertainty quantification~\cite{panda2025local}.

The appropriately chosen underlying RKHS for the robust operator-theoretic action of Koopman operators to perform eDMD naturally becomes a justifiable option, which is called \emph{Kernel Extended DMD} (KeDMD) \cite{baddoo2022kernel,klus2018kernel,klus2020kernel,rosenfeld2022dynamic} and is central to this paper. Gaussian Radial Basis Function (GRBF) kernel \cite[Page 4639]{steinwart2006explicit}, \cite{singh2024machine} is one of the most popular choices of kernel functions, which finds application not only in data-driven modeling \cite{degennaro2019scalable}, but also in various ML platforms such as support vector machines \cite{steinwart2006explicit} etc. 

In general, an arbitrary RKHS fails to exhibit invariance under the action of the Koopman operator. Therefore, the choice of kernel function and its corresponding RKHS plays a pivotal role. Recent review in the scope of KeDMD in \cite[Page 44]{colbrook2023multiverse,COLBROOK2024127} recommends that one should choose the kernel function space so that the Koopman operator on the underlying RKHS is not just densely defined but also `\emph{closable}'; however, coming up with such kernel spaces is non-trivial \cite{ikeda2022koopman,ikeda2022boundedness}. 

Further, the problem of ST-reconstruction for data predictions while collecting data snapshots irregularly is still under-investigated. Much of data decomposition techniques rely on the collection of data snapshots with regular and non-sparse time sampling, which is further accompanied by having full-rank data matrix. Ideally, \emph{full dataset} is needed to initiate DMD and is synthesis by stacking the data snapshots collected in a regular time-sampled manner. The experimental analysis in \cite{bevanda2024koopman,tu2014dynamic,klus2016numerical,rosenfeld2022dynamic,brunton2022data,williams2015data,schmid2010dynamic,schmid2022dynamic,baddoo2022kernel} (and references therein as well) seems to enjoy such a hypothetical ideal situation, which seems to be quite far from the actual reality. 
\subsection{Related Work}
Notable work among various others in this particular direction includes the work of \cite{freeman2023data}. Here, they utilized the quantum operation induced by the Koopman operator of the dynamical system allowing the Bayesian observational update in the operator algebra setting \cite{takesaki2003theory}. Nevertheless, describing and forecasting the time evolution of dynamical systems still remains a challenging problem in the context of limited irregularly placed data snapshots. Lack of investigation in this direction motivates this paper, which provides the solution to this problem by leveraging kernel-extended DMD and the novelty of the Laplacian kernel.

\subsection{Our contribution}
This paper aims to provide a solution to the challenge of recovering data dynamics when we have irregular and sparse data together with a limited dataset. Here, recovering data dynamics implies having the reconstruction of data or its ST-features by dominant Koopman modes determined over the RKHS that arises through the Laplacian kernel as a $L^2-$measure.
We now itemize our main contribution as follows.
\begin{enumerate}
\item We provide the KeDMD algorithm in \Cref{section_KeDMD} powered by the application of data-driven Koopman operators and the Laplacian kernel, which we refer to as Lap-KeDMD to recover the ST modal information from irregular and sparsely sampled data. We considered \emph{seven data sets} (see \Cref{tab:dataset}) to test the effectiveness of the Lap-KeDMD algorithm and compared it with the standard GRBF kernel, which is the standard benchmark opponent of the Laplacian kernel. 
\item
We develop an RKHS by the Laplacian kernel in \Cref{subsec_LapKernelL2} embedded as the Lebesgue measure $L^2$ by generating the orthonormal basis and then applying the Moore-Aronszajn Theorem in \Cref{theorem_RKviaOrthonormalbasis}.
    \footnote{Once the RKHS is constructed then we provide the Koopman operator theoretic quantification (boundedness and then compactness)
    over this in the \textbf{Supplementary Materials}. Achieving the compactness of the Koopman operator over this RKHS is crucial because it allows us to represent it as a finite rank representation of it for its approximation, which we further use in the ST reconstruction.} 
    We show that the Koopman operators over this RKHS space is \emph{closable} in \Cref{closable_Koopman_theorem}
    while showing that the current GRBF kernel leads to a failure of closability in \Cref{theorem_GRBF_Koop_closable_failure}. 
    \item Lastly, we provide theoretical justification that the ST reconstruction via the Laplacian kernel is faithful according to the definition given in \Cref{def_0msnapdiff}, in contrast to the application of the GRBF kernel in \Cref{theorem_spectralapprox}. 
    After establishing this result, we also show that the Koopman mode decomposition needed for the ST reconstruction then aligns with that performed on the RKHS of the Laplacian kernel in \cref{theorem_partialsum}.
\begin{table}[htbp]
\footnotesize
\caption{Dataset consider to test the Lap-KeDMD algorithm for their ST reconstruction.}\label{tab:dataset}
\begin{center}
  \begin{tabular}{c|c|c|c} \hline
   \# &  Dataset & Dataset's nature & References \\ \hline
    1 & Burger's Equation & Non-linear PDE & {\cite{brunton2022data}} \\
    2 & Fluid flow past cylinder & Navier-Stokes PDE  & \cite{brunton2022data}\\ 
    3 & Duffing's Oscillator & Chaotic (2-D) ODE & \cite{colbrook2023multiverse} \\ 
    4 & Seattle I-5 freeway traffic data & {\tt{Non-deterministic System}} & {\cite{chen2021low}} \\ 
    5 & 1963 Lorenz Attractor & Chaotic Non-linear (3-D) ODE  &\cite{lorenz1963deterministic,brunton2017chaos}\\
    6 & Rössler Attractor & Chaotic Non-linear (3-D) ODE & \cite{rossler1976equation,brunton2017chaos} \\
    7 & NOAA Sea Surface Temp. Anomaly & {\tt{Non-deterministic System}} & \cite{zhangNOAA}\\
    \hline
  \end{tabular}
\end{center}
\end{table}
\end{enumerate}
\Cref{tab:dataset_sensors} provides the number of spatial and temporal sensors for the dataset in \Cref{tab:dataset}.
\begin{table}[H]
\footnotesize
\caption{Spatial and temporal sensors information for the datasets in \Cref{tab:dataset}.}\label{tab:dataset_sensors}
\begin{center}
  \begin{tabular}{c|c|c|c|c} \hline
   \# &  Dataset & \# Spatial Sensors & \# Temporal Sensors & Total data-values \\ \hline
    1 & Burger's Equation & $256$ & $101$ & $25,856$\\
    2 & Fluid flow past cylinder & $89,351$ & $151$ & $13,492,001$\\ 
    3 & Duffing's Oscillator & $2$  & $5\times10^4$ & $1\times10^5$ \\ 
    4 & Seattle I-5 freeway traffic data & $75$ & $72$ & $5,400$\\ 
    5 & 1963 Lorenz Attractor & $3$  & $2\times10^5$ & $6\times10^5$\\
    6 & Rössler Attractor&$3$ & $6.4\times10^4$ & $192,000$\\
    7 & NOAA Sea Surface Temp. Anomaly & $2,592$ & $2,101$ & $5,445,792$\\
    \hline
  \end{tabular}
\end{center}
\end{table}
\begin{note}
    We provide only direct results pertaining to the execution of Kernel Extended Dynamic Mode Decomposition and also the compactness characterization of the Koopman operators of RKHS yielded by the Laplacian kernel as a $L^2-$measure. All other helper results and theorem are provided with details in the supplement part of this article.
\end{note}
\subsection{Notations}
Most of the mathematical symbols in this paper are standard; however, we still provide a tabular description of the set of notation in \Cref{tab:notation}.
\begin{table}[htbp]
\footnotesize\caption{List of frequently used mathematical symbols.}
    \centering
    \begin{tabular}{c|c}
        \textsc{Symbols} & \textsc{Meaning} \\\hline
        $X,\mathcal{X},\mathbb{X}$& Non-empty input space\\
        $f_1\otimes f_2(x,x')$& Tensor product notation i.e. $\coloneqq f_1(x)f_2(x')$\\
        $\mathbb{R,C,W}$&Field of real, complex and whole numbers respectively\\
        $\mathcal{K}_{\mathbf{F}_t}$ & Infinite dimensional linear Koopman Operator  \\
        $\mathcal{K}$ & Finite rank representation of $\mathcal{K}_{\mathbf{F}_t}$ (short notation of $\mathcal{K}_{\mathbf{F}_t}$)\\
        $\mathcal{K}_{\varphi}$ & Koopman operator induced by multi-variable complex holomorphic function $\varphi$\\
        $\mathcal{K}_{\varphi}|_{H}$ & Restriction of $\mathcal{K}_{\varphi}$ on $H$\\
        $\mathbf{K}$& Single dimension field space such as $\mathbb{R}$ or $\mathbb{C}$\\
        $\mI_D$&Identity Matrix of dimension $D\times D$\\
        $\vx, \left(\vx^\top\right)$& Real or Complex Vector (transpose)\\
        $\langle\vx,\bmz\rangle$&Dot product of vectors, i.e. $\langle\vx,\bmz\rangle=\vx\cdot\bmz^*=\vx\cdot\overline{\bmz}^\top=\sum_ix_i\overline{z}_i$\\
        $\langle f,g\rangle_{\mathcal{H}}$&RKHS $\mathcal{H}$ Inner product between $f$ and $g$\\
        $\|\vx\|_2$&Euclidean length of vector, i.e. $\left(\sum_i|x_i|^2\right)^{1/2}$\\
        $ K_{\text{exp}}^{\gamma,\sigma}(\bm{x},\bm{z})$& Exponential Power Kernel with parameters $\gamma$ and $\sigma$, both $>0$\\
        $dA,dV_D$&Lebesgue area measure, $D-$dimensional Lebesgue volume measure respectively\\
        $\delta_{NM}$&Kronecker delta function, i.e. $\delta_{NM}=1$ if $N=M$ and $\delta_{NM}=0$ if $N\neq M$\\
        $\Gamma(z)$&Gamma Function, i.e. $\Gamma(z)=\int_0^\infty t^{z-1}e^{-t}dt$\\
    \end{tabular}
    \label{tab:notation}
\end{table}
\section{Background}
In reconstruction of ST data, we have the Sampling-flow assumption that there is an underlying dynamical system and the data snapshots are a noisy observable of that dynamical system collected in certain sense. This particular assumption is common and standard to initiate the data-driven methods discussion \cite{bevanda2024koopman,brunton2022data,colbrook2023multiverse,giannakis2020extraction} and etc. 
\begin{assumption}[Sampling-flow assumption]\label{assumption_1} Let $\Omega\in\Cn$ be a metric space and $\mathbf{F}_t:\Omega\to\Omega$ be the flow as defined in \eqref{eq_1} along with the Borel-probability measure $\mu$ whose support $\operatorname{supp}\mu=\Omega$. Let the system be sampled at a fixed-time-instant, say $\Delta t~(>0)$ such that $\mathbf{F}_{N\Delta t}:\Omega\to \Omega$:
\begin{align}\label{eq_1}
\frac{d}{dt}{\bm{x}}(t)=\mathbf{f}\left({\vx}(t)\right)
    \implies\mathbf{F}_t\left({\vx}(t_0)\right)={\vx}\left(t_0\right)+\int_{t_0}^{t_0+t}\mathbf{f}\left({\vx}(\mathfrak{t})\right)d\mathfrak{t}.
\end{align}
\end{assumption}
For system \eqref{eq_1}, we can generate (finite) number of snapshots $\left\{\mathbf{x}_i\right\}_{i=0}^N\in\Omega$ that satisfy \eqref{eq_1}.
\begin{definition}[Koopman Operators]\label{def_KoopmanOperators}
Under \cref{assumption_1}, 
the dynamical flow $\mathbf{F}_t$ induces a linear map $\mathcal{K}_{\mathbf{F}_t}:\mathcal{F}\to\mathcal{F}$ 
defined as 
\begin{align}\label{eq_9r}
    \mathcal{K}_{\mathbf{F}_t}g\coloneqq g\circ\mathbf{F}_t
\end{align}
where $\mathcal{F}$ is the observable function space  invariant under the action of $\mathcal{K}_{\mathbf{F}_t}$. 
\end{definition}
In principle, we associate with $\mathcal{F}\coloneqq L^2(\Omega,\mu)$, where $\mu$ is a Borel measure, which may not necessarily be finite in whole $\Omega$. We often suppress the Borel measure $\mu$ when discussing the obvious observable space $L^2(\Omega,\mu)$ and use $\LOmega$ to denote the same. Equipping with a Borel measure $\mu$, then one can have following for a Borel measurable set $\mathds{B}\subset\Omega$
\begin{align}\label{eq_BorelSet}
    \mu\left(\mathds{B}\right)=\mu\left(\mathbf{F}{_t}^{-1}(\mathds{B})\right),\quad\text{where}\quad\mathbf{F}_t^{-1}(\mathds{B})\coloneqq\left\{\mathbf{x}:\mathbf{F}_t(\mathbf{x})\in\mathds{B}\right\}.
\end{align}
These details are crucial to initiate the decomposition of the Koopman mode after invoking the spectral theory of the unitary operator \cite{hall2013quantum}.

\begin{figure}[H]
    \centering
    \includegraphics[width=0.75\linewidth]{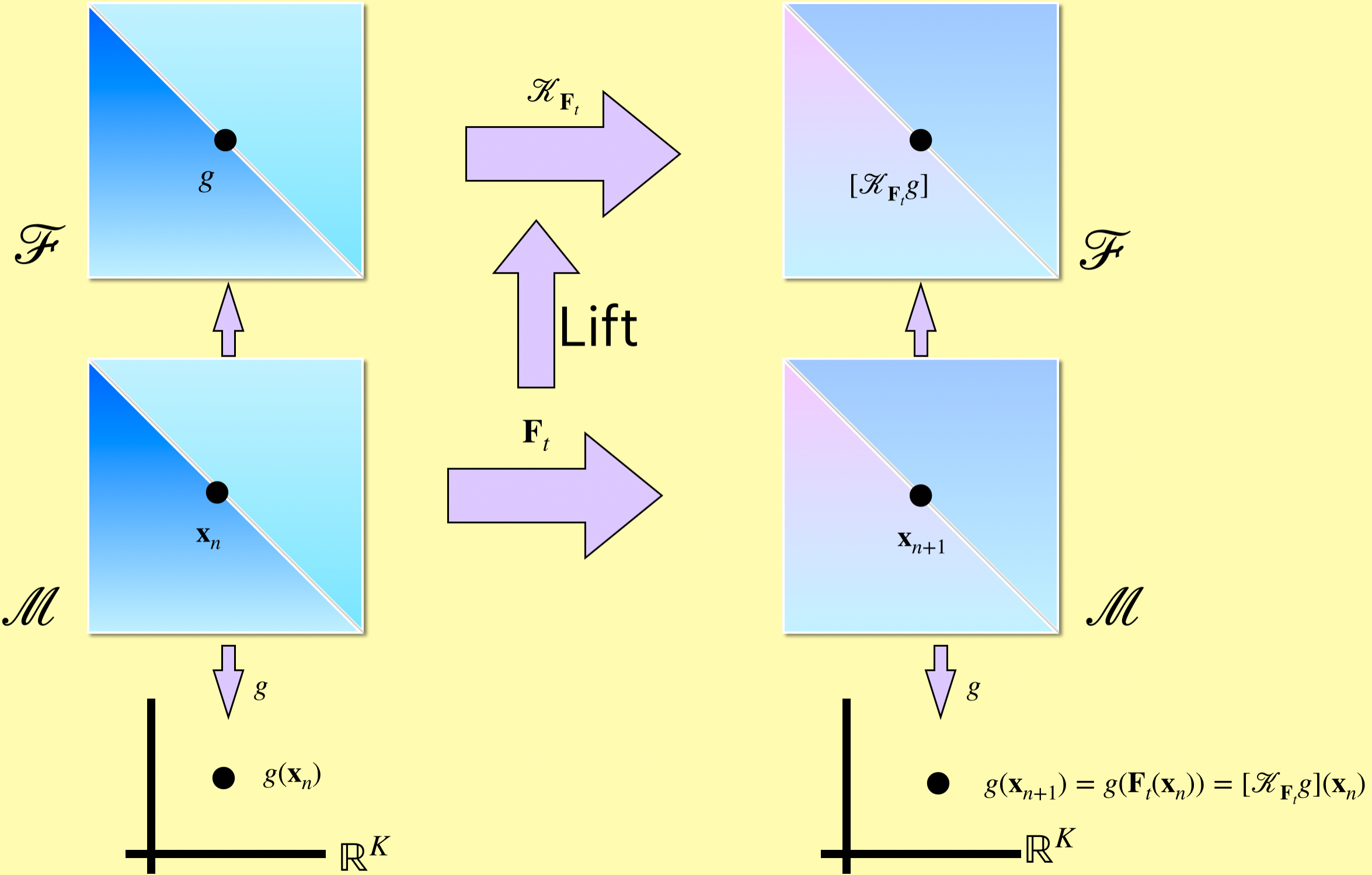}
    \caption{The lifting perspective of the Koopman operators.}
    \label{fig:LiftKoopman}
\end{figure}

The Koopman operators are linear despite being maybe defined over nonlinear dynamical system. However, there is an unavoidable trade-off that is the Koopman operators are \emph{infinite} dimensional operators that lift the finite dimensional nonlinear dynamical system into infinite dimensional linear dynamical system on the observable function space, follow \cref{fig:LiftKoopman}.
\subsection{Koopman Mode Decomposition}
Recall the Borel measurable set against the Borel measure $\mu$ in the observable space $\mathcal{F}=L^2\left(\Omega,\mu\right)=\LOmega$. These systems, in fact, are general and common in practice due to the phenomena \emph{energy preserving} and \emph{volume preserving} with respect to $\mu$ over $\Omega$ \cite{arnol2013mathematical,dubrovin2012modern,halmos1977bernoulli,hill1986introduction,walters2000introduction}. Furthermore, we can also assume that the dynamical system considered is invertible in the complement of the sets in $\Omega$ where the measure $\mu$ of that set is $0$ \cite{eisner2015operator}. With this assumption, we can have the Koopman operator as \emph{unitary} $\mathcal{K}:\LOmega\to\LOmega$ against the Borel measure $\mu$. In essence, if the adjoint of $\mathcal{K}$ is denoted by $\mathcal{K}^*$, and $I$ as an identity operator, then $\mathcal{KK}^*=\mathcal{K}^*\mathcal{K}=I$ over the Hilbert space $\LOmega$.
\subsubsection{Koopman eigenvalues and eigenfunctions}
Let an observable $\zeta$ be chosen from $\LOmega$ and if $\zeta$ is an eigenfunction of $\koop:\LOmega\to\LOmega$ with its eigenvalue $\lambda$, then $\forall N\in\mathbb{W}$ following is true:
\begin{align}
    \zeta\left(\mathbf{x}_N\right)=\left[\mathcal{K}^N\zeta\right](\mathbf{x}_0)=\lambda^N\zeta(\mathbf{x}_0).
\end{align}
The eigenvalue $\lambda\in\mathbb{C}$ governs the oscillatory and growth/decay behavior of the observable $\zeta\in\LOmega$. In particular, the oscillation behavior of $\zeta$ is determined by the complex argument of $\lambda$ and, on the other hand, the growth or decay of $\zeta$ is dictated by the absolute value of $\lambda$. Therefore, the eigen-observables pair $(\zeta,\lambda)$ potentially encodes the key insights of the dynamical system \cite{mezic2004comparison,mezic2005spectral,mezic2020spectrum,mezic2021koopman}. With this in mind, we are now ready to provide a brief description of Koopman mode decomposition.

Let the Koopman eigenfunction $\zeta$ be indexed by $\left\{\zeta_n\right\}_{n\in\mathbb{N}}$ and similarly for the eigenvalues as well, that is, $\left\{\lambda_n\right\}_{n\in\mathbb{N}}$ and finally $\mathcal{K}\zeta_n=\lambda_n\zeta_n$. If $\mathcal{K}$ has an invariant subspace spanned out by $\left\{\zeta_n\right\}_{n\in\mathbb{N}}$, then the \emph{Koopman mode decomposition} expands the whole tray of observables $\mathbf{g}\coloneqq\left[g_1,g_2,\ldots,g_\ell\right]^\top$ with respect to the \emph{Koopman modes} $\mathbf{c}_n$ in this subspace \cite{mezic2005spectral,rowley2009spectral} as follows:
\begin{align}\label{eq_KoopmanModes}
    \mathbf{g}\left(\mathbf{x}\right)=\sum_{n\in\mathbb{N}}\mathbf{c}_n\zeta_n(\mathbf{x})\quad\text{implying}\quad\mathbf{g}(\mathbf{x}_N)=\sum_{n\in\mathbb{N}}\mathbf{c}_n\lambda_n^N\zeta_n(\mathbf{x}_0).
\end{align}
The Koopman modes $\mathbf{c}_n$ depend on the choice of the observables $\zeta_n$ \cite{rowley2009spectral}. The following set of definitions plays an important role in the paper that sets the key stage shared between the data-driven Koopman operator and that of the operator-theoretic Koopman operator. 
\subsection{\texorpdfstring{$0$}{}\textsuperscript{th} Snapshot \texorpdfstring{$m$}{}\textsuperscript{th} Koopman mode difference}\label{def_0msnapdiff}
We define an important notion of the Koopman mode decomposition difference between the observable space $\LOmega$ and an arbitrarily Hilbert space $\mathcal{H}$ when considering such an observable that lives in both Hilbert spaces.
\begin{definition}\label{def_0thmthdifference}
    Let $\mathcal{K}_{\mathbf{F}_t}:\LOmega\to\LOmega$ be the Koopman operator corresponding to the dynamical system \eqref{eq_1}. Let $\varphi:\Cn\to\Cn$ be a holomorphic function defined over a Hilbert space $\mathcal{H}$ and consider the Koopman operator induced by $\varphi$ as $\mathcal{K}_{\varphi}:\mathcal{H}\to\mathcal{H}$. 
    
    For the Koopman eigenfunction observable $\zeta_n\in\LOmega\cap\mathcal{H}$ along the $0$\textsuperscript{th} snapshot $\mathbf{x}_0$, the `$0$\textsuperscript{th} snapshot $m$\textsuperscript{th} Koopman mode decomposition difference' with respect to the Hilbert space $\mathcal{H}$ denoted by $\Delta_{\mathcal{H}}^{(0,m)}$, given as
    \begin{align}\label{eq_0thmthdiff}
        \Delta_{\mathcal{H}}^{(0,m)}\coloneqq\sum_{n\in\mathbb{N}}\mathbf{c}_n\left[\mathcal{K}_{\mathbf{F}_t}^m\zeta_n-\mathcal{K}_{\varphi}^m\zeta_n\right](\mathbf{x}_0).
    \end{align}
    We can associate with the partial sum of $\Delta_{\mathcal{H}}^{(0,m)}$, denoted by $\mathcal{P}_N\Delta_{\mathcal{H}}^{(0,m)}$, given as
    \begin{align}\label{eq_0thmthdiff_partial}
        \mathcal{P}_N\Delta_{\mathcal{H}}^{(0,m)}\coloneqq\sum_{n=1}^N\mathbf{c}_n\left[\mathcal{K}_{\mathbf{F}_t}^m\zeta_n-\mathcal{K}_{\varphi}^m\zeta_n\right](\mathbf{x}_0).
    \end{align}
    Lastly, $\Delta_{\mathcal{H}}^{(0,m)}=\lim_{N\to0}\mathcal{P}_N\Delta_{\mathcal{H}}^{(0,m)}$.
\end{definition}
The specialized sum introduced above can be regarded as the \emph{weighted deviation} of the spatial mode reconstruction through the Koopman mode decomposition with respect to when the observable lies in different Hilbert spaces; here, the \emph{weighted} word signifies the presence of the terms of the Koopman modes $\mathbf{c}_n$. This definition can be made richer and robust, which can help justify the operator-theoretic investigation of the Koopman operators to its data-driven aspect counterpart. We do the same in the following part of the definition by making this difference \emph{faithful}; this will be essential in drawing the conclusion for the spatial temporal reconstruction over the RKHS of the Laplacian kernel.
\begin{definition}\label{def_0thmthdifference_faithful}
    Under the initial configuration of \Cref{def_0thmthdifference}, we say that $\Delta_{\mathcal{H}}^{(0,m)}$ with respect to the Hilbert space $\mathcal{H}$ is faithful if there exists an invertible $\mathcal{A}\in\mathbb{C}^{D\times D}$ whose Frobenius norm satisfies $0<\|\mathcal{A}\|_{\text{Frob}}<1$ such that $\varphi(\bmz)\coloneqq\mathcal{A}\bmz+\bm{b}$ when closable $\mathcal{K}_{\varphi}:\mathcal{H}\to\mathcal{H}$ is compact over $\mathcal{H}$. The faithful $\Delta_{\mathcal{H}}^{(0,m)}$ with respect to $\mathcal{H}$ is denoted by $\digamma\Delta_{\mathcal{H}}^{(0,m)}$, is explicitly given as:
    \begin{align}\label{eq_faithful0thmth}
        \digamma\Delta_{\mathcal{H}}^{(0,m)}\coloneqq\sum_{n\in\mathbb{N}}\mathbf{c}_n\left[\mathcal{K}_{\mathbf{F}_t}^m\zeta_n(\mathbf{x}_0)-\zeta_n\left(\mathcal{A}^m\mathbf{x}_0+\bm{b}\sum_{i=0}^{m-1}\mathcal{A}^i\right)\right].
    \end{align}
We can associate with the partial sum of $\digamma\Delta_{\mathcal{H}}^{(0,m)}$, denoted by $\mathcal{P}_N\digamma\Delta_{\mathcal{H}}^{(0,m)}$ for some $N\in\mathbb{N}$, given as 
\begin{align}\label{eq_faithful0thmth_faithful}
    \mathcal{P}_N\digamma\Delta_{\mathcal{H}}^{(0,m)}\coloneqq\sum_{n=1}^N\mathbf{c}_n\left[\mathcal{K}_{\mathbf{F}_t}^m\zeta_n(\mathbf{x}_0)-\zeta_n\left(\mathcal{A}^m\mathbf{x}_0+\bm{b}\sum_{i=0}^{m-1}\mathcal{A}^i\right)\right].
\end{align}
Lastly, $\digamma\Delta_{\mathcal{H}}^{(0,m)}=\lim_{N\to\infty}\mathcal{P}_N\digamma\Delta_{\mathcal{H}}^{(0,m)}$. 
When such a specific operator-theoretic quantification of the Koopman operator is determined subject to its existence over the Hilbert space $\mathcal{H}$, then
\begin{align}\label{eq_faithful0thmth_faithful0}
    \digamma\Delta_{\mathcal{H}}^{(0,m)}=\lim_{N\to\infty}\mathcal{P}_N\digamma\Delta_{\mathcal{H}}^{(0,m)}\equiv0
\end{align}
\end{definition}
\Cref{def_0thmthdifference_faithful} sheds the light on the existence of those closable Koopman operators which admit compactification on the Hilbert space $\mathcal{H}$ that can justify the existence of both \eqref{eq_faithful0thmth} and \eqref{eq_faithful0thmth_faithful}. 
Obviously, this depends on the choice of the Hilbert space $\mathcal{H}$; the conclusion of this article provides an application of this definition.

\subsection{Spectrum, spectral measure and diagonalization}\label{subsection_KoopmanSpectralMeasure}
For the Koopman operator $\mathcal{K}:\LOmega\to\LOmega$, we denote its \emph{spectrum} by $\rho(\mathcal{K})$ and it is defined as follows:
\begin{align}
    \rho(\mathcal{K})\coloneqq\left\{z\in\mathbb{C}:\frac{1}{\mathcal{K}-zI}\notin\text{bounded operator}\right\}\subset\mathbb{C}.
\end{align}
Since the Koopman operator $\mathcal{K}:\LOmega\to\LOmega$ is a unitary operator, therefore $\rho(\mathcal{K})\subset\mathbb{D}$, where $\mathbb{D}=\left\{z\in\mathbb{C}:|z|=1\right\}$.
The spectral theorem \cite{hall2013quantum} provides the \emph{diagonalization} of the Koopman operator $\mathcal{K}$ by employing the \emph{projection-valued measure} $\mathds{V}$ on $\rho(\mathcal{K})\subset\mathbb{D}$. For all Borel measurable set $\mathds{B}\subset\mathbb{D}$, the projection valued measure $\mathds{V}$ links an orthogonal projection $\mathds{V}(\mathds{B})$ which eventually yields following spectral decomposition over $g\in\LOmega$:
\begin{align}
    g=\left(\int_{\mathbb{D}}1d\mathds{V}(\lambda)\right)g\quad\text{with}\quad\mathcal{K}g=\left(\int_{\mathbb{D}}\lambda d\mathds{V}(\lambda)\right)g.
\end{align}
The nature of the projection value measure confirms $\mathds{V}(\mathbb{D})=1$. We can enrich our results by incorporating functional calculus as follows
\begin{align}
    g\left(\mathbf{x}_N\right)=\left[\mathcal{K}^Ng\right](\mathbf{x}_0)=\left[\left(\int_{\mathbb{D}}\lambda^N d\mathds{V}(\lambda)\right)g\right](\mathbf{x}_0),
\end{align}
where $\mathbf{x}_N$ is the $N$\textsuperscript{th} data snapshot delivered by the system \eqref{eq_1} given under \Cref{assumption_1}.
\subsection{Kernel Space}
\begin{definition}[Kernel Function]
    Let $X=\emptyset$, then a function $k:X\times X\to\mathbf{K}$ is called the kernel on $X$ if there exists a $\mathbf{K}-$Hilbert space $\left(H,\langle\cdot,\cdot\rangle_H\right)$ accompanied by a map $\Phi:X\to H$ such that $\forall x,x'\in X$, we have
    \begin{align}
        k(x,x')=\langle\Phi(x'),\Phi(x)\rangle_H.\label{eq_kernel_H}
    \end{align}
    We consider $\Phi$ the feature map and $H$ the feature space of $k$.
\end{definition}
We provide the definition of radial functions (cf. \cite{fasshauer2007meshfree}) and then some examples of that remain in direct connection with the present article.
\begin{definition}
    A function $\Phi:\mathbb{R}^D\to\mathbb{R}$ is called as \emph{radial} if there exists a \emph{uni-variate} function $\phi:\{0\}\cup\mathbb{R}_+\to\mathbb{R}$ such that 
        \begin{align*}
            \Phi\left(\bm{x}\right)=\phi(r),~\text{where $r=\|\bm{x}\|$}.
        \end{align*}
\end{definition}
\begin{example}
    The most common kernel function that arises from the class of radial basis function 
    used in modern machine learning and artificial intelligence routine such as speech enhancement is the class of \emph{exponential power kernels} \cite{giraud2005positive,hui2018kernel} given as 
\begin{align}\label{eq_expkernel}
    K_{\text{exp}}^{\gamma,\sigma}(\bm{x},\bm{z})\coloneqq\exp\left(-\frac{\|\bm{x}-\bm{z}\|_2^\gamma}{\sigma}\right);~\text{where $\gamma,\sigma>0$ and $\bm{x},\bm{z}\in\mathbb{C}^D$}.
\end{align}
The constant $\sigma$ in \eqref{eq_expkernel} is referred to as \emph{kernel bandwidth} and $\gamma$ is called \emph{shape parameter}. When $\gamma=1$ in \eqref{eq_expkernel}, we get $K_{\text{exp}}^{1,\sigma}(\bm{x},\bm{z})\coloneqq\exp\left(-\frac{\|\bm{x}-\bm{z}\|_2}{\sigma}\right)$, which is referred to as \emph{Laplace Kernel}. When $\gamma=2$ in \eqref{eq_expkernel}, we get $K_{\text{exp}}^{2,\sigma}(\bm{x},\bm{z})\coloneqq\exp\left(-\frac{\|\bm{x}-\bm{z}\|_2^2}{\sigma}\right)$, which is referred to as \emph{Gaussian Radial Basis Function (GRBF) Kernel}. We direct the interested readers \cite{rasmussen2006gaussian} to learn more about the feature map of the GRBF kernel. It should be noted that the choice of kernel functions can dramatically change the performance of the (supervised) ML routine \cite{geifman2020similarity,NEURIPS2020_1006ff12}, specifically in those situations where a shorter training time or limited information is available.
\end{example}
\begin{figure}[H]
    \centering
    \includegraphics[width=0.23\linewidth]{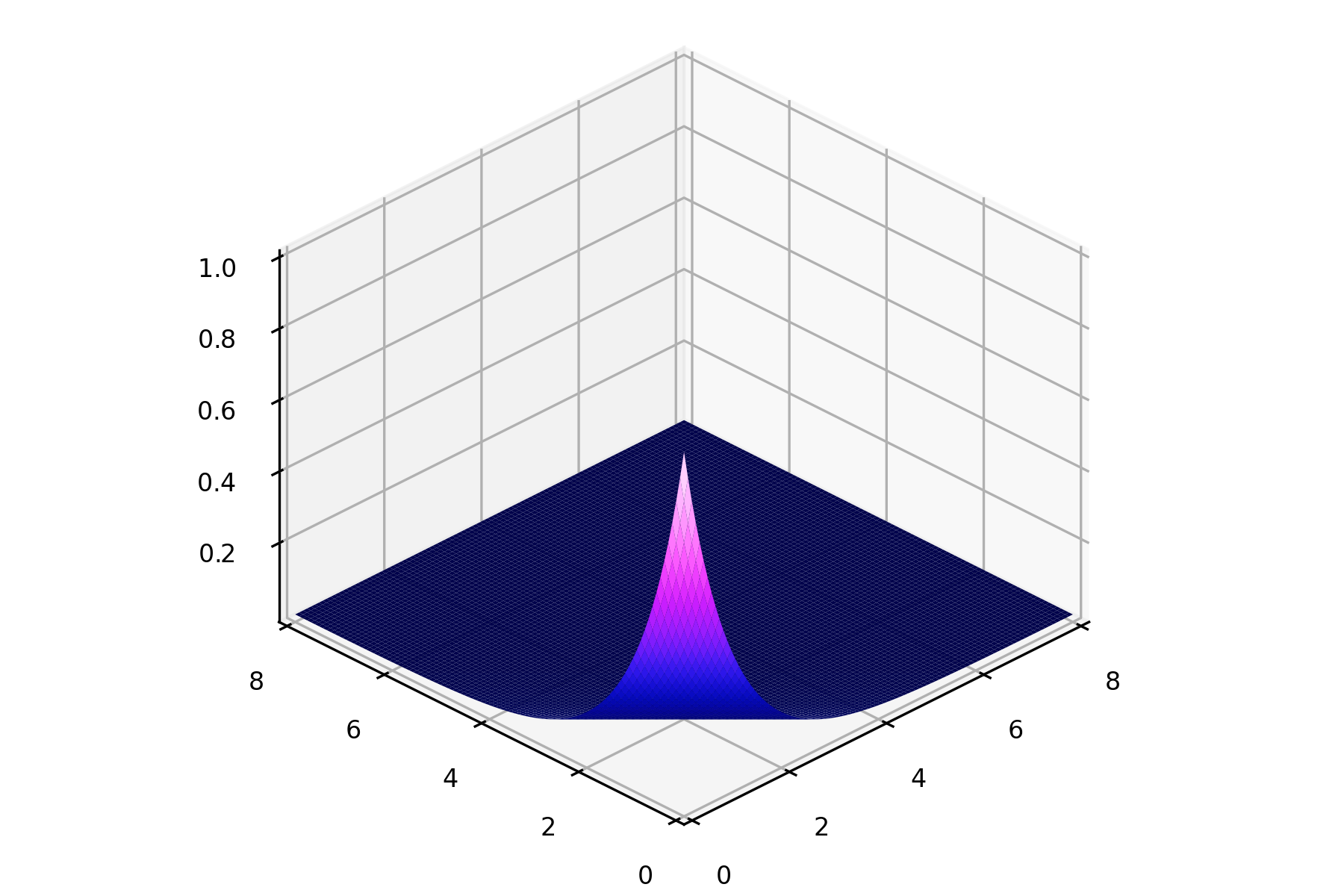}
    \includegraphics[width=0.1557\linewidth]{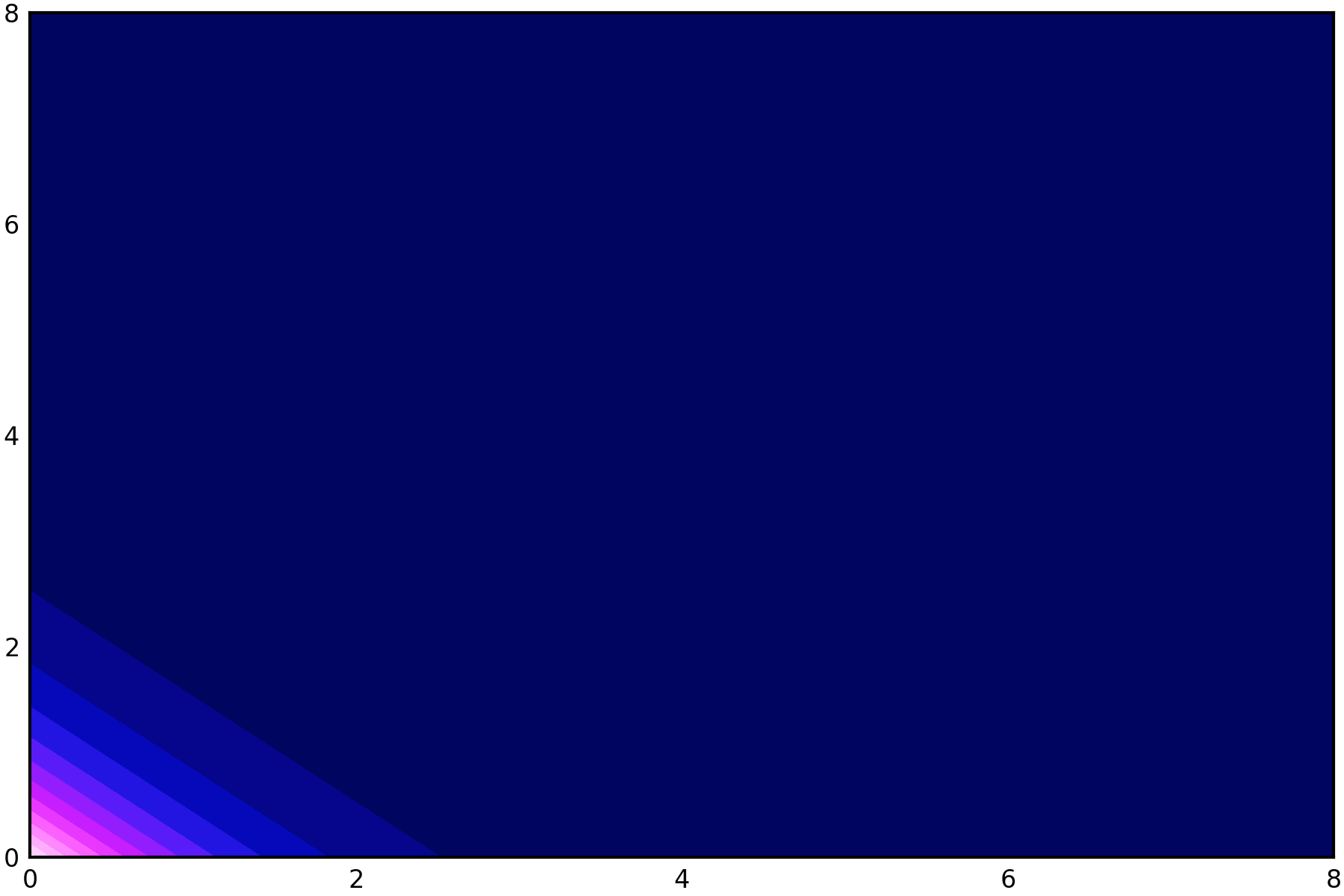}
    \includegraphics[width=0.23\linewidth]{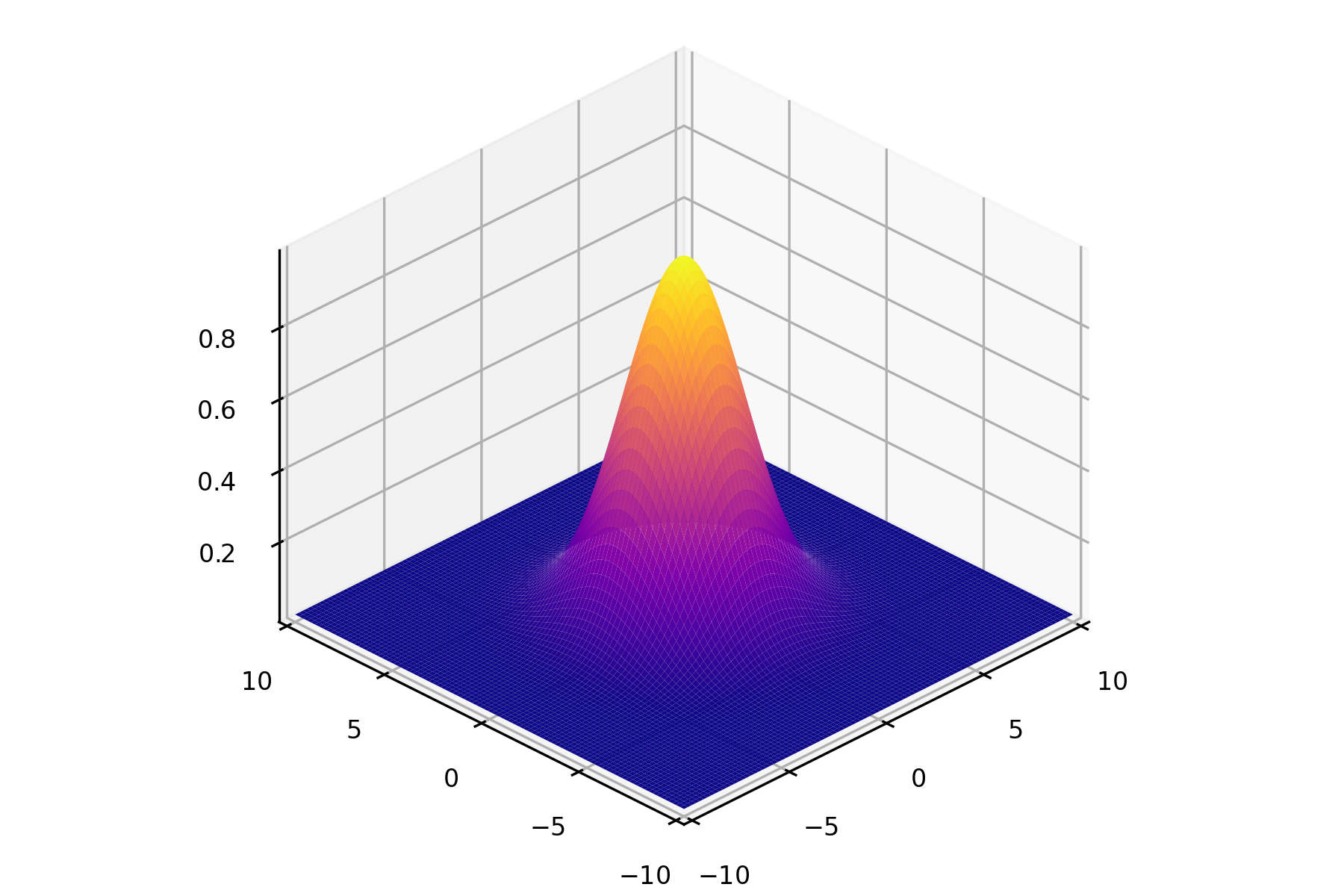}
    \includegraphics[width=0.1557\linewidth]{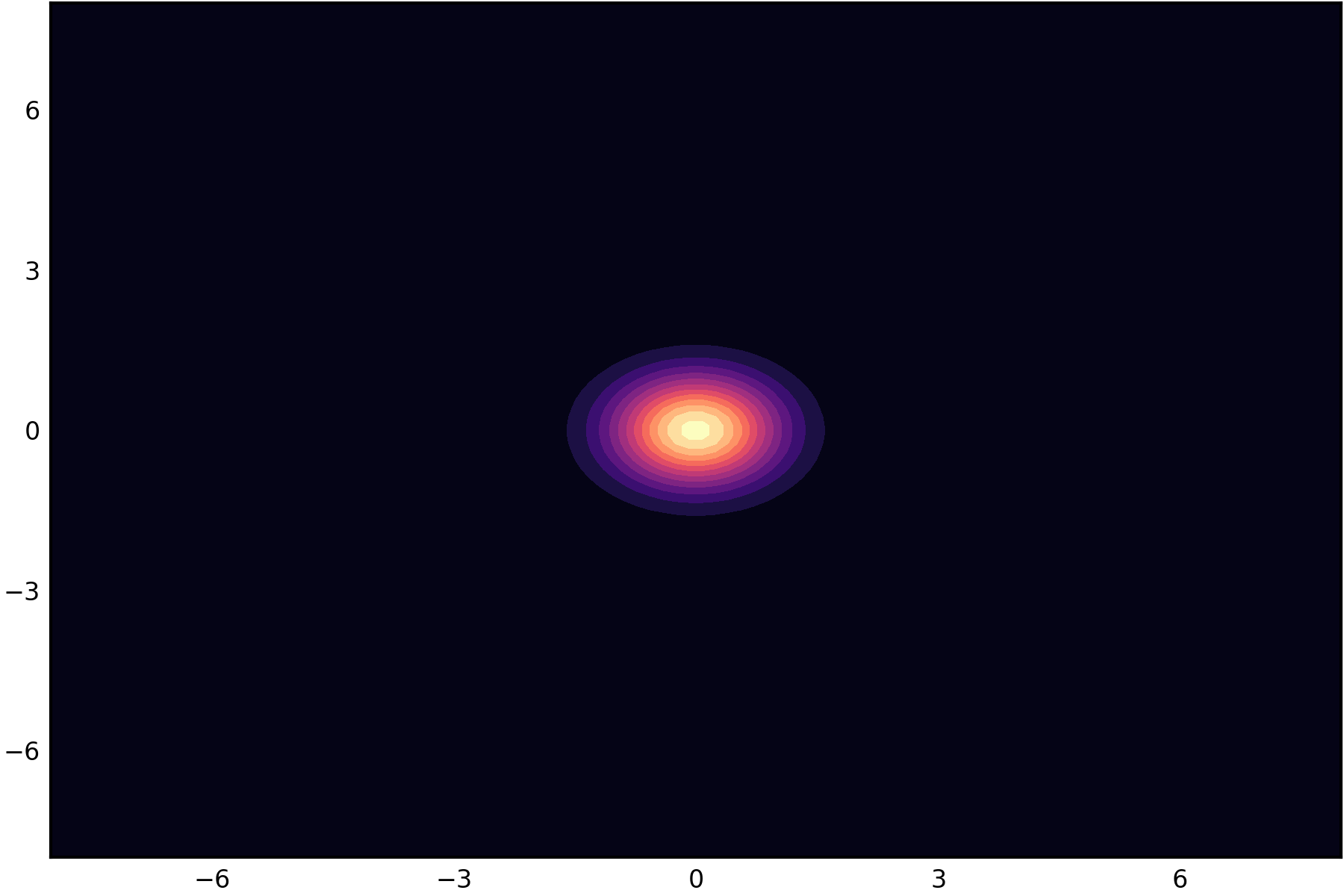}
    \caption{Pairwise respective surface and contours plot of kernel functions
    $K_{\text{exp}}^{1,\sigma}$ and $K_{\text{exp}}^{2,\sigma}$.}
    \label{fig:LAPGRBF}
\end{figure}
\begin{definition}[Reproducing Kernel Hilbert Space]\label{dfn_RKHS}
    Let $X=\emptyset$ and $\left(H,\langle\cdot,\cdot\rangle_H\right)$ be the Hilbert function space on $X$.
    \begin{enumerate}
        \item The space $H$ is called as the \emph{\textbf{reproducing kernel Hilbert space}} (RKHS) if $\forall x\in X$, the evaluation functional $\mathcal{E}_x:H\to\mathbf{K}$ defined as $\mathcal{E}_x(f)\coloneqq f(x),~f\in H$ is continuous.
        \item A function $k:X\times X\to\mathbf{K}$ is called \emph{reproducing kernel} of $H$ if we have:
        \begin{enumerate}
            \item $k(\cdot,x)\in H~\forall x\in X$, that is $\|k(\cdot,x)\|_H<\infty$, and
            \item $k(\cdot,\cdot)$ has the reproducing property; that is $f(x)=\langle f,k(\cdot,x)\rangle_H~\forall f\in H\text{~and~} x\in X.$
        \end{enumerate}
        \end{enumerate}
\end{definition}
Note that the norm convergence above yields the point-wise convergence inside RKHS. We also recall the Moore-Aronszajn Theorem which we will use to construct the RKHS based on the Laplacian kernel as a $L^2-$measure.
\begin{theorem}[Moore-Aronszajn Theorem \cite{aronszajn1950theory}]\label{theorem_aronsjan}
    Let $H$ be an RKHS on an empty set $X$, then $k:X\times X\to\mathbf{K}$ defined as $k(x,x')\coloneqq \langle\mathcal{E}_x,\mathcal{E}_{x'}\rangle_H$ for $x,x'\in X$ is the only reproducing kernel of $H$. Furthermore, for some index set $\mathcal{I}$, if we have $\left\{\mathbf{e}_i\right\}_{i\in\mathcal{I}}$ as an orthonormal basis, then for all $x,x'\in X$, we have the following with absolute convergence
    \begin{align}\label{eq_5theorem1.4}
        k(x,x')=\sum_{i\in\mathcal{I}}\mathbf{e}_i(x)\overline{\mathbf{e}_i(x')}.
\end{align} 
    \end{theorem}
\section{Kernel Extended Dynamic Mode Decomposition}\label{section_KeDMD}
The dynamic mode decomposition (DMD) algorithm can easily and accurately characterize both periodic and quasi-periodic behaviors in nonlinear dynamical systems. Unfortunately, it fails to portray the true nonlinear phenomena when it is relying on linear observables. To overcome this very limitation, another variant of DMD was introduced in \cite{williams2015data} which takes advantage of \emph{kernel trick} \cite{scholkopf2000kernel} which also depends on the implicit choice of kernel functions. 

\begin{figure}[H]
    \centering
    \includegraphics[width=0.32\linewidth]{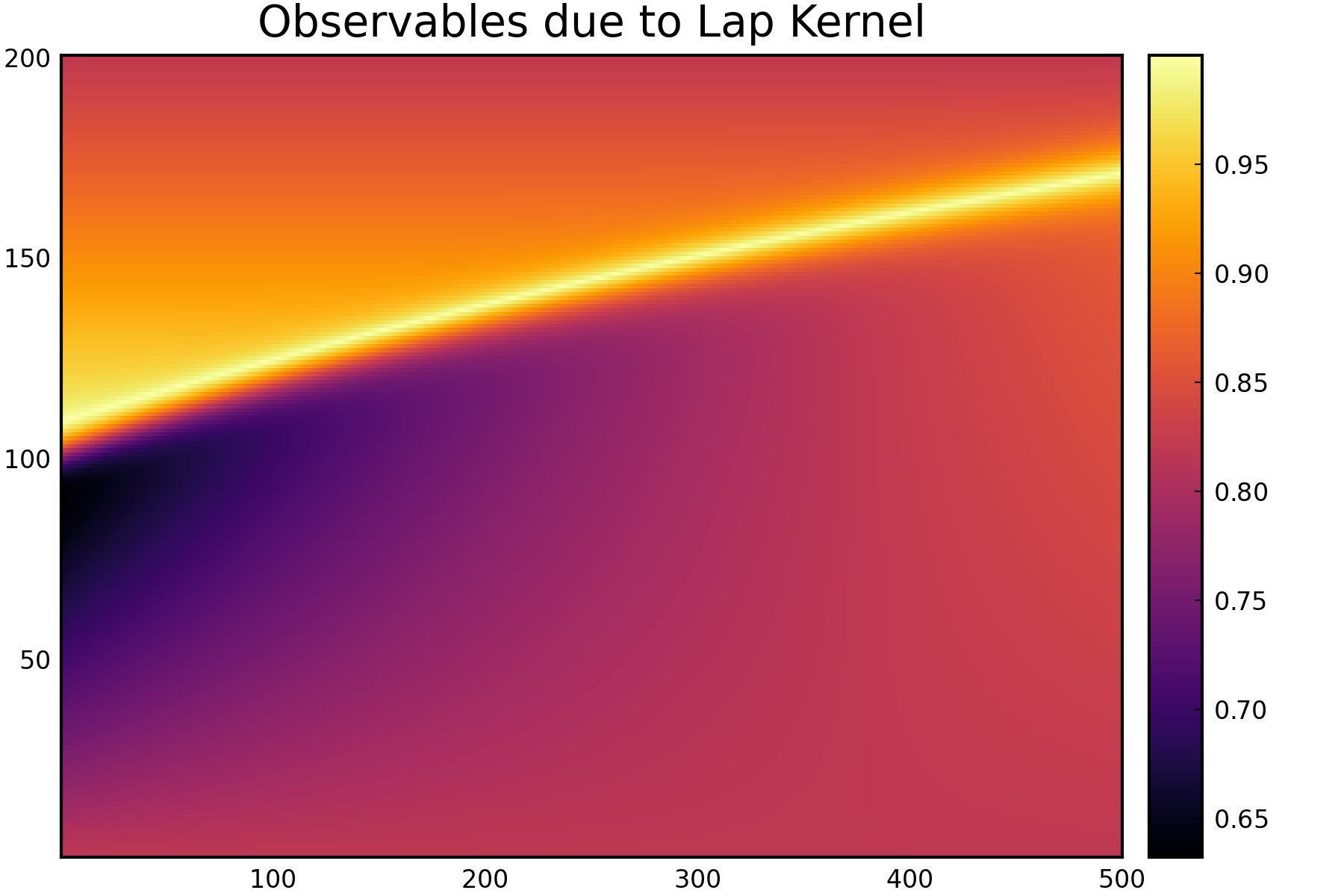}
    \includegraphics[width=0.32\linewidth]{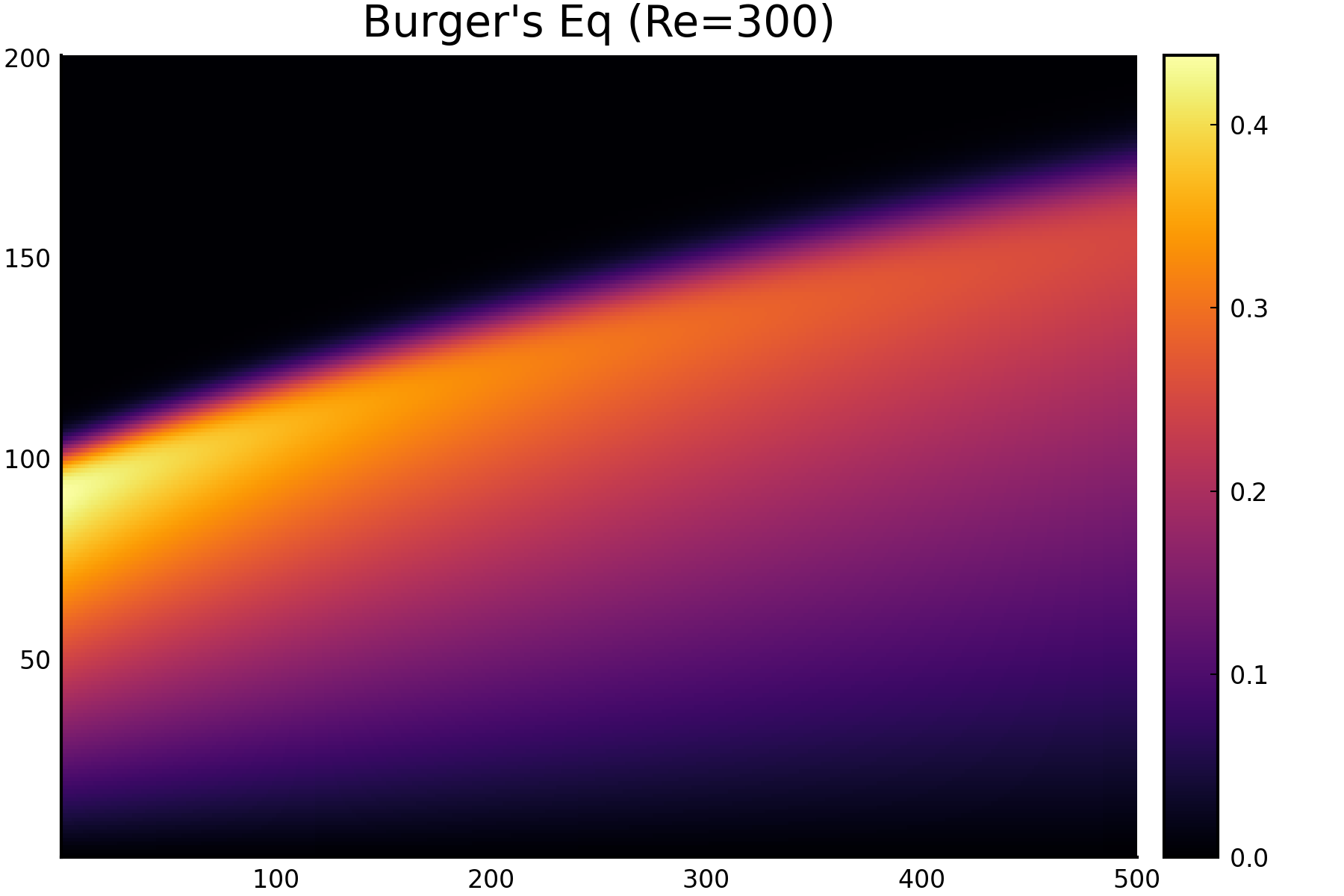}
    \includegraphics[width=0.32\linewidth]{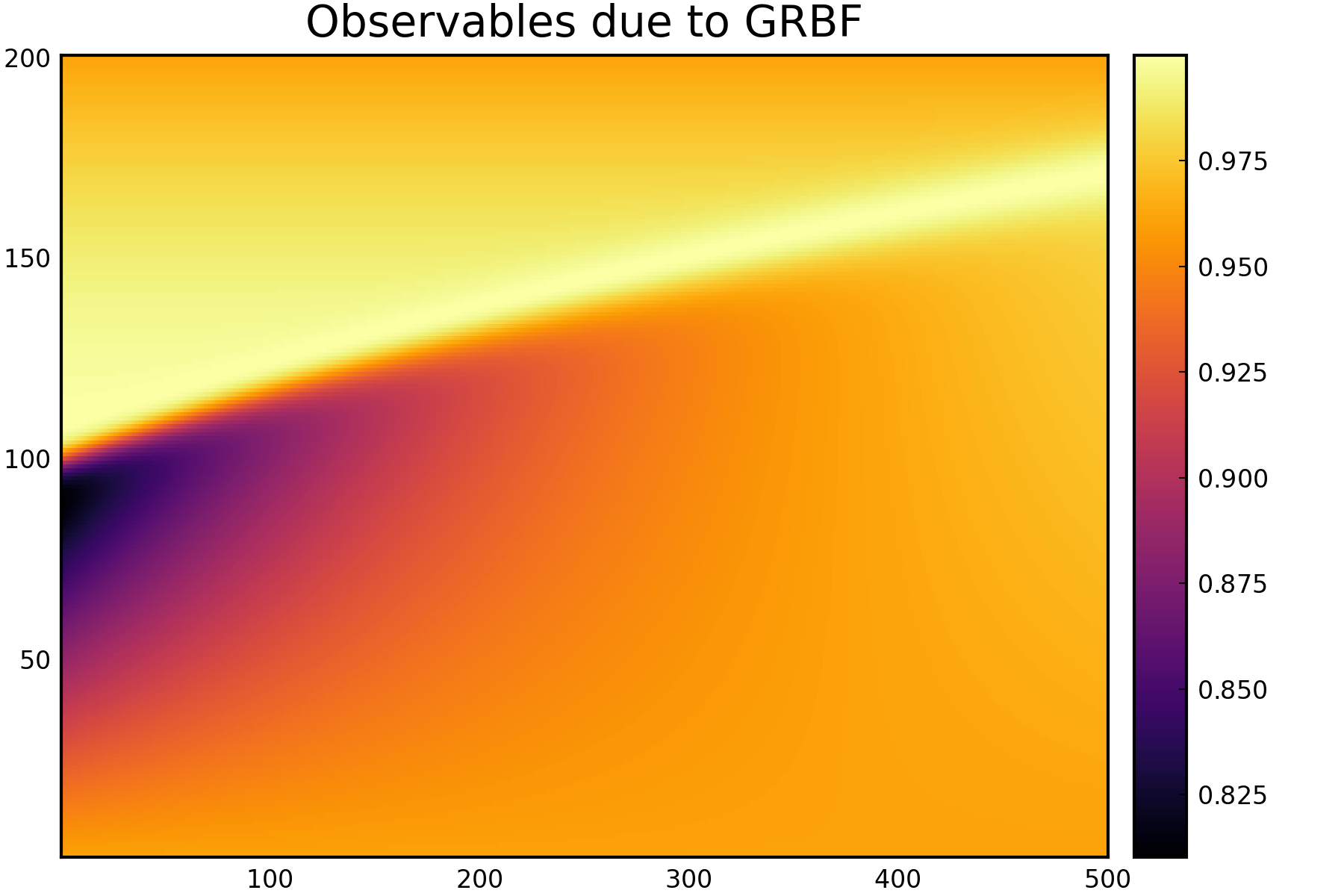}
    \caption{An illustration of notion of observables for the Burger Equation's data \emph{(middle)} when we chose the dictionary of Laplacian kernel \emph{(left)} and GRBF kernel \emph{(right)}.}
    \label{fig:Kernel_Observables}
\end{figure}

We briefly recall the kernel-extended DMD algorithm as follows. Suppose 
\begin{align*}
    \left[\vx_1,\vx_2,\ldots,\vx_M\right]\in\mathbb{R}^{N\times M}
\end{align*} be $M$ data points from a $N$-dimensional dynamical system corresponding to $\vx_i\to\mathbf{F}_t\left(\vx_i\right)$ which evolves on $\mathbb{R}^N$ with respect to time. Let 
\begin{align*}
    \bm{\psi}\left(\vx\right)=\left[\psi_1(\vx),\psi_2(\vx),\ldots,\psi_N(\vx)\right]:\mathbb{R}^N\to\mathbb{R}^K
\end{align*}
be the set of vector-valued Koopman observables such that $\forall i$, $\psi_i:\mathbb{R}^N\to\mathbb{R}$ is square-integrable. Let 
\begin{align}
    \mathbf{G}&=\frac{1}{M}\sum_{i=1}^{M-1}\bm{\Psi}\left(\vx_i\right)^\top\bm{\Psi}\left(\vx_i\right)\label{eq_Gram}\\
    \mathbf{A}&=\frac{1}{M}\sum_{i=1}^{M-1}\bm{\Psi}\left(\vx_i\right)^\top\bm{\Psi}\left(\vx_{i+1}\right)\label{eq_Interaction}.
\end{align}
Then, the finite rank representation of the Koopman operator $\mathcal{K}_{\mathbf{F}_t}\approx\mathcal{K}\in\mathbb{R}^{K\times K}$ is obtained as the solution of the optimization problem 
\begin{align}\label{KoopmanOptimization}
    \min_{\mathcal{K}}\|\mathbf{G}\mathcal{K}-\mathbf{A}\|\overset{\text{\cite{williams2015data}}}{=}{\mathcal{K}\triangleq\mathbf{G}^{+}\mathbf{A}.}
\end{align}
\subsection{Results}
Now, that we have recalled the Kernel Extended DMD, in this section, we provide the results for the spatial-temporal reconstruction results both by the Laplacian kernel and the GRBF kernel. 
\begin{details}
We now provide the governing equations for the data set considered, and these are given below in \Cref{tbl:table_of_governing_equations}.
\begin{table}[H]
        \centering
        \caption{Governing equations of experiments based on its availability (refer \Cref{tab:dataset}, if needed).}
        \begin{tabular}{cM{110mm}
        }
           \toprule
            \small{\textsc{\textbf{\#}}} 
            & \small{\textbf{\textsc{Governing Equation(s)}}} 
            \\
            \midrule
            1. 
            & $\dfrac{\partial u}{\partial t}+u\dfrac{\partial u}{\partial x}=\nu\dfrac{\partial^2 u}{\partial x^2}$ 
            \\
            2. 
            & $\dfrac{\partial}{\partial t}\mathbf{u}(x,y,t)+\mathbf{u}(x,y,t)\cdot\nabla\mathbf{u}(x,y,t)+\nabla p(x,y,t)-\dfrac{1}{\textsc{Re}}\nabla^2\mathbf{u}(x,y,t)=0.$  
            \\
            3. 
            & $
    \Dot{\vx}=\vy,\hspace{.2em}
    \Dot{\vy}=-0.5\vy+\vx-\vx^3.$ 
             \\
            4. 
            & \textcolor{red}{{\tt{n/a}}} 
             \\
             5. & $\Dot{x}=10(y-x),~\Dot{y}=x(28-z)-y,~\Dot{z}=xy-8/3z$\\
             6. &$\Dot{x}=-y-z,~\Dot{y}=x+0.2y,~\Dot{z}=0.2+z(x-5.7)$\\
             7. & \textcolor{red}{{\tt{n/a}}} 
             \\
            \bottomrule
        \end{tabular}
        \label{tbl:table_of_governing_equations}
    \end{table}
    \end{details}
\begin{datasetinfo}
The following table provides details of the experiments in terms of data matrix.
\begin{table}[H]
\caption{Details of experiments for recovering spatial-temporal modes.}
\begin{center}
\begin{tabular}{llll}
\multicolumn{1}{c}{\bf \small{\textsc{Experiment}}}  &\multicolumn{1}{c}{\bf \small{\textsc{True Data}}}
&\multicolumn{1}{c}{\bf \small{\textsc{Partial Rank}}}  &\multicolumn{1}{c}{\bf \small{\textsc{Reshape Data}}}
\\ \midrule \\
Burger's Equation  &$256\times101$ &$256\times\bm{40}$  & $-$ 
\\
Fluid flow across cylinder  &$89,351\times151$  & $89,351\times\bm{100}$ & $-$
\\
Duffing's Oscillator  &$2\times50,000$ & $2\times\bm{35,000}$ & $-$
\\
Seattle I5 highway traffic data &$72\times75$ & $72\times\bm{75}$ & $-$
\\
Lorenz 1963 Attractor & $3\times200,000$ &$3\times\bm{200,000}$& $\simeq20,000\times30$
\\
Rössler Attractor & $3\times64,000$&$3\times\bm{60,000}$ & $\simeq10,000\times18$
\\
NOAA Sea Surface Temp. Anomaly & $2,592\times2,101$ & $2,592\times\bm{10}$& $-$
\end{tabular}
\label{table:experiment_with_data_matrix_references}
\end{center}
\end{table}
\begin{note}\label{note_for_Exp4}
We have provided an extra column of information about the data matrix adjustment on which the Kernel Extended DMD is applied for rich visualizations of experimental results. In addition, this adjustment is necessary to make the result more achievable and comprehensible, in contrast to having unregulated and degenerate results. In the following figure, we show one such datamatrix to contrast the irregular snapshot with the ground truth.
\end{note}
\end{datasetinfo}
\begin{datamatrix}The considered data set are provided as follows:
\begin{figure}[H]
\centering
\setkeys{Gin}{width=\linewidth}
\begin{subfigure}{0.3\textwidth}
\includegraphics{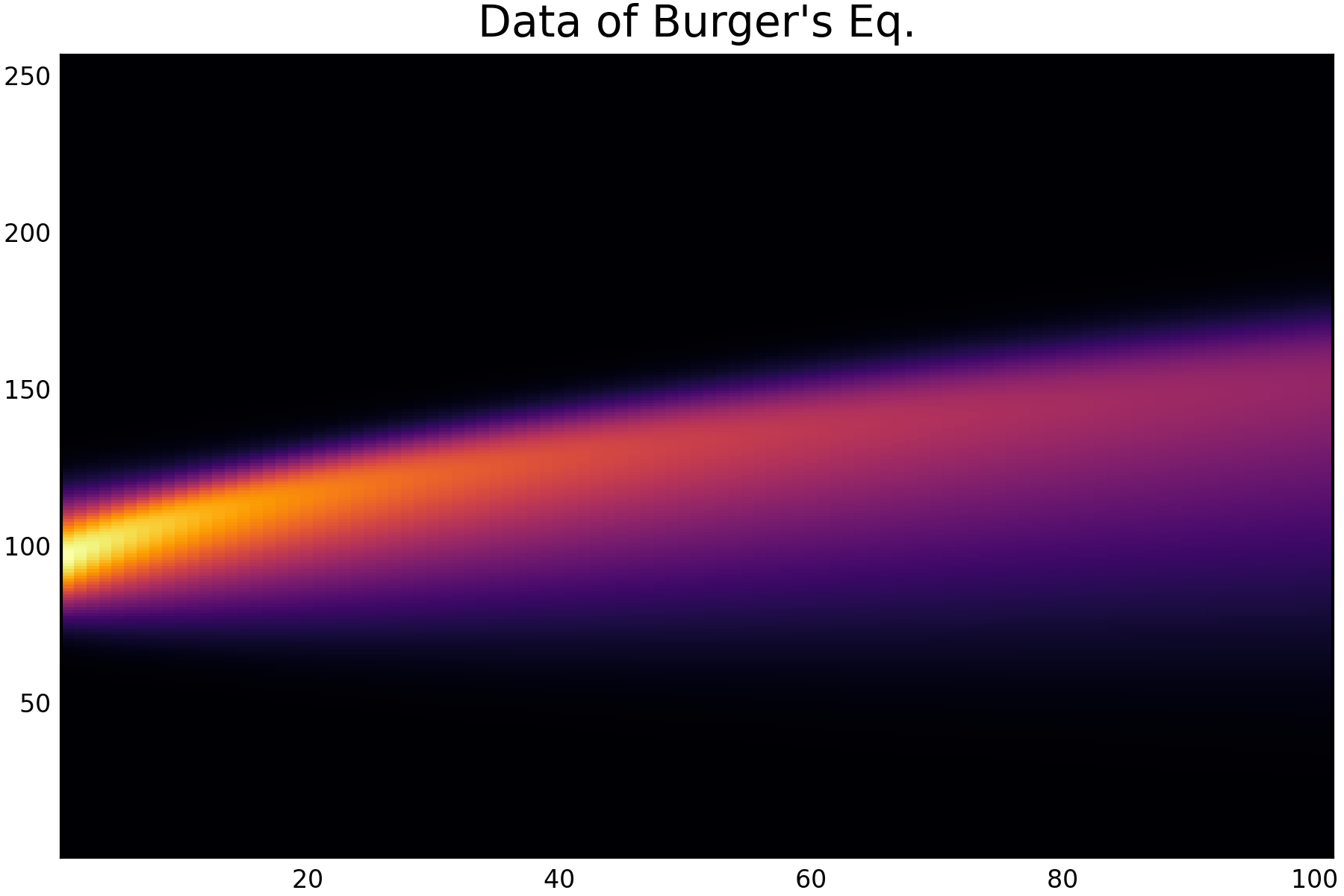}
\caption{Ground Truth}
\label{fig:Exp1_GroundTruth}
\end{subfigure}%
\hfil
\begin{subfigure}{0.3\textwidth}
\includegraphics{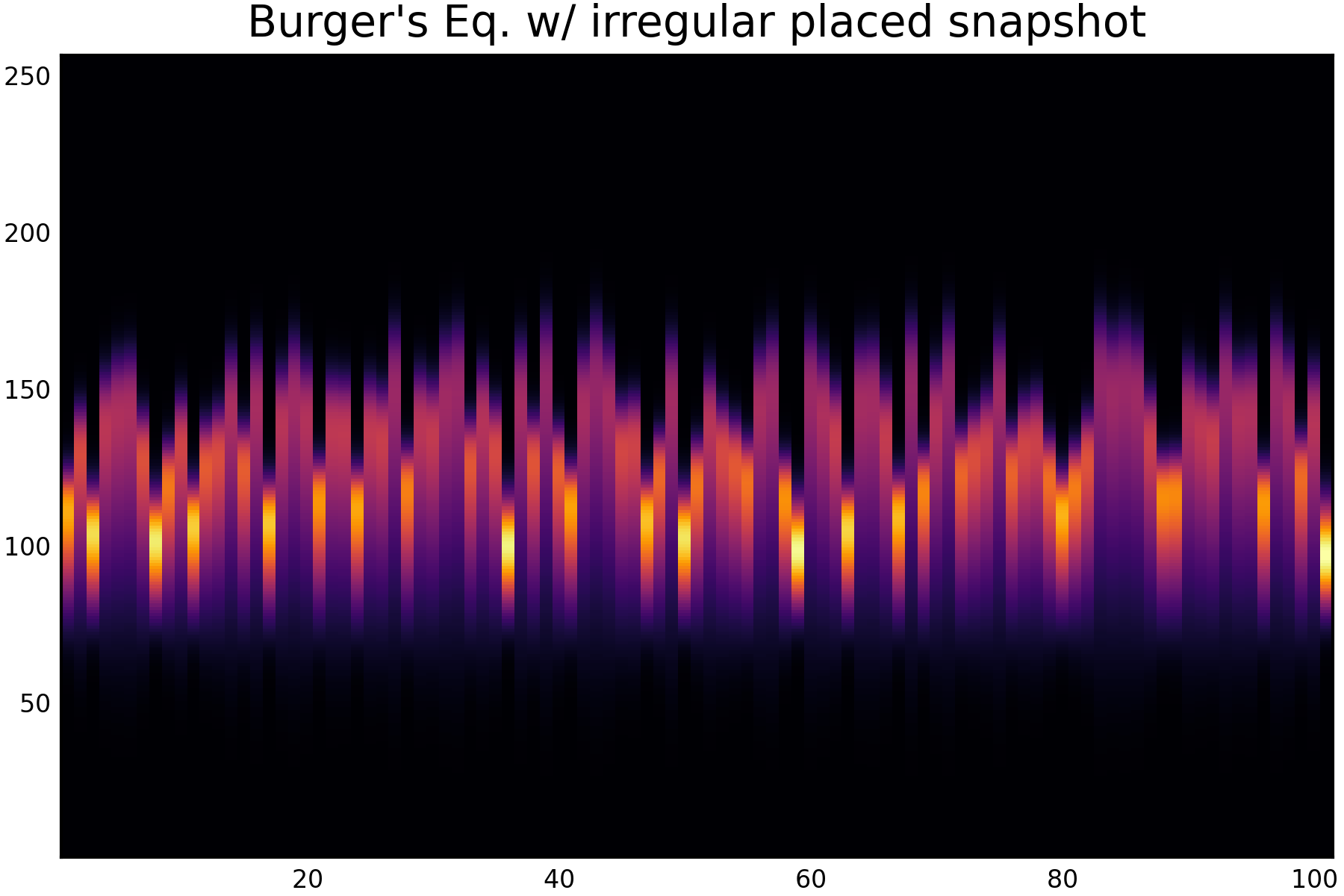}
\caption{Irregular \& Sparse}
\label{fig:Exp1_IrrSparse}
\end{subfigure}
\hfil
\begin{subfigure}{0.3\textwidth}
\includegraphics{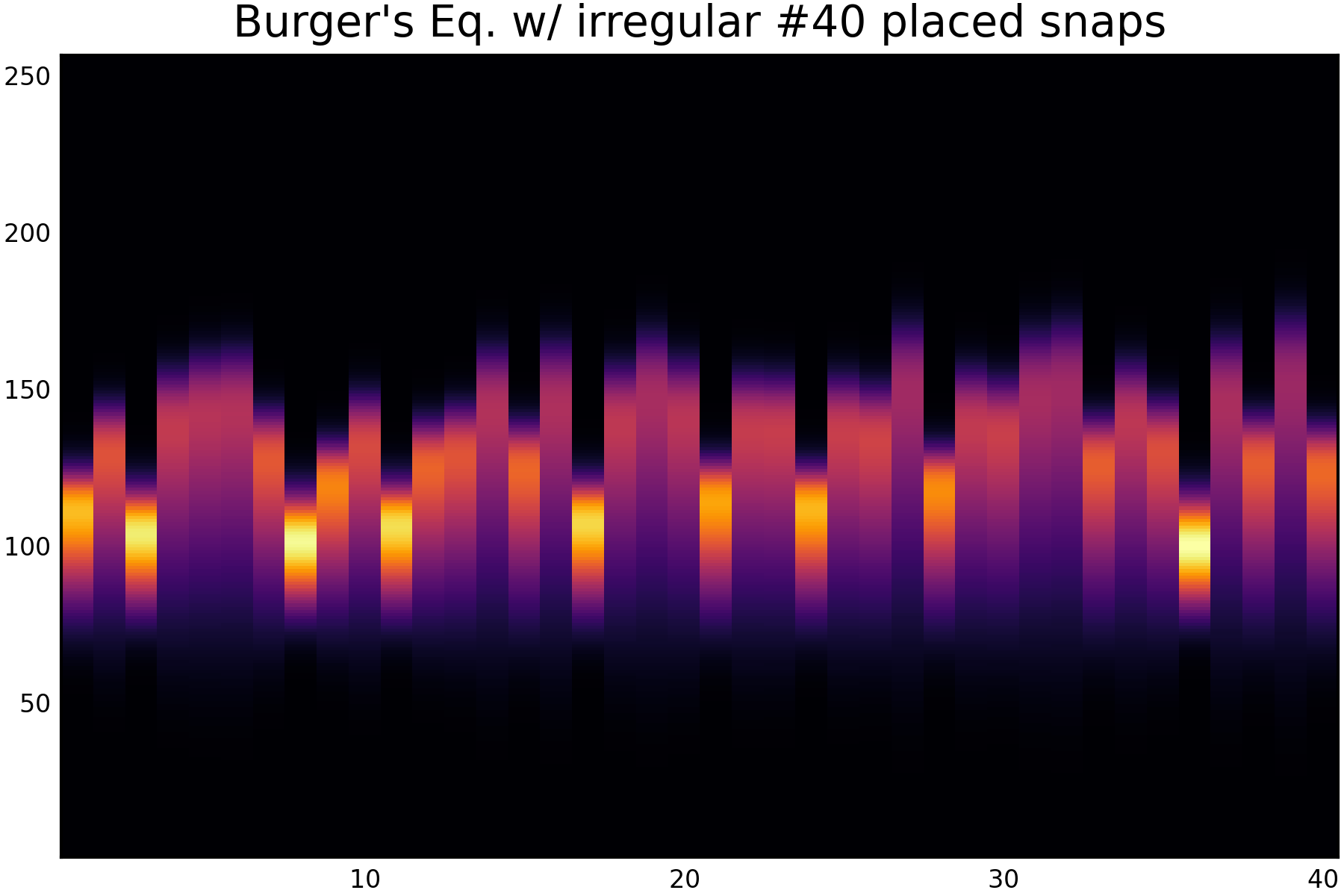}
\caption{Partial rank}
\label{fig:Exp1_PartailRank}
\end{subfigure}
\caption{\emph{\textbf{Experiment 1}}: Data matrix of non-linear Burger's Equation.}
\label{fig:DataMatrixExperiment1}
\end{figure}

\end{datamatrix}

\subsection{Experimental evaluations}
For all seven experiments under the assumption of \Cref{aim}, we provide the following details:
\begin{enumerate}
    \item For the experimental simulation, the irregular and sparse format of the data set of the first three experiments is delivered by operating {\tt{shuffle}} command available in {\tt{Julia}} on the Apple MacBook M2 MAX 32 GB machine.
    \item We provide the snapshot information in \emph{four formats} for Experiments 1, 2 \& 3 and \emph{three formats} for Experiment 4. These are categorized as \emph{Actual, Irregular \& Sparse, Reconstruction by $\Lap$ and Reconstruction by $\GRBF$}. Follow \Cref{note_for_Exp4} if needed!
    \item We compare the results reconstructed by $\Lap$ or $\GRBF$ with the snapshot of \emph{Irregular \& Sparse} in the form of a heat-map, which is an industry standard practice. We say reconstructed result \emph{poor} if there is a presence of (but not limited to) some unwanted or undesirable fringe pattern in the reconstructed result; we do not want poor results.
    \item We also compare the results reconstructed by $\Lap$ or $\GRBF$ with the snapshot of \emph{Irregular \& Sparse} in the form of spatial values and then also provide the error plot. We adopt the following mathematical definition for the consideration of \textbf{element-wise error (EWE)} between the \emph{reconstructed} $\mathbf{F}$ with respect to the \emph{actual} $\mathbf{A}$ for a total of $N-$data points:
    \begin{align*}
        \text{Element Wise Error}\coloneqq N^{-1}\left|1-\frac{\mathbf{F}}{\mathbf{A}}\right|.\tag{EWE}
    \end{align*}
    \item The data snapshot information is made available in the opening of the respective subsection for all four experiments. The snapshots considered are those for which the real part of the Koopman operator eigenvalue $(\in\mathbb{C})$ based on the provided irregular and sparse data set is close to $1$.
    \item After showcasing the results of all four experiments in every subsection of this section, we close this section providing detailed insights in the last part of this section.
\end{enumerate}
We now present the experimental evaluations of all four experiments.
\subsubsection{Experiment 1}
For this experiment, we provided only $40$ snapshots out of the available $101$ snapshots for the spatial-temporal reconstruction via kernel-extended DMD using the Laplacian kernel and the GRBF kernel. After the reconstruction through both kernels, we then inspect and compare with each other for one single snapshot; here we have the snapshot $\#39$.
\begin{figure}[H]
\centering
\setkeys{Gin}{width=\linewidth}
\begin{subfigure}{0.35\textwidth}
\includegraphics{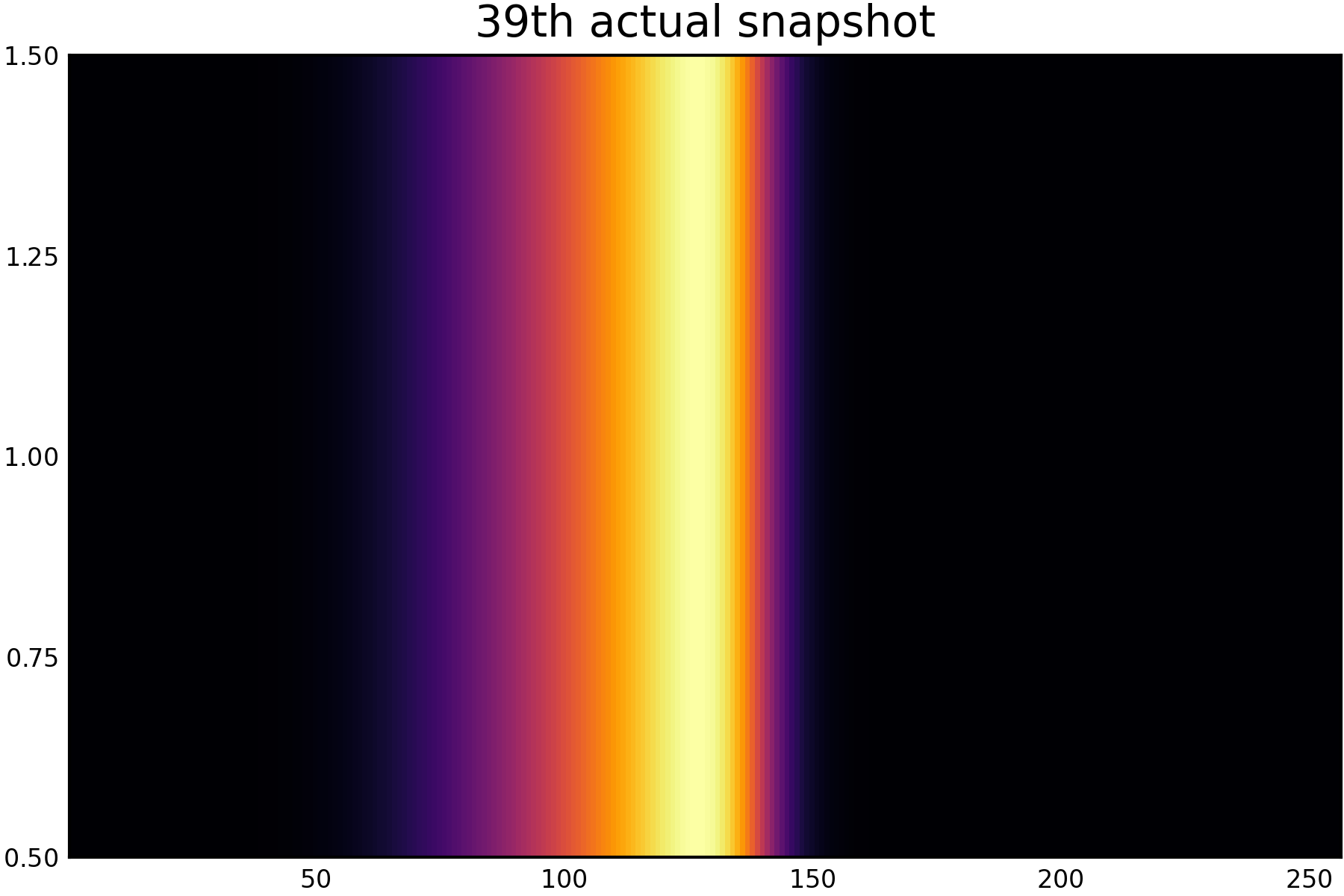}
\caption{Actual}
\label{fig:GroundTruth_Burger}
\end{subfigure}%
\hfil
\begin{subfigure}{0.35\textwidth}
\includegraphics{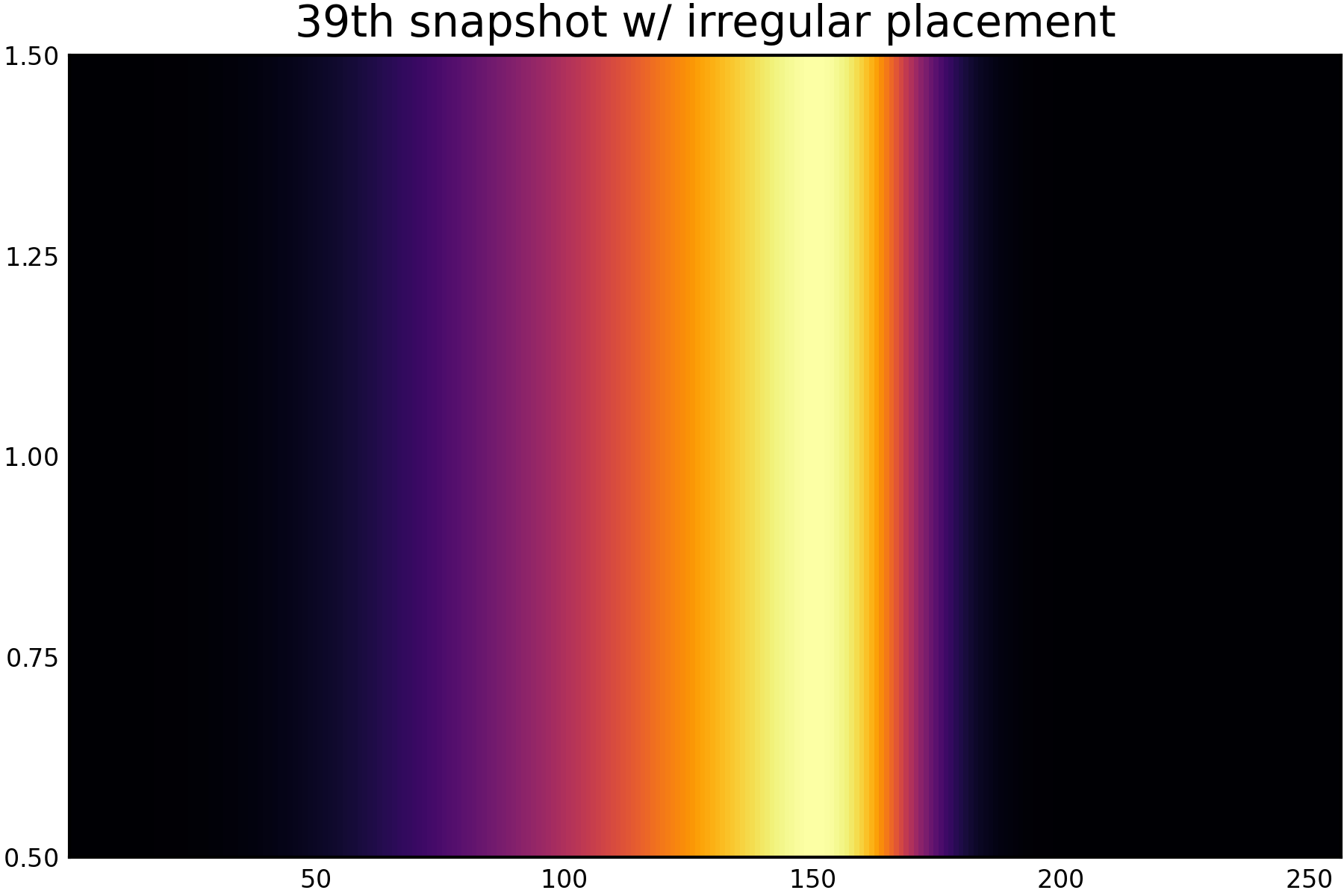}
\caption{Irregular and Sparse}
\label{fig:IrrSparse_Burger}
\end{subfigure}
\hfil
\begin{subfigure}{0.35\textwidth}
\includegraphics{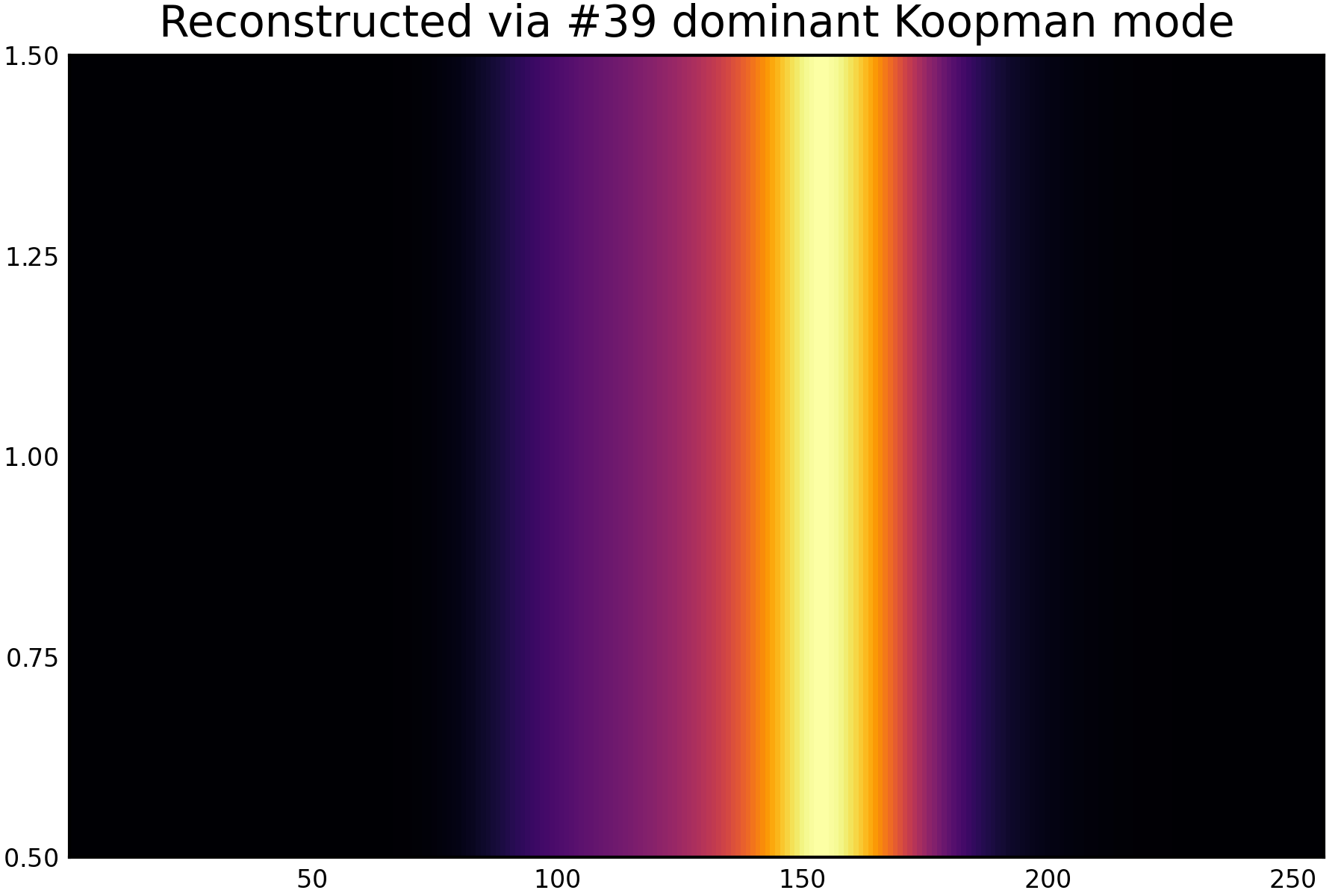}
\caption{Reconstruction by $\Lap$}
\label{fig:ReconLap_Burge}
\end{subfigure}
\hfil
\begin{subfigure}{0.35\textwidth}
\includegraphics{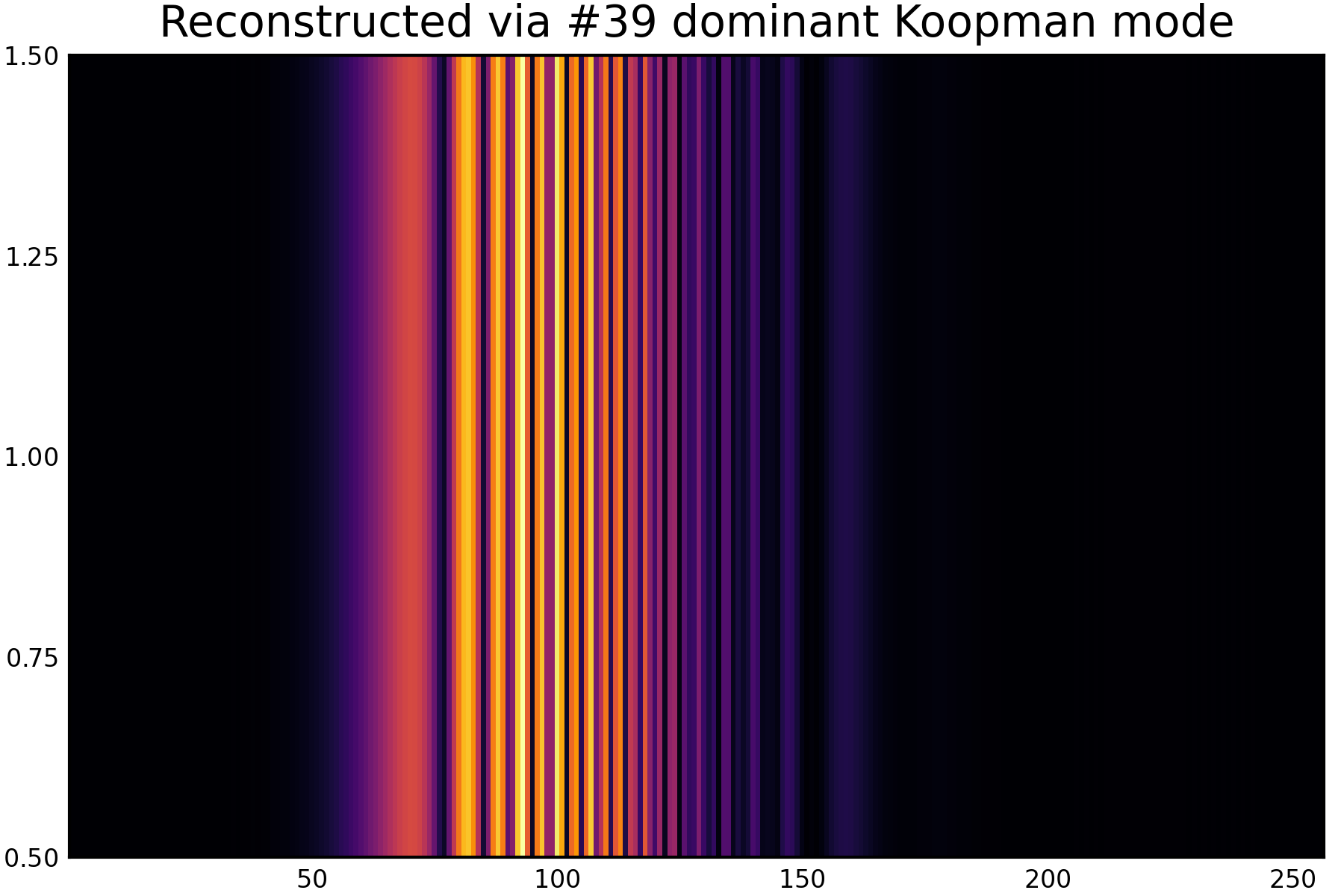}
\caption{Reconstruction by $\GRBF$}
\label{fig:ReconGRBF_Burger}
\end{subfigure}
\caption{Spatial-temporal reconstruction of the snapshot $\#39$ in the format of data matrix for non-linear Burger's Equation.}
\label{fig:DataMatrixComparisonExp1}
\end{figure}
After this, the following figure shows the spatial values of the total counts $256$, and the reconstruction for the snapshot $\#39$ by both kernels.  
\begin{figure}[H]
\centering
\setkeys{Gin}{width=\linewidth}
\begin{subfigure}{0.35\textwidth}
\includegraphics{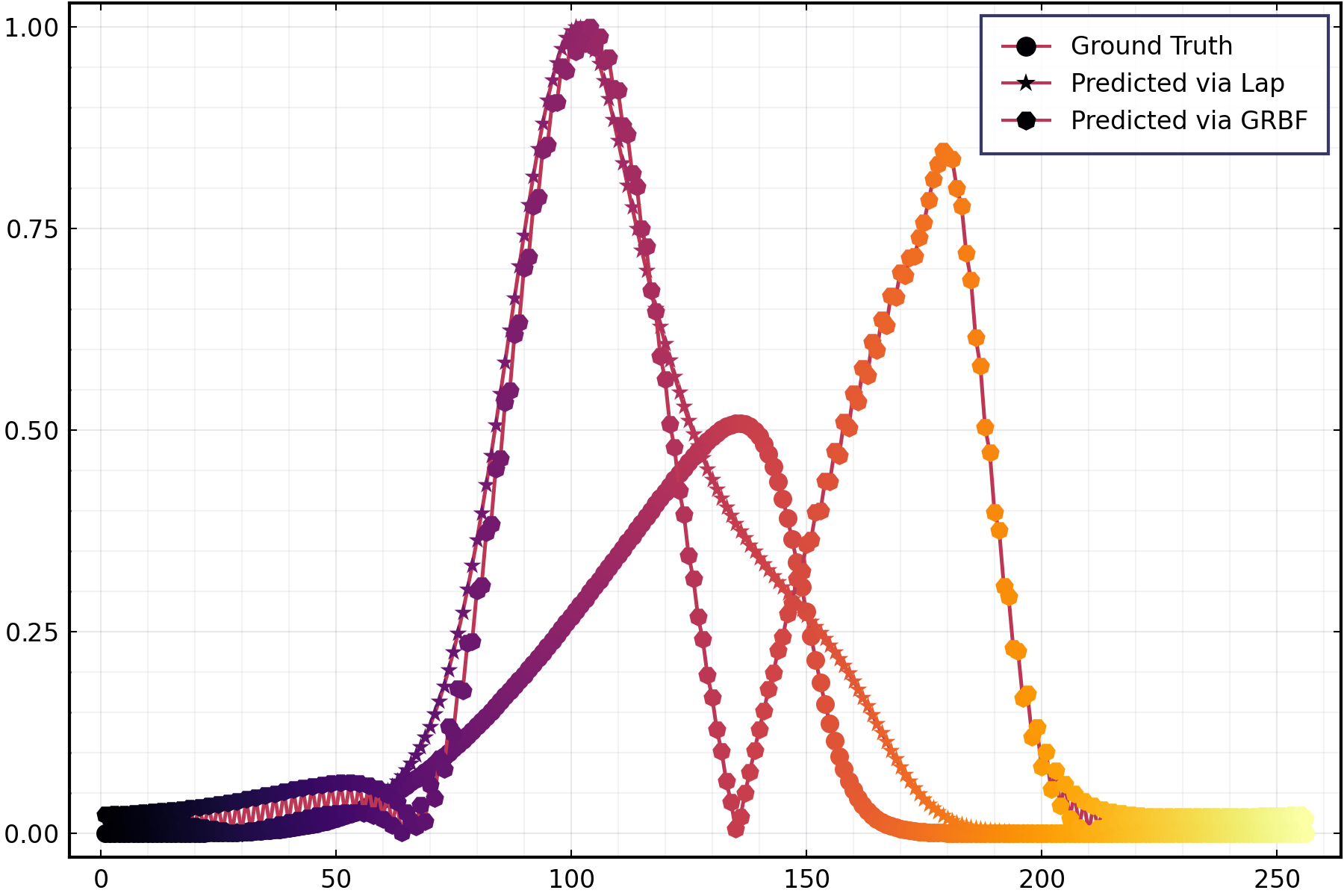}
\caption{Spatial values}
\label{fig:SpatialValues_Burger}
\end{subfigure}%
\hfil
\begin{subfigure}{0.35\textwidth}
\includegraphics{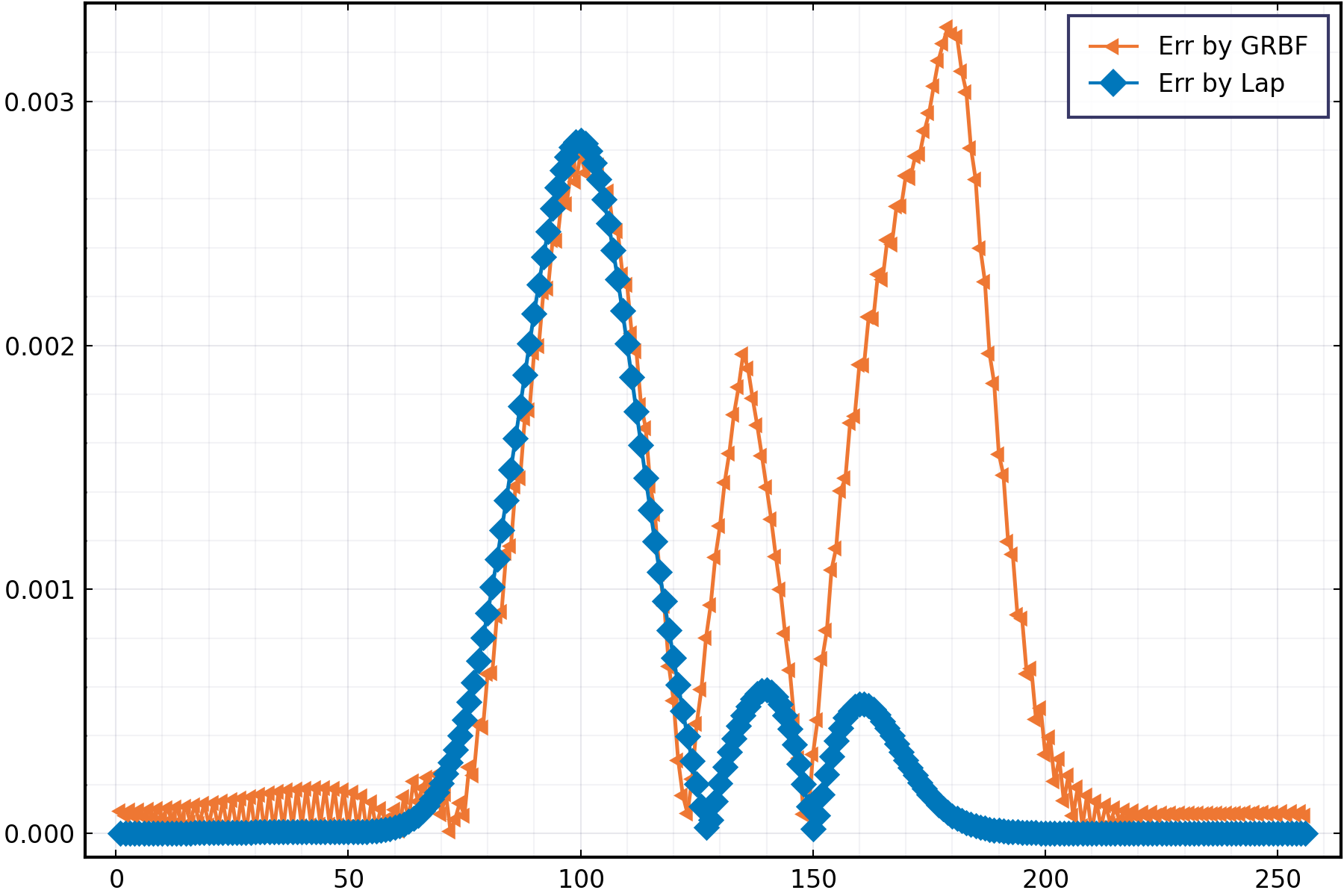}
\caption{EWE plot}
\label{fig:Errs_Burger}
\end{subfigure}
\caption{Spatial-values comparison of the snapshot $\#39$ for non-linear Burger's Equation along with EWE plot.}
\label{fig:SpatialValuesComparisonExp1}
\end{figure}
\subsubsection{Experiment 2}
For this experiment, we provided only $100$ snapshots out of the available $151$ snapshots for the spatial-temporal reconstruction via kernel-extended DMD using the Laplacian kernel and the GRBF kernel. After the reconstruction through both kernels, we then inspect and compare with each other for one single snapshot; here we have the snapshot $\#99$.
\begin{figure}[H]
\centering
\setkeys{Gin}{width=\linewidth}
\begin{subfigure}{0.3\textwidth}
\centering
\includegraphics[scale=1.1]{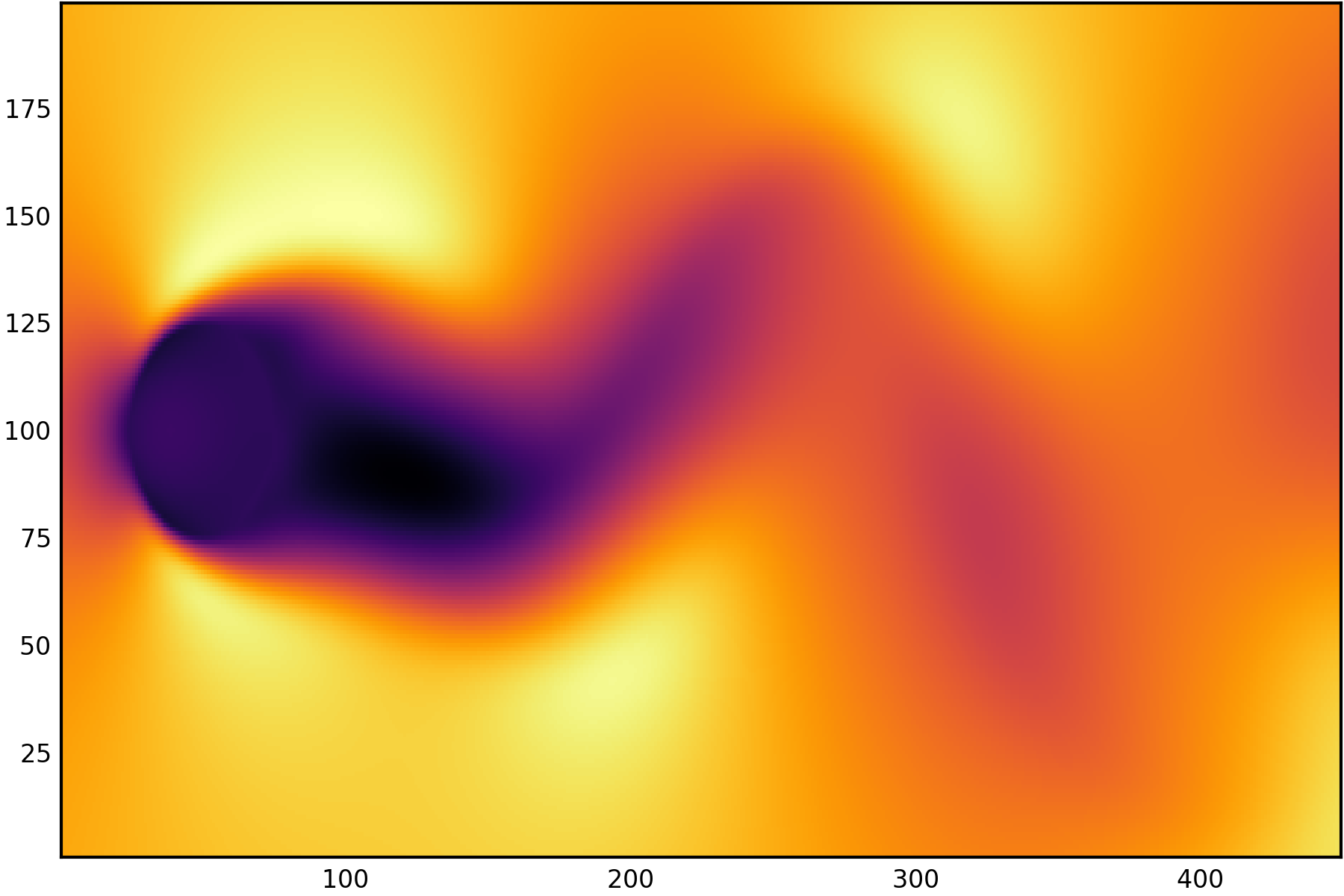}
\caption{Actual}
\label{fig:GroundTruth_fluidflow}
\end{subfigure}%
\hfil
\begin{subfigure}{0.3\textwidth}
\centering
\includegraphics[scale=1.1]{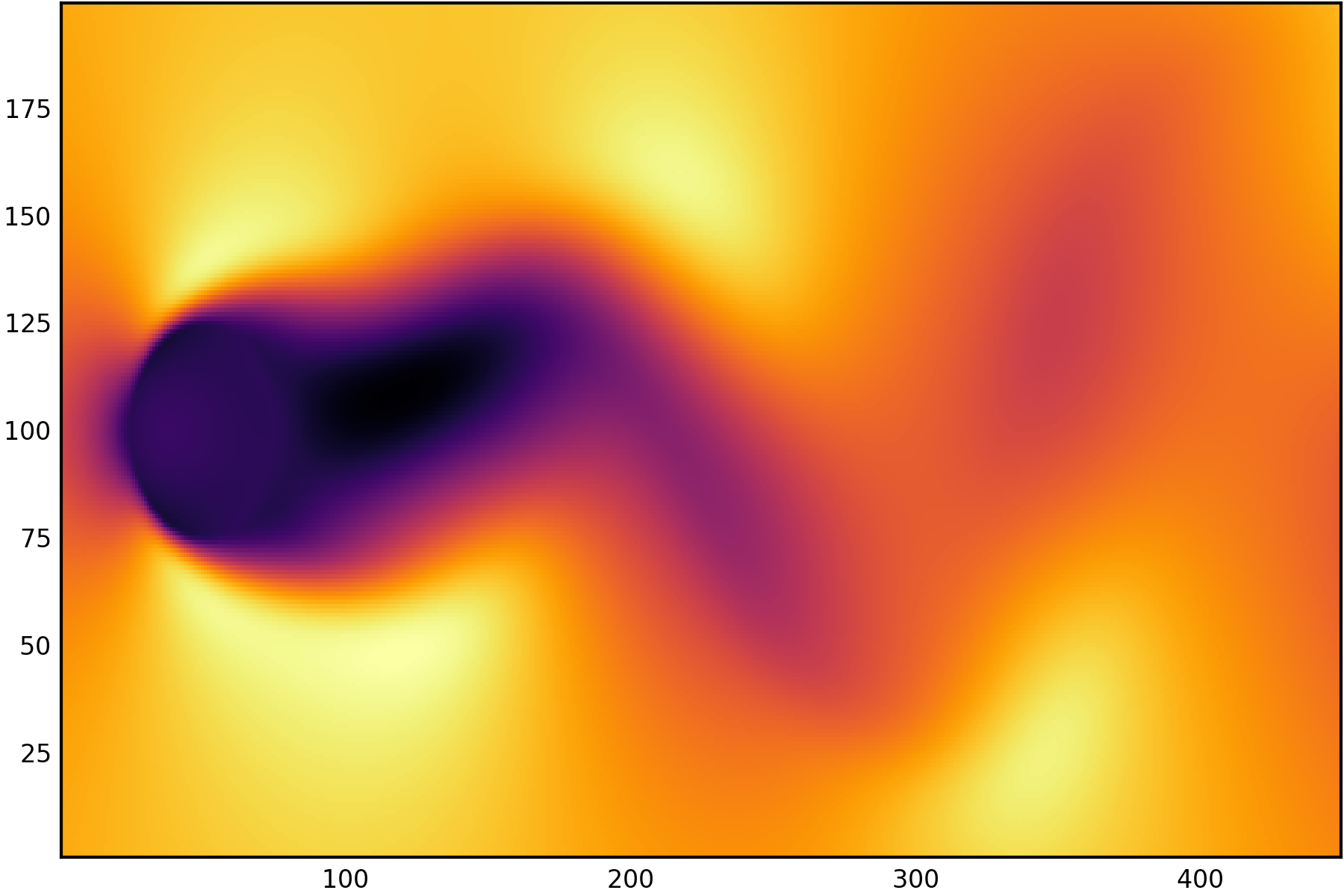}
\caption{Irregular \& Sparse}
\label{fig:IrrSparse_fluidflow}
\end{subfigure}
\hspace{5em}
\begin{subfigure}{0.3\textwidth}
\centering
\includegraphics[scale=1.1]{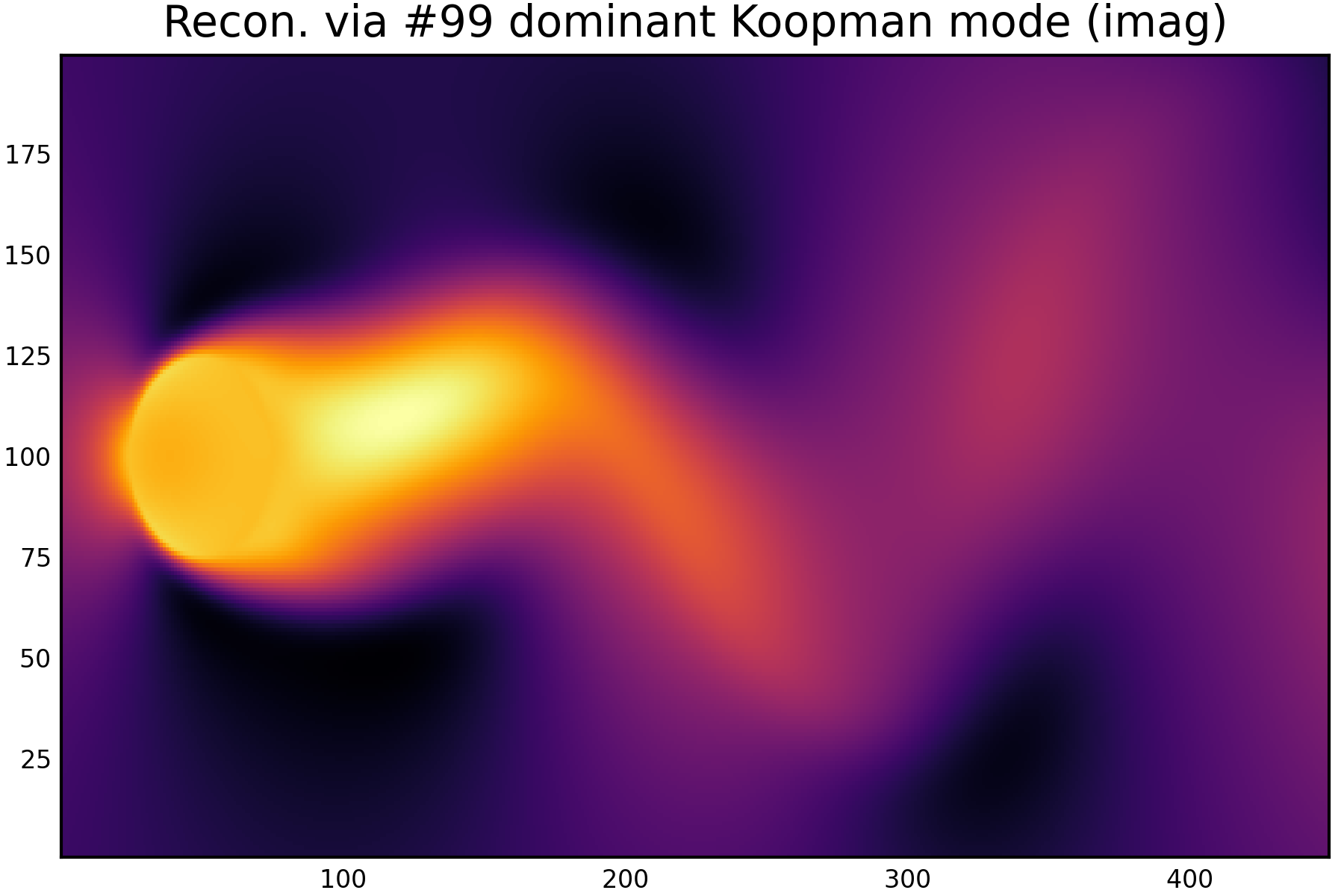}
\caption{Reconstruction by $\Lap$ 
}
\label{fig:ReconLap_fluidflow_imag}
\end{subfigure}
\hfil
\begin{subfigure}{0.3\textwidth}
\centering
\includegraphics[scale=1.1]{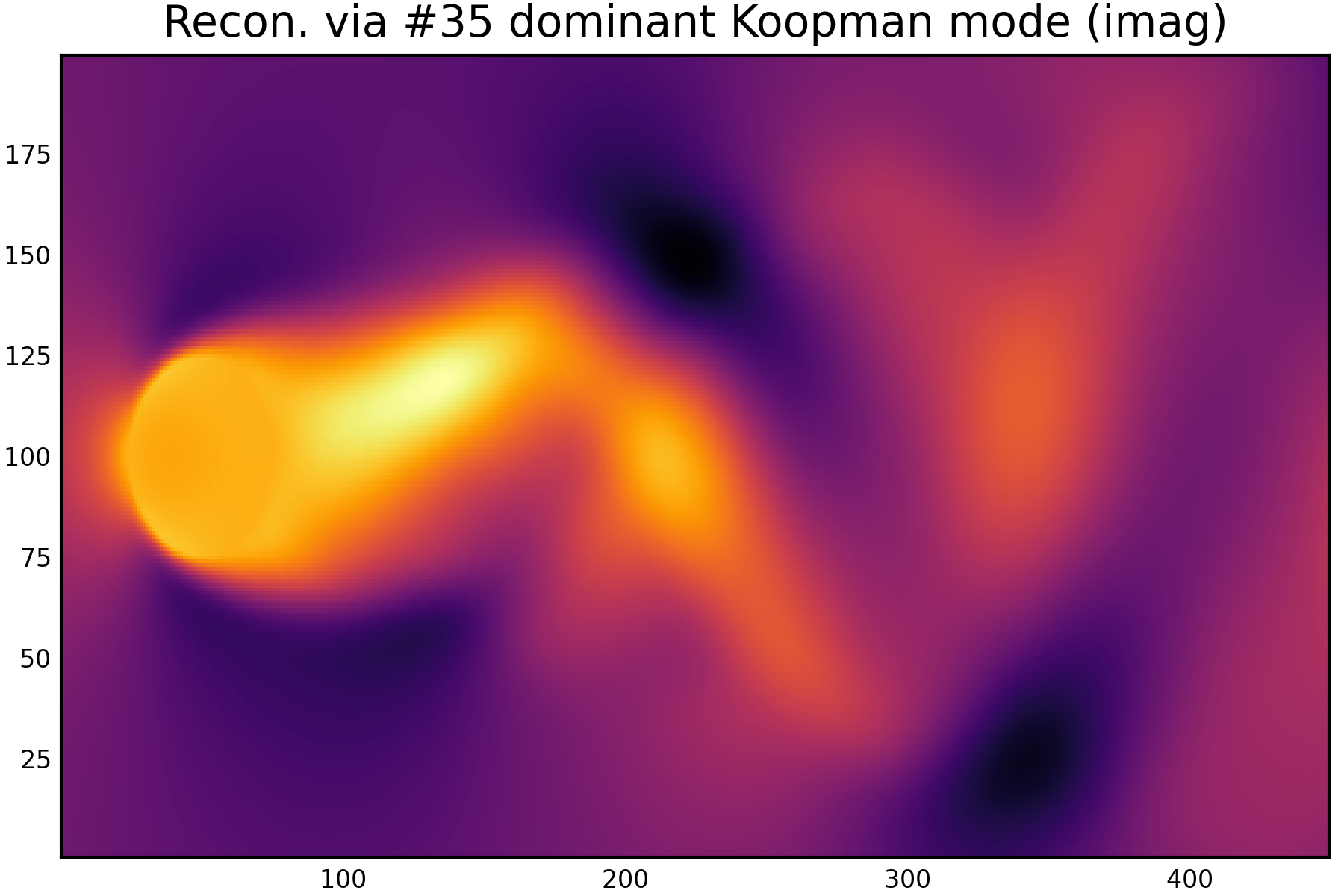}
\caption{Reconstruction by $\GRBF$ 
}
\label{fig:ReconGRBF_fluidflow_imag}
\end{subfigure}
\caption{Spatial-temporal reconstruction of the snapshot $\#99$ in the format of data matrix for periodic fluid flow across cylinder.}
\label{fig:DataMatrixComparisonExp2}
\end{figure}
After this, the following figure shows the spatial values of the considered counts $100$ out of the total counts available $89,351$, and the reconstruction for the snapshot $\#99$ by both kernels. We did this to get better visualization of our plot results.
\begin{figure}[H]
\centering
\footnotesize
\setkeys{Gin}{width=\linewidth}
\begin{subfigure}{0.35\textwidth}
\includegraphics[width=1\linewidth]{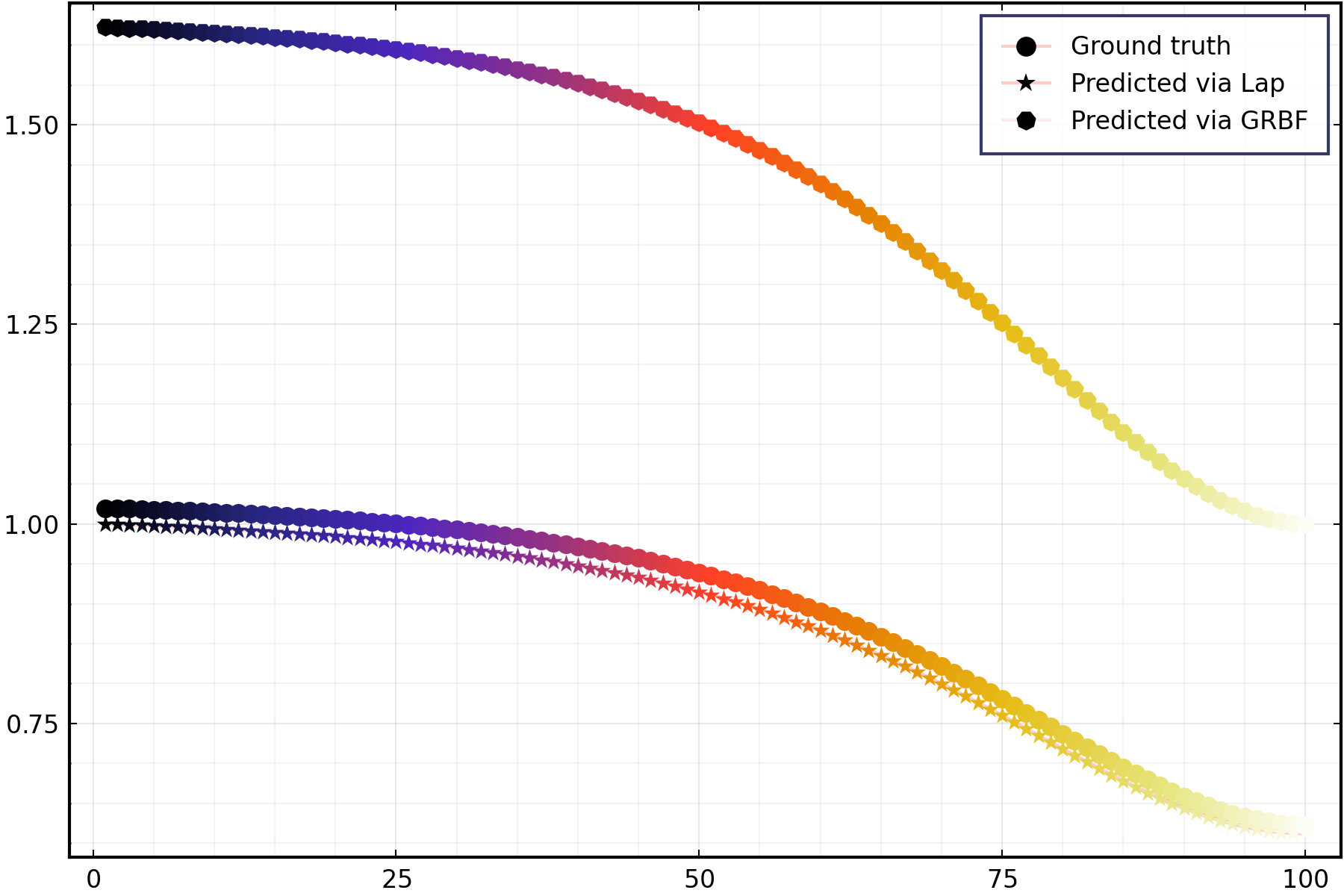}
\caption{$100-$ Spatial values}
\label{fig:SpatialValues_fluid_flow}
\end{subfigure}%
\hfil
\begin{subfigure}{0.35\textwidth}
\includegraphics[width=1\linewidth]{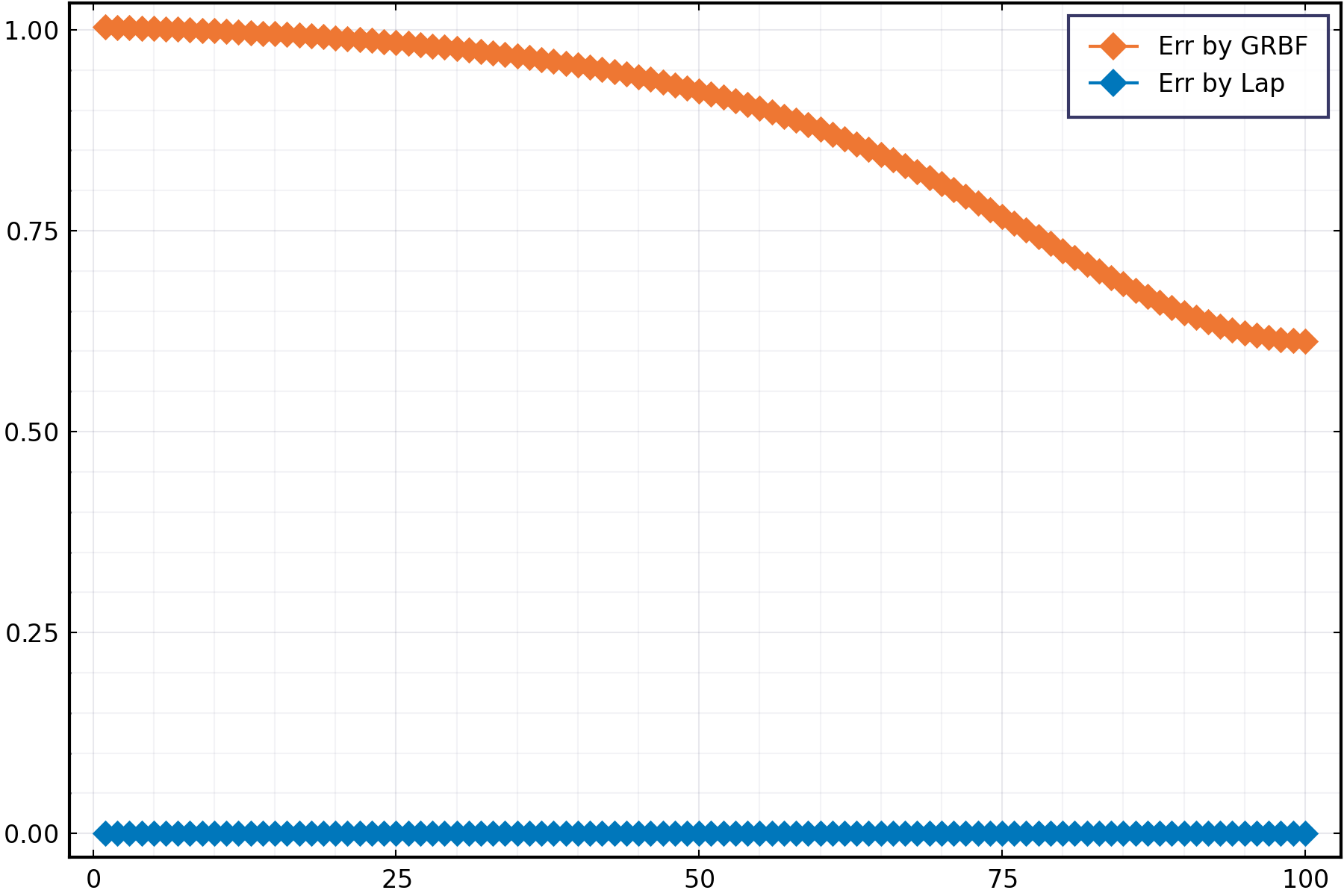}
\caption{EWE plot}
\label{fig:Errs_fluid_flow}
\end{subfigure}
\caption{Spatial-values comparison of the snapshot $\#99$ for the fluid flow across cylinder.}
\label{fig:SpatialValuesComparisonExp2}
\end{figure}
\subsubsection{Experiment 3}
Dataset for this experiment built with an initial condition given as $[-1.8760,1.7868]^\top$. This experiment belongs to the natural setting of the phase-portrait reconstruction problem which we have adopted to test the proposed Lap-KeDMD algorithm in terms of spatial modes reconstruction. Follow the supplementary materials as well.
\begin{figure}[H]
\centering
\setkeys{Gin}{width=\linewidth}
\begin{subfigure}{0.225\textwidth}
\includegraphics[scale=1]{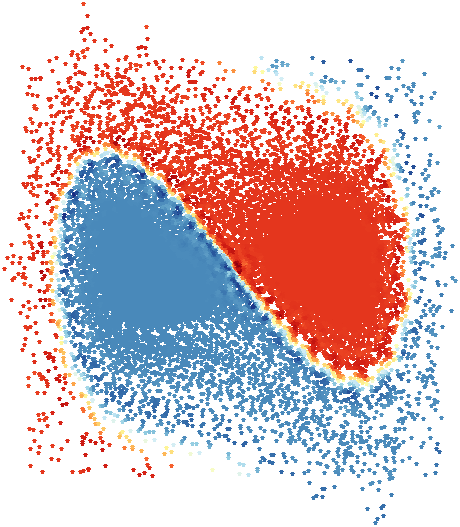}
\caption{Actual}
\label{fig:GroundTruth_Duffing}
\end{subfigure}%
\begin{subfigure}{0.225\textwidth}
\includegraphics{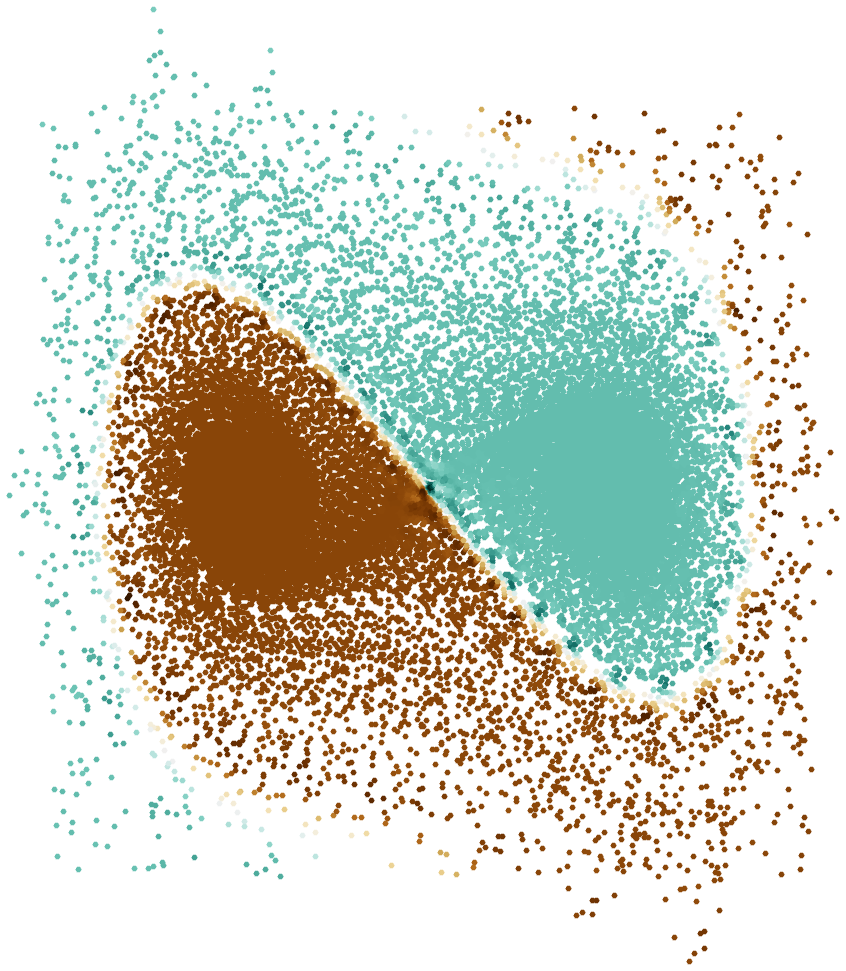}
\caption{Irregular \& Sparse}
\label{fig:IrrSparse_Duffing}
\end{subfigure}
\begin{subfigure}{0.225\textwidth}
\includegraphics{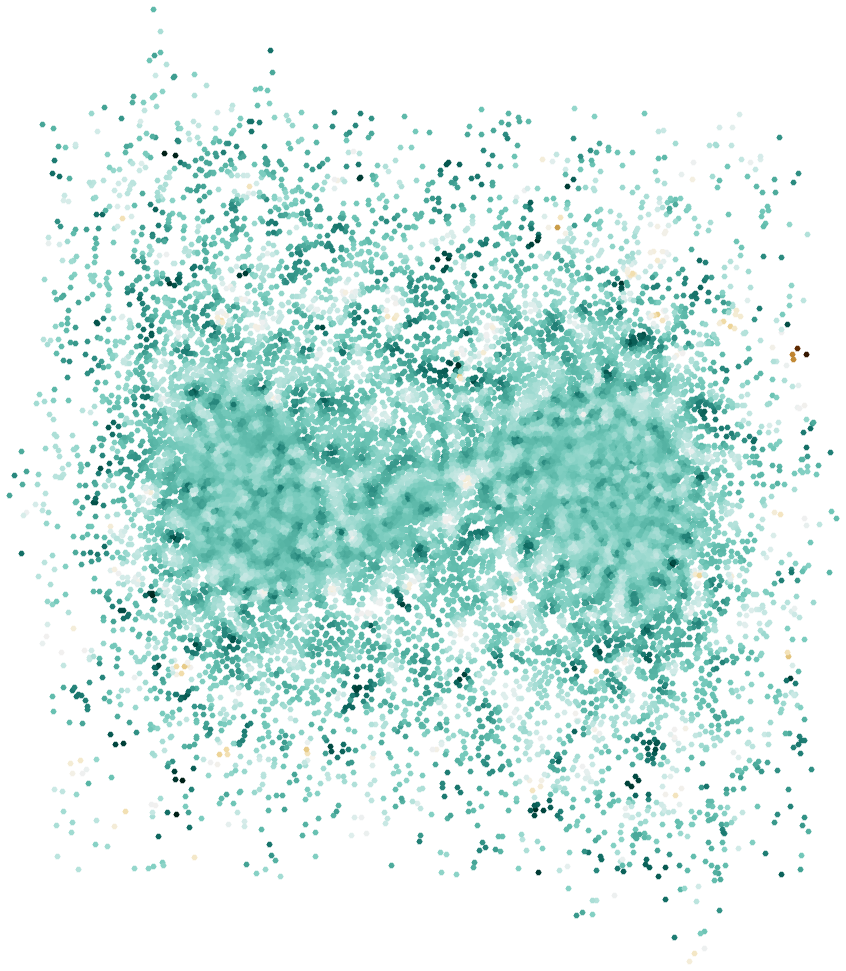}
\caption{Recon. by $\Lap$}
\label{fig:ReconLap_Duffing}
\end{subfigure}
\begin{subfigure}{0.225\textwidth}
\includegraphics{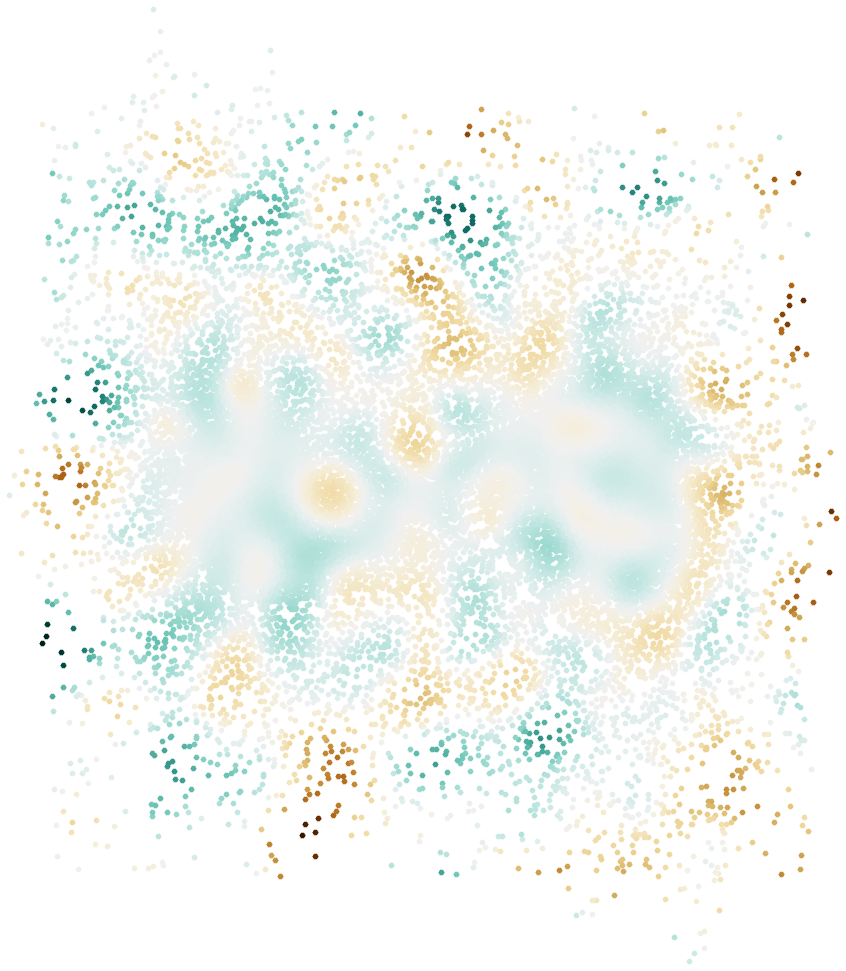}
\caption{Recon. by $\GRBF$}
\label{fig:ReconGRBF_Duffing}
\end{subfigure}
\caption{Reconstruction of the snapshot $\#39$ for the chaotic Duffing's Oscillator.}
\label{fig:DataMatrixComparisonExp3}
\end{figure}
\subsubsection{Experiment 4}
\begin{figure}[H]
    \centering
    \begin{subfigure}{0.5\textwidth}
    \centering
    \includegraphics[width=0.85\linewidth]{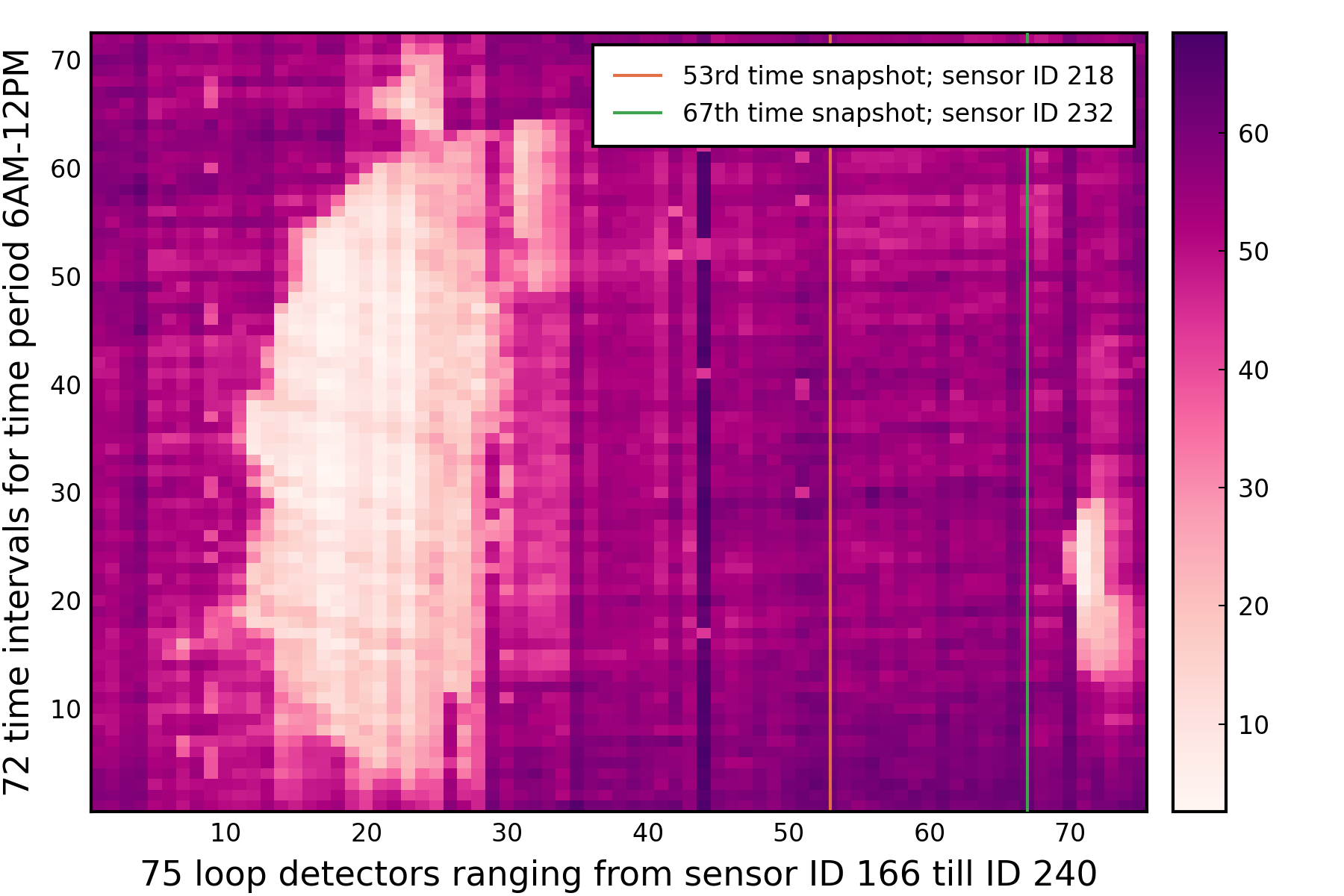}
    \caption{Seattle I-5 Freeway Traffic data}
    \label{fig:SeattleI5dataheatmap}
    \end{subfigure}
    \hfil
    \begin{subfigure}{0.45\textwidth}
    \centering
    \includegraphics[width=1.\linewidth]{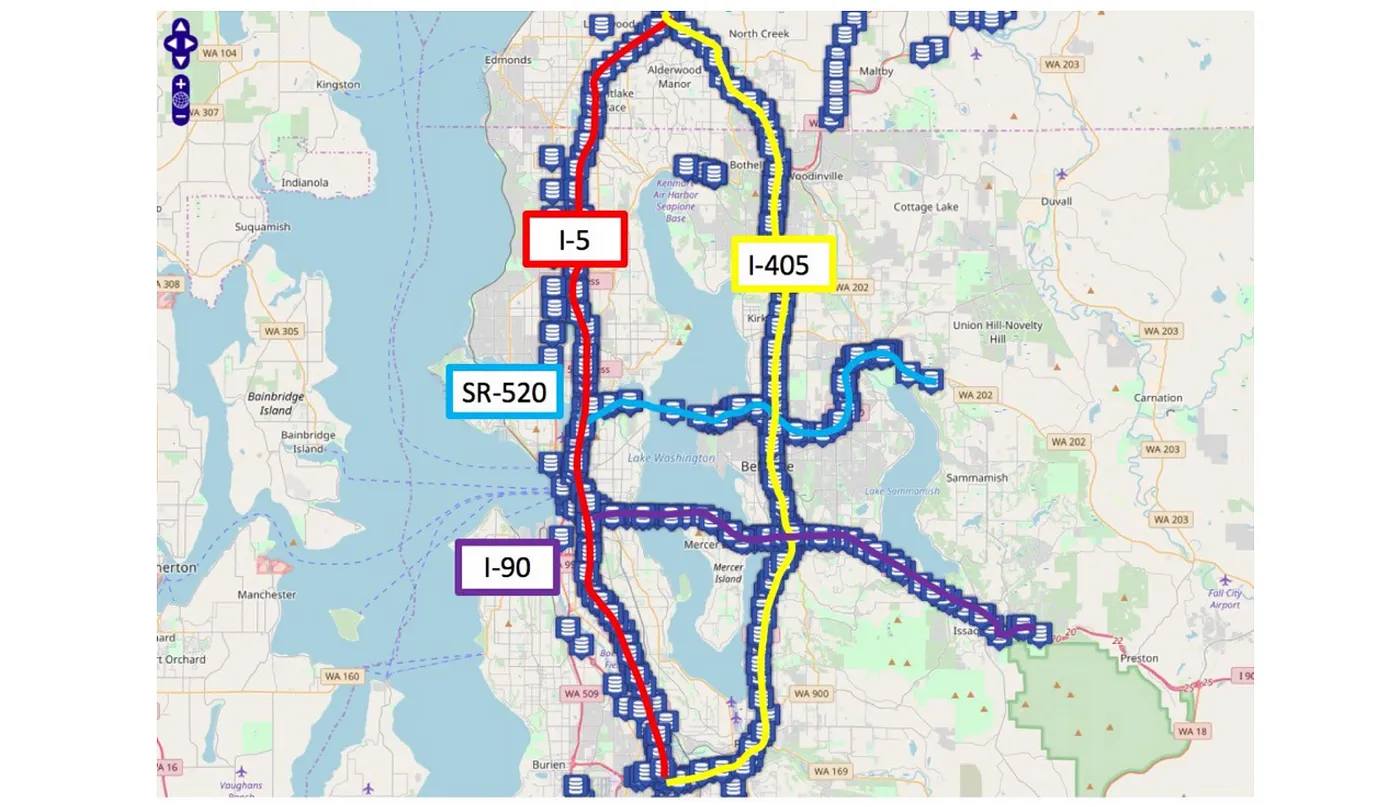}
    \caption{Seattle I-5 geographical map}
    \label{fig:GeograhicMapSeattle}
    \end{subfigure}
    \caption{Experimental information for the Seattle I-5 Traffic data set.
    }
    \label{fig:Seattle_218+232_r}
\end{figure}

For the Seattle I-5 freeway traffic speed data, there are $75$ loop detectors that record vehicle speed in a $72$ time interval frame during a rush hour period from 6:00 am to 12:00 pm; follow \cite{cui2018deep,cui2019traffic} for more details. These $75$ loop detectors present on the highway have sensor IDs as $166$ to $240$ sequentially; follow \Cref{fig:SeattleI5dataheatmap} for the traffic speed distribution from this data set. In addition, we also provide the geographic map of the Seattle I-5 freeway that we collected from \cite{cui2018deep,cui2019traffic} in \Cref{fig:GeograhicMapSeattle}. In this map, we highlight the inductive locations of each loop detector in \textbf{\emph{\textcolor{blue}{blue icon}}} that indicate the loop detectors at all the traffic sensors present at all mileposts. 
We provide the histogram distribution of speed data for the Seattle I5 freeway collected through loop detector IDs $166$\textsuperscript{th}, $186$\textsuperscript{th} and $226$\textsuperscript{th} followed by their respective graphs.
\begin{figure}[H]
    \centering
    \setkeys{Gin}{width=\linewidth}
    \begin{subfigure}{0.3\textwidth}
    \centering
    \includegraphics[width=\linewidth]{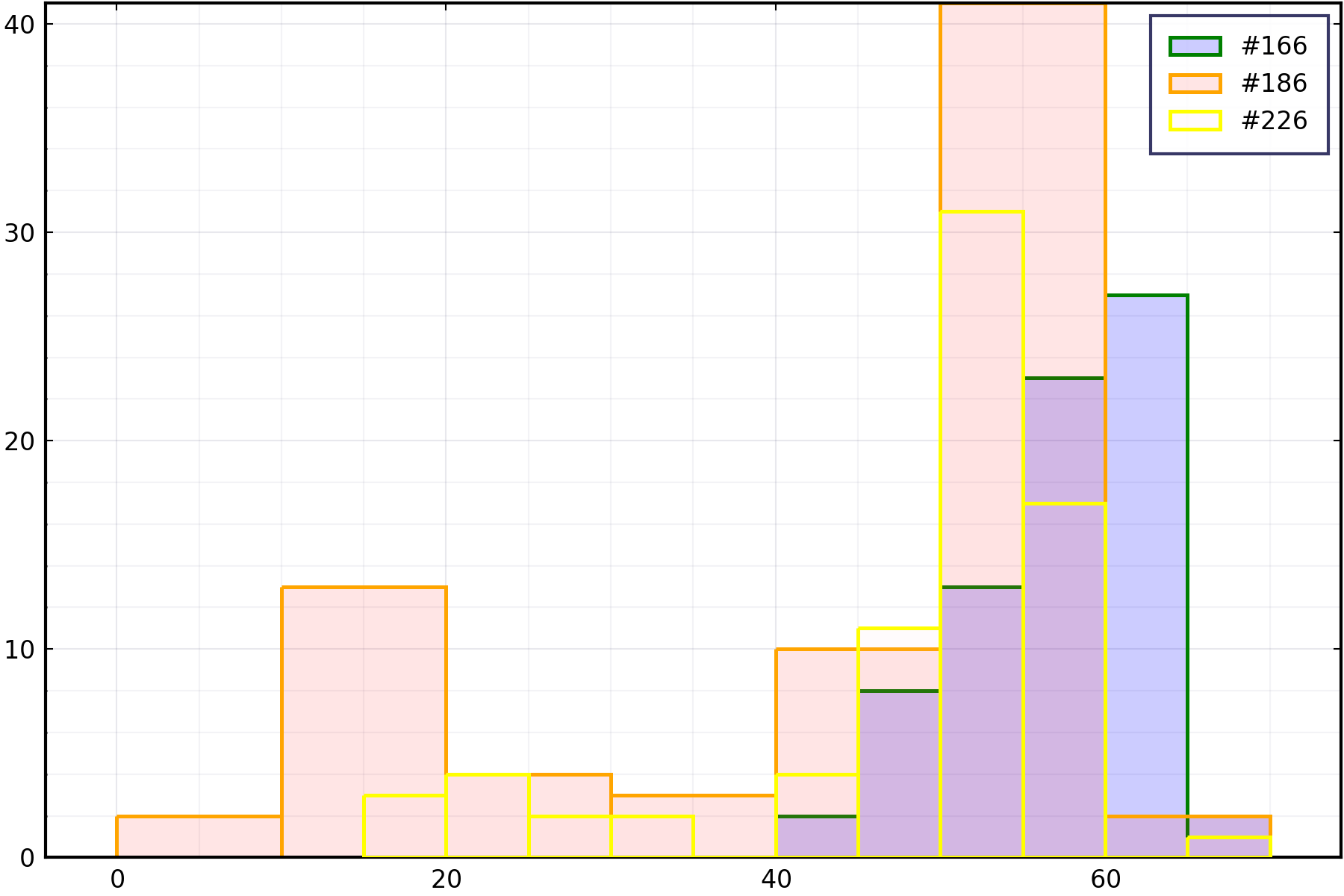}
    \caption{Speed vs. histogram bins}
    \end{subfigure}
    \hfil
    \begin{subfigure}{0.3\textwidth}
    \centering
    \includegraphics[width=\linewidth]{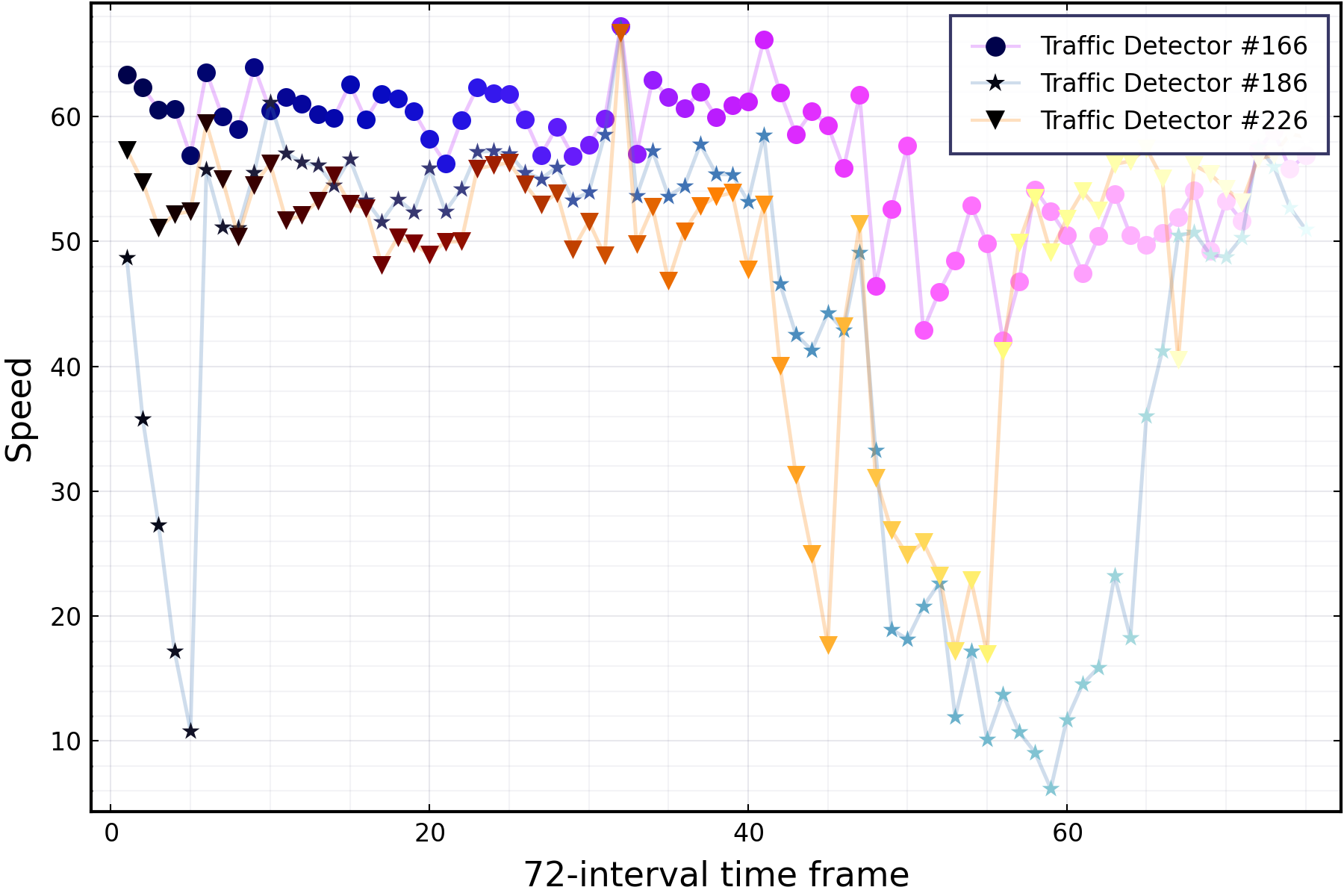}
    \caption{Speed Plot}
    \end{subfigure}
    \caption{Histogram distribution and plot of speed data on three loop detector ID's on Seattle I5 freeway.
    }
    \label{fig:Seattle_traffic_data}
\end{figure}
We now present the results of the spatial-temporal reconstruction and their error plots of Seattle I5 traffic data for the sensor ID's $218$\textsuperscript{th}
through the kernel-extended DMD by both kernels. The speed data collected in the ID of the Seattle I-5 highway sensor $218$\textsuperscript{th} corresponds to the $53$\textsuperscript{rd} time snapshot.
\begin{figure}[H]
        \centering
        \begin{subfigure}{0.3\textwidth}
        \centering
        \includegraphics[width=1\linewidth]{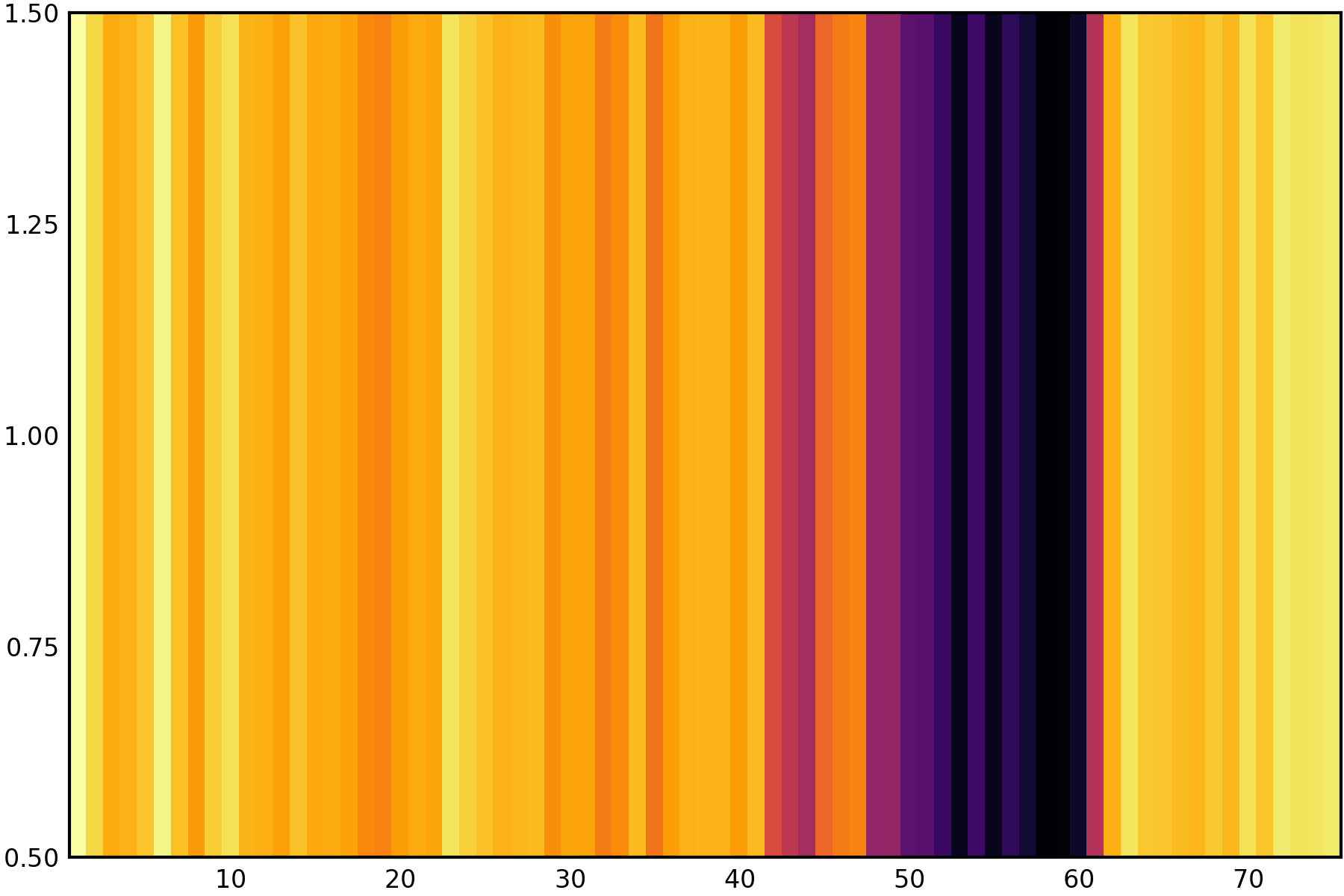}
        \caption{Actual}
        \label{fig:actual_snap_53}
        \end{subfigure}
        \begin{subfigure}{0.3\textwidth}
        \centering
        \includegraphics[width=1\linewidth]{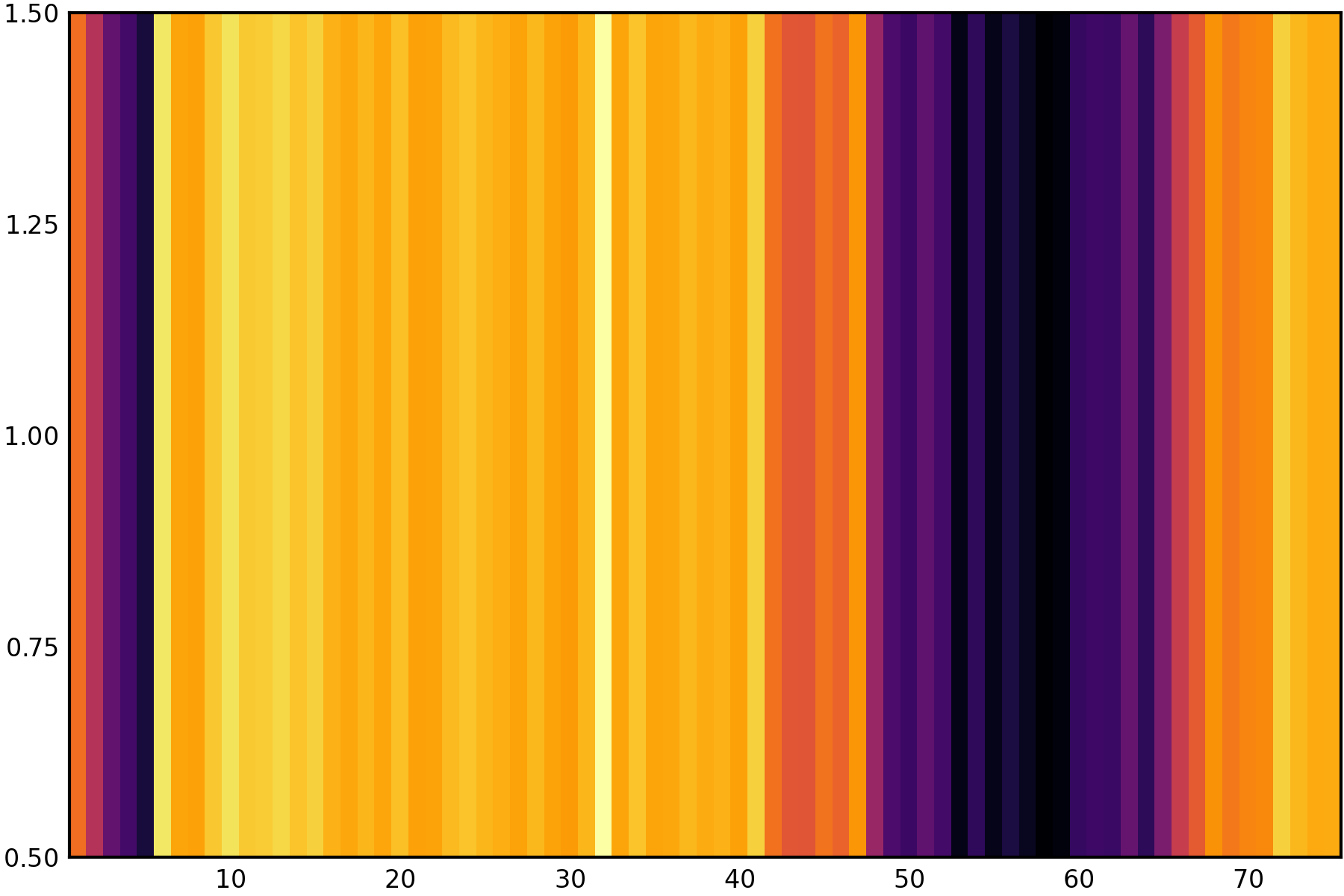}
        \caption{Reconstruction by $\Lap$}
        \label{fig:dominant_53_lap}
        \end{subfigure}
        \begin{subfigure}{0.3\textwidth}
        \centering
        \includegraphics[width=1\linewidth]{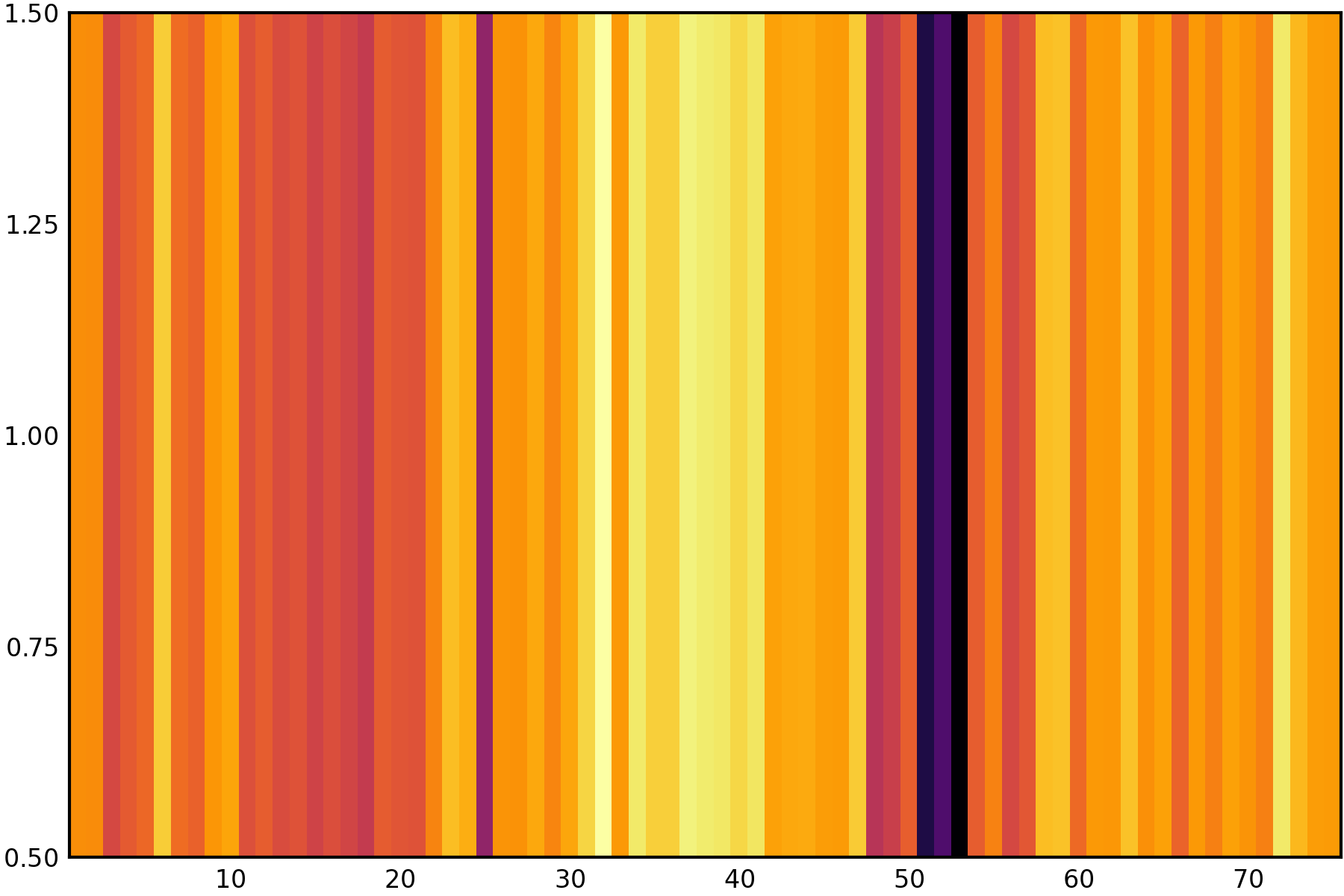}
        \caption{Reconstruction by $\GRBF$}
        \label{fig:dominant_53_grbf}
        \end{subfigure}
        \caption{Reconstruction for the Seattle I-5 highway sensor $218$\textsuperscript{th}. 
        }
    \label{fig:STRecon_Traffic_53}
    \end{figure}
\begin{figure}[H]
    \centering
    \begin{subfigure}{0.3\textwidth}
    \centering
    \includegraphics[width=1\linewidth]{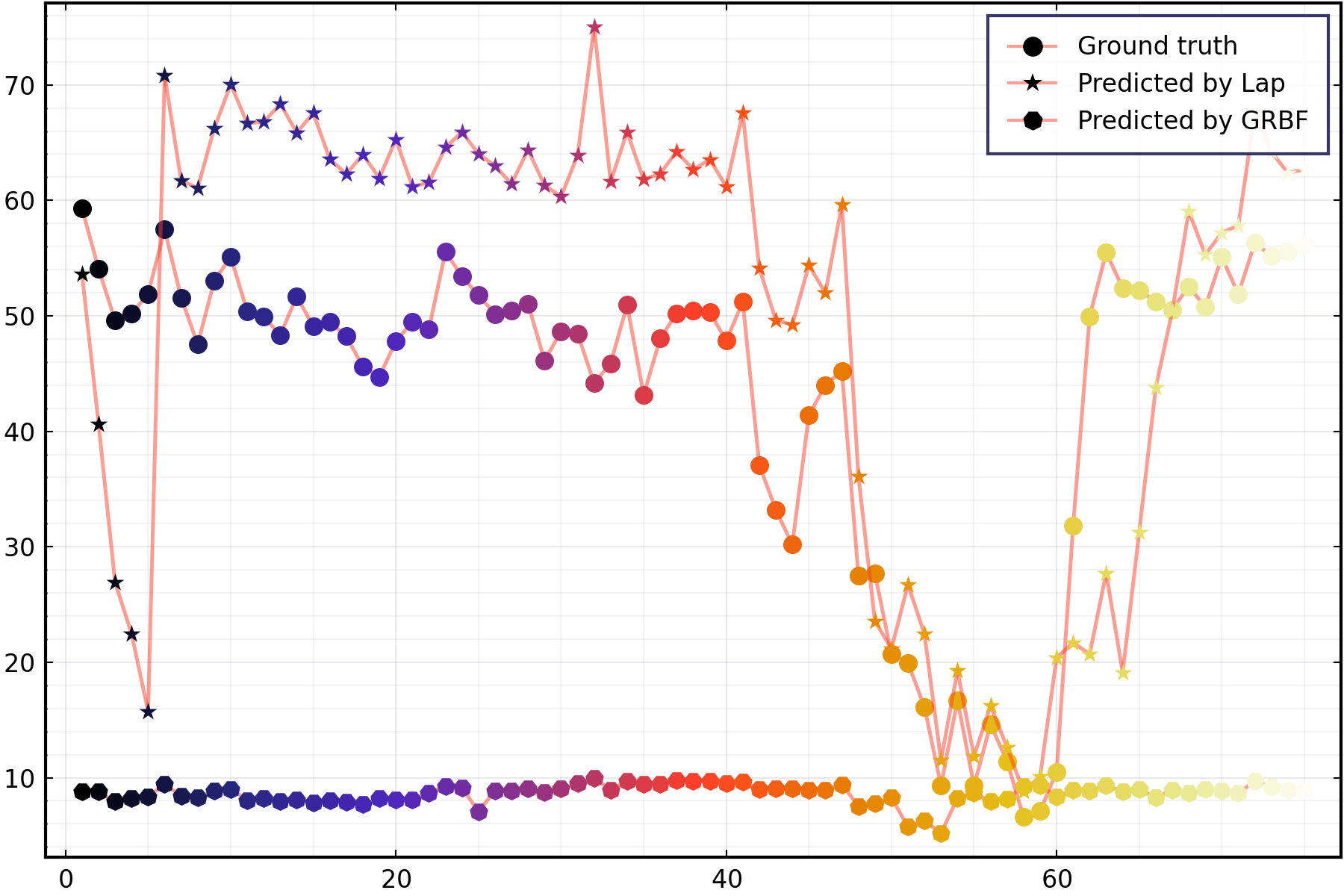}
    \caption{Spatial values}
    \end{subfigure}
    \hfil
    \begin{subfigure}{0.3\textwidth}
    \centering
    \includegraphics[width=1\linewidth]{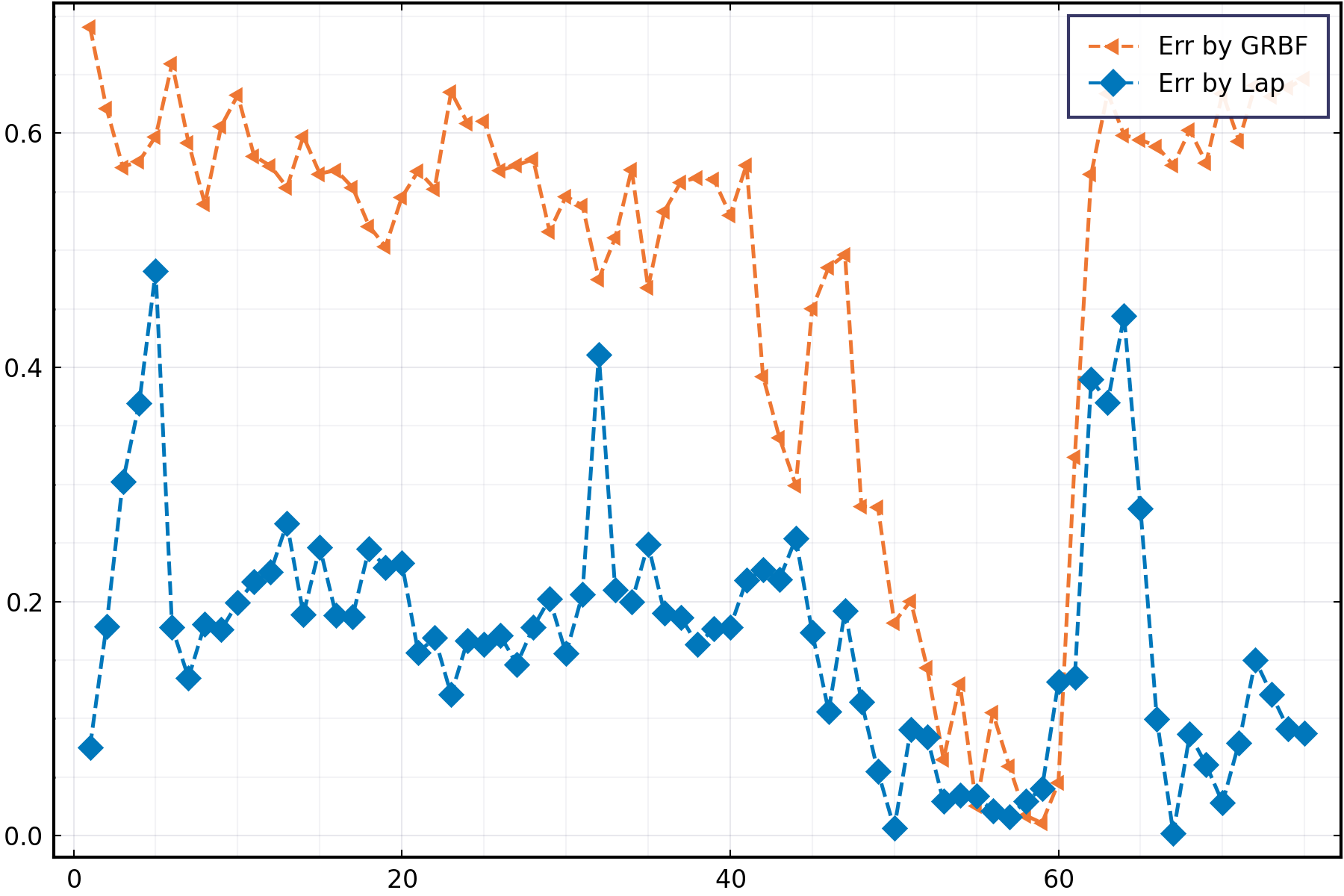}
    \caption{EWE plot}
    \end{subfigure}
    \caption{Reconstruction result for the Seattle I-5 freeway traffic sensor ID $218$\textsuperscript{th}.
    }
    \label{fig:Seattle_53_comparison+err}
\end{figure}
\subsubsection{Experiment 5}
For the following experiment of the Lorenz 1963 three-dimensional dynamical system, we now provide the first three important characteristic plots of the system.
\begin{figure}[H]
    \centering
    \begin{subfigure}{0.3\textwidth}
    \centering
        \includegraphics[width=.7\linewidth]{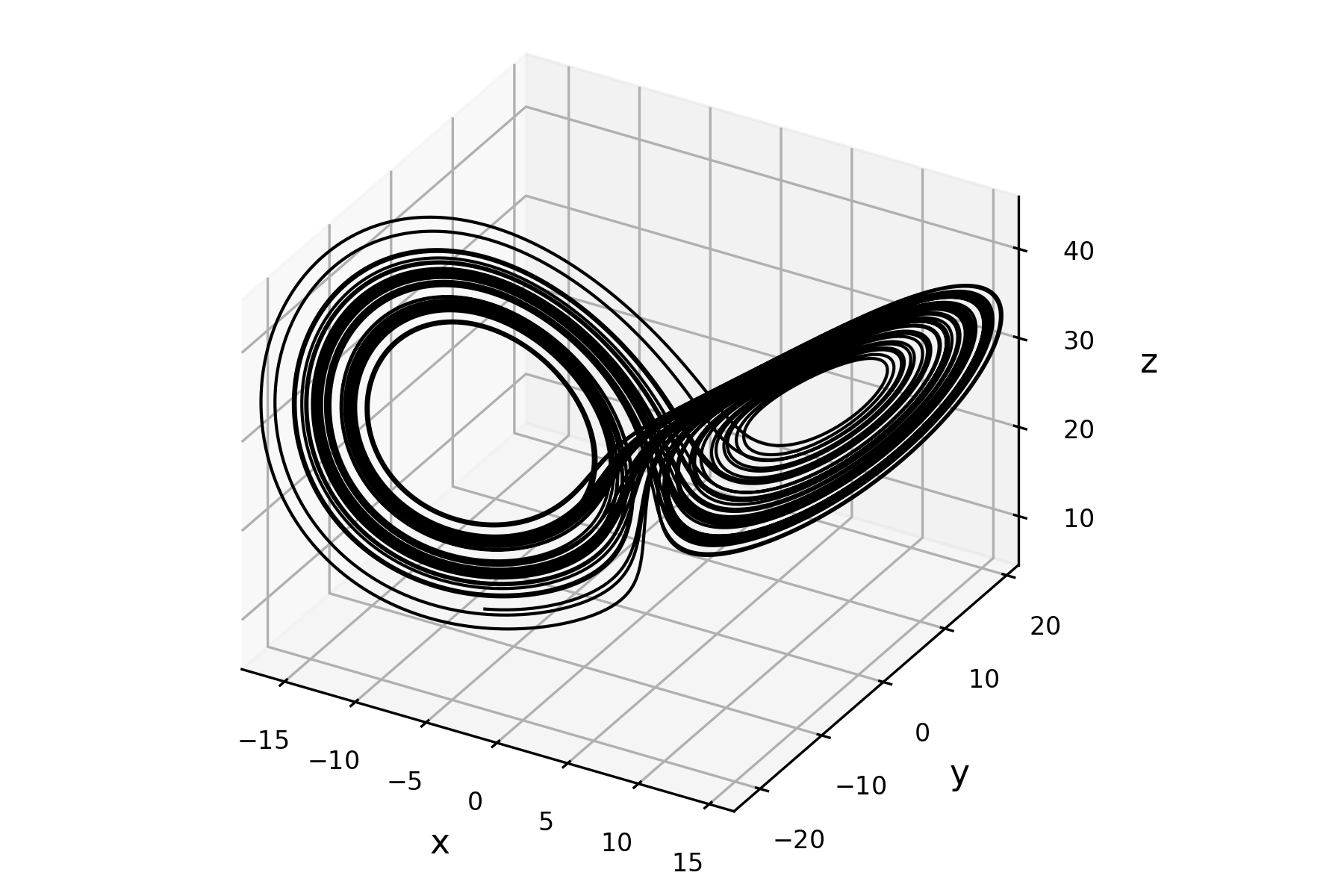}
        \caption{Lorenz 1963}
    \end{subfigure}
    \begin{subfigure}{0.3\textwidth}
    \centering
        \includegraphics[width=.7\linewidth]{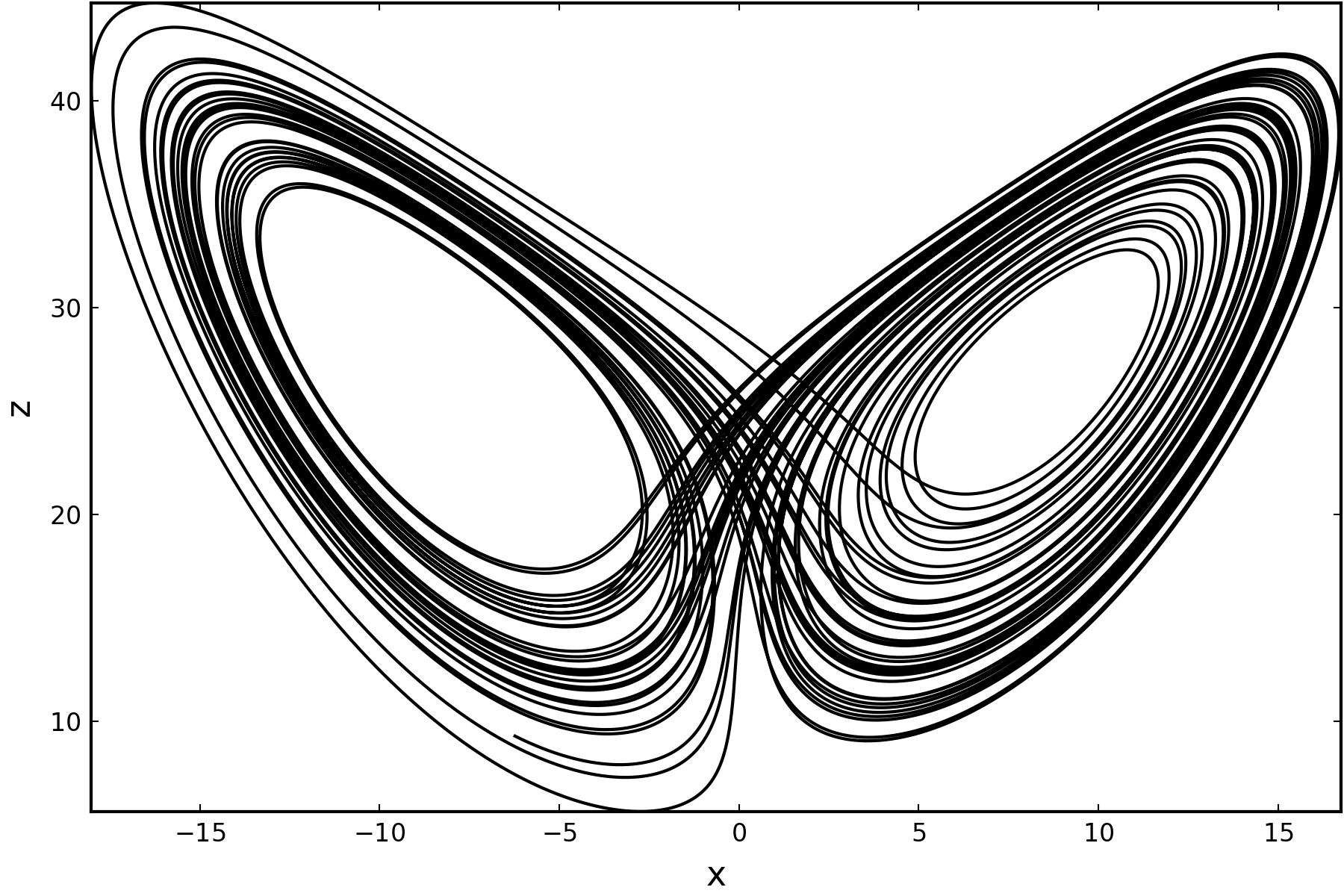}
        \caption{2D $xz$ attractor}
    \end{subfigure}
    \begin{subfigure}{0.3\textwidth}
    \centering
        \includegraphics[width=.7\linewidth]{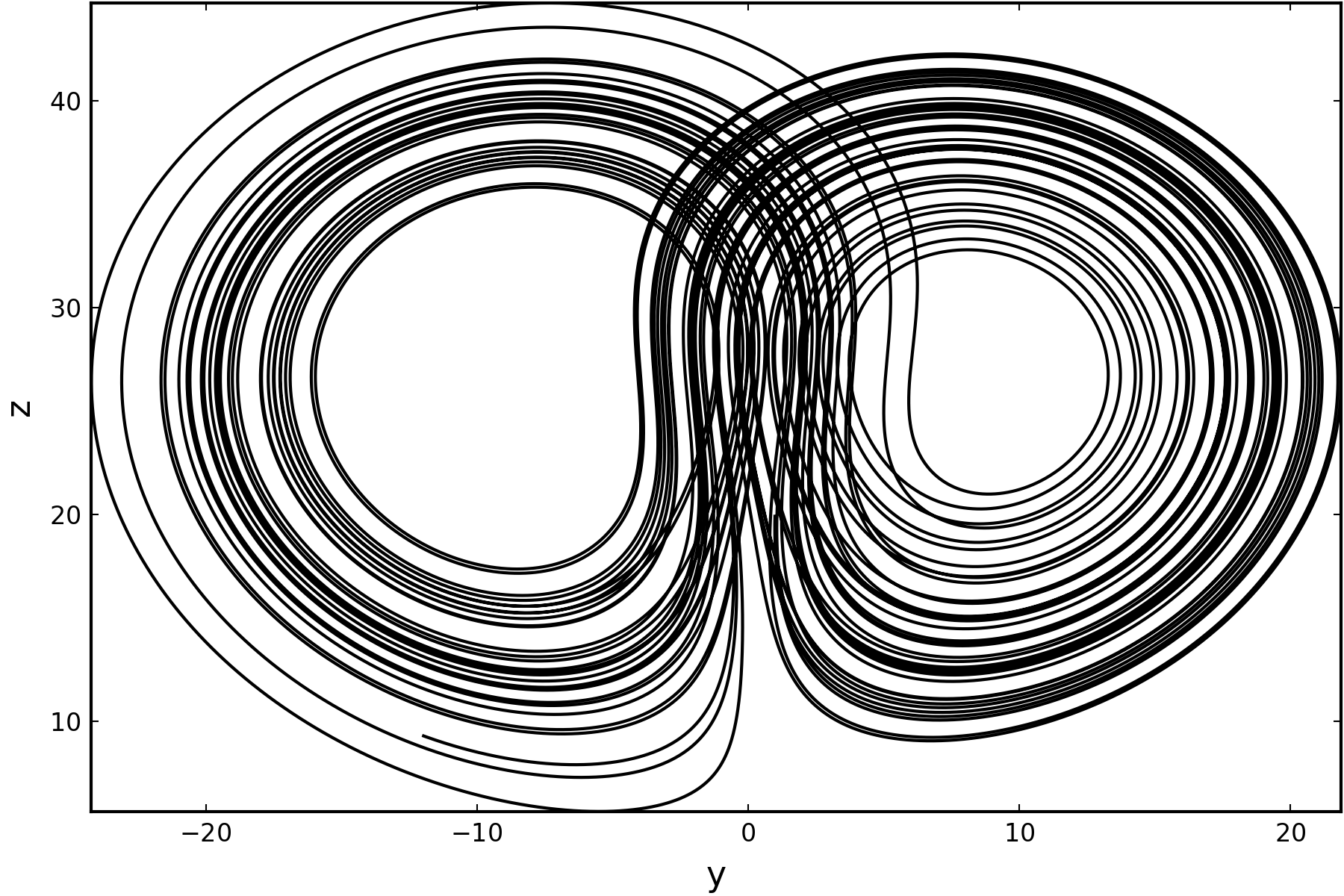}
        \caption{2D $yz$ attractor}
    \end{subfigure}
    \begin{subfigure}{0.3\textwidth}
    \centering
        \includegraphics[width=.7\linewidth]{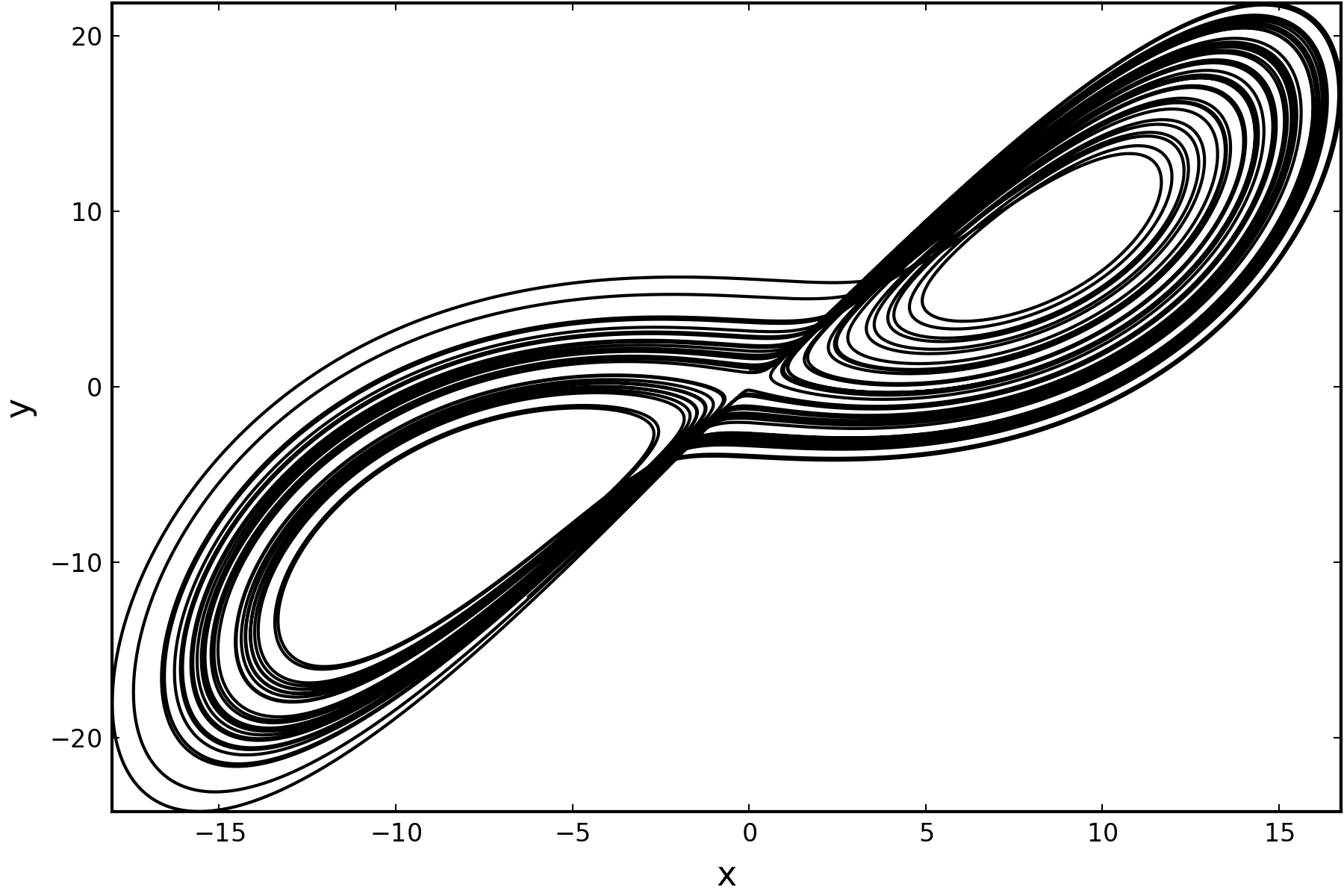}
        \caption{2D $xy$ attractor}
    \end{subfigure}
    \label{fig:LorenzData}
    \caption{Lorenz 1963 Attractor followed by its three two dimensional attractor.}
\end{figure}
In the following collection of figures, we provide the time series information of the Lorenz 1963 data.
\begin{figure}[H]
    \centering
    \begin{subfigure}{0.3\textwidth}
    \centering
        \includegraphics[width=\linewidth]{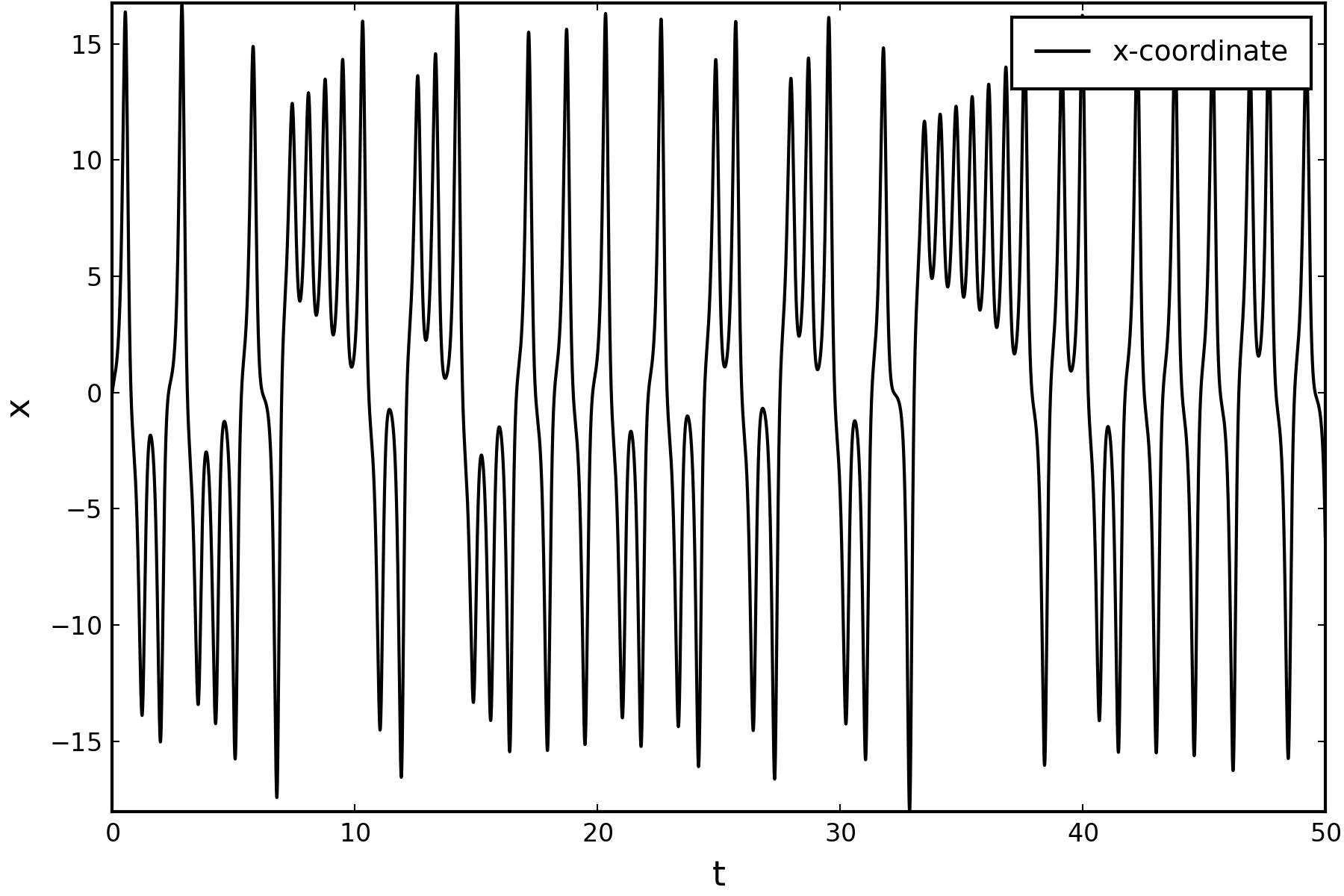}
        \caption{$x$ dimension}
    \end{subfigure}
    \begin{subfigure}{0.3\textwidth}
    \centering
        \includegraphics[width=\linewidth]{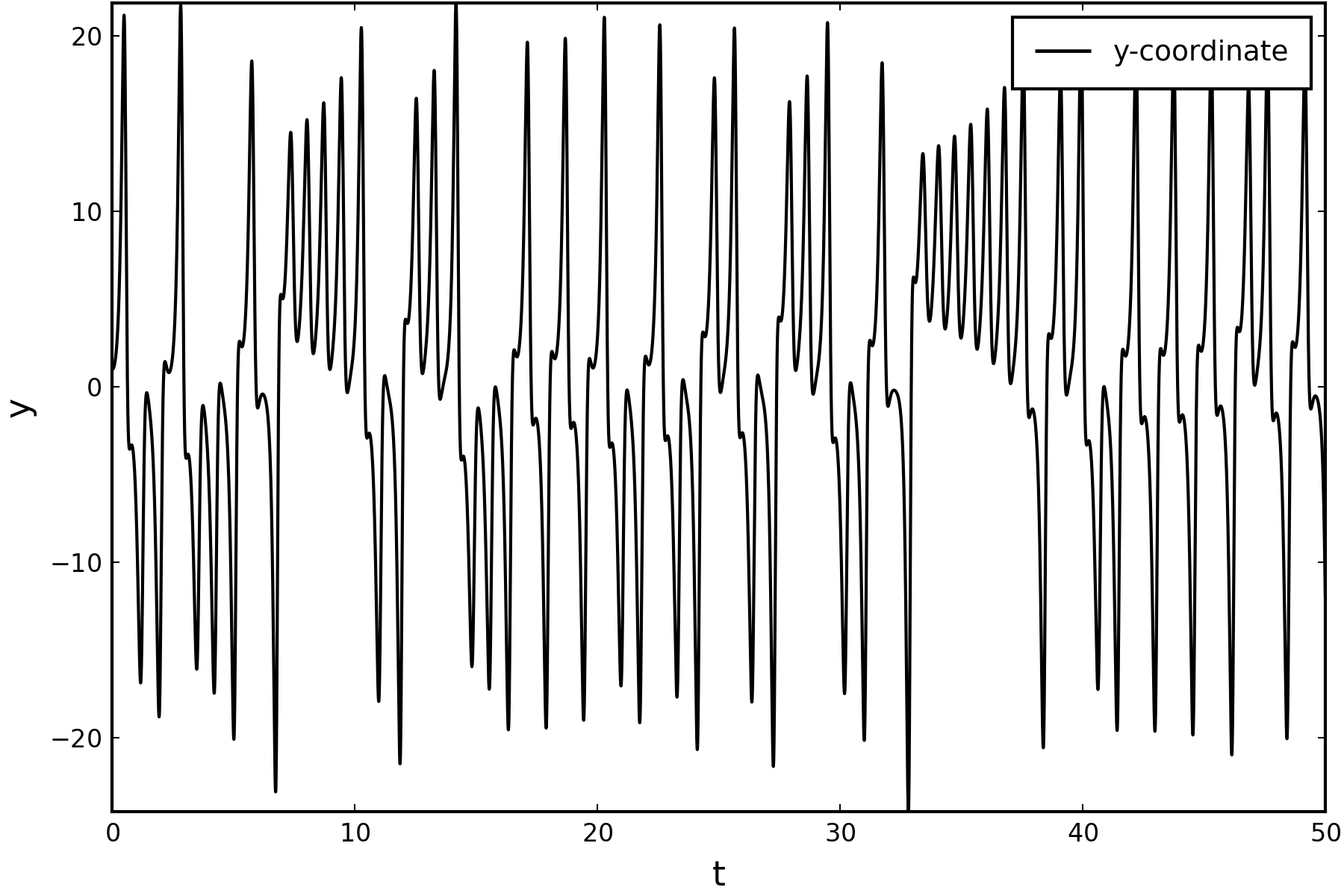}
        \caption{$y$ dimension}
    \end{subfigure}
    \begin{subfigure}{0.3\textwidth}
    \centering
        \includegraphics[width=\linewidth]{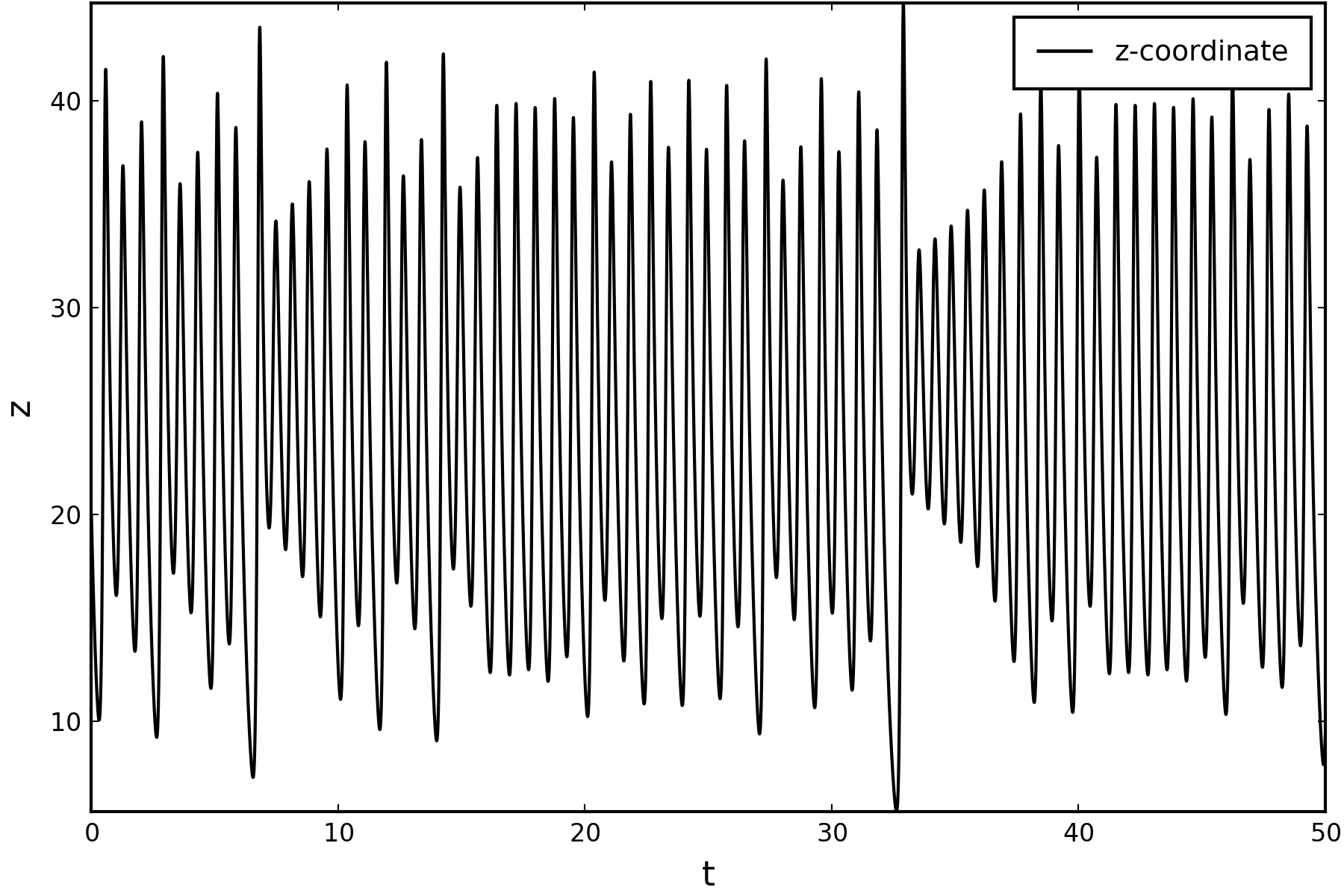}
        \caption{$z$ dimension}
    \end{subfigure}
    \caption{Lorenz 1963 time series all across three dimension.}
    \label{fig:LorenzTimeSeries}
\end{figure}
We now present the results of the spatial-temporal reconstruction and their error plots of the Lorenz 1963 dynamical system through the kernel-extended DMD by both kernels. The results shown below correspond to the snapshot $16$\textsuperscript{th} out of the total $30$ snapshots available.
\begin{figure}[H]
    \centering
    \begin{subfigure}{0.3\textwidth}
    \centering
    \includegraphics[width=\linewidth]{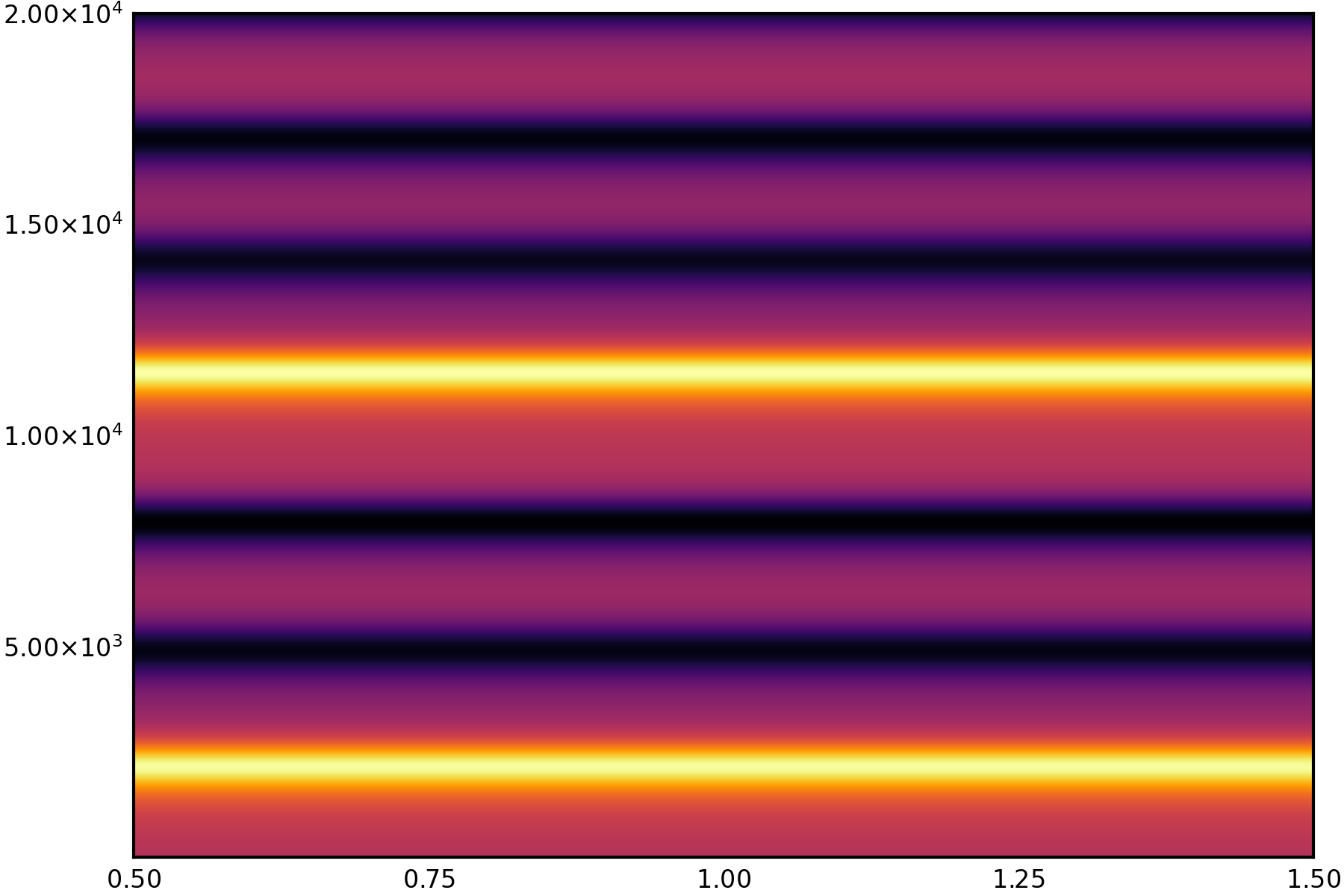}
    \caption{Actual}
    \end{subfigure}
    \begin{subfigure}{0.3\textwidth}
    \centering
    \includegraphics[width=\linewidth]{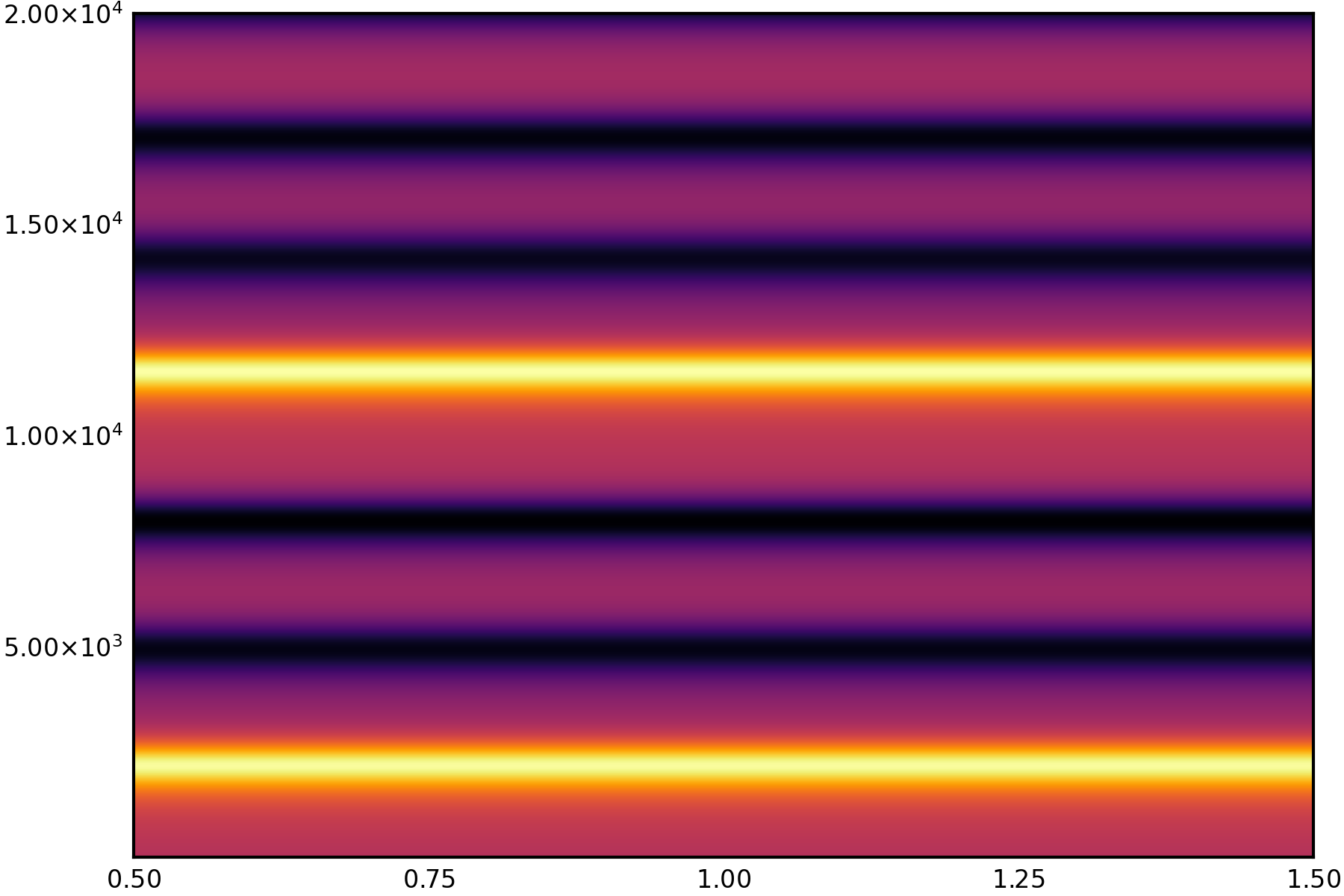}
    \caption{Reconstruction by $\Lap$}
    \end{subfigure}
    \begin{subfigure}{0.3\textwidth}
    \centering
    \includegraphics[width=\linewidth]{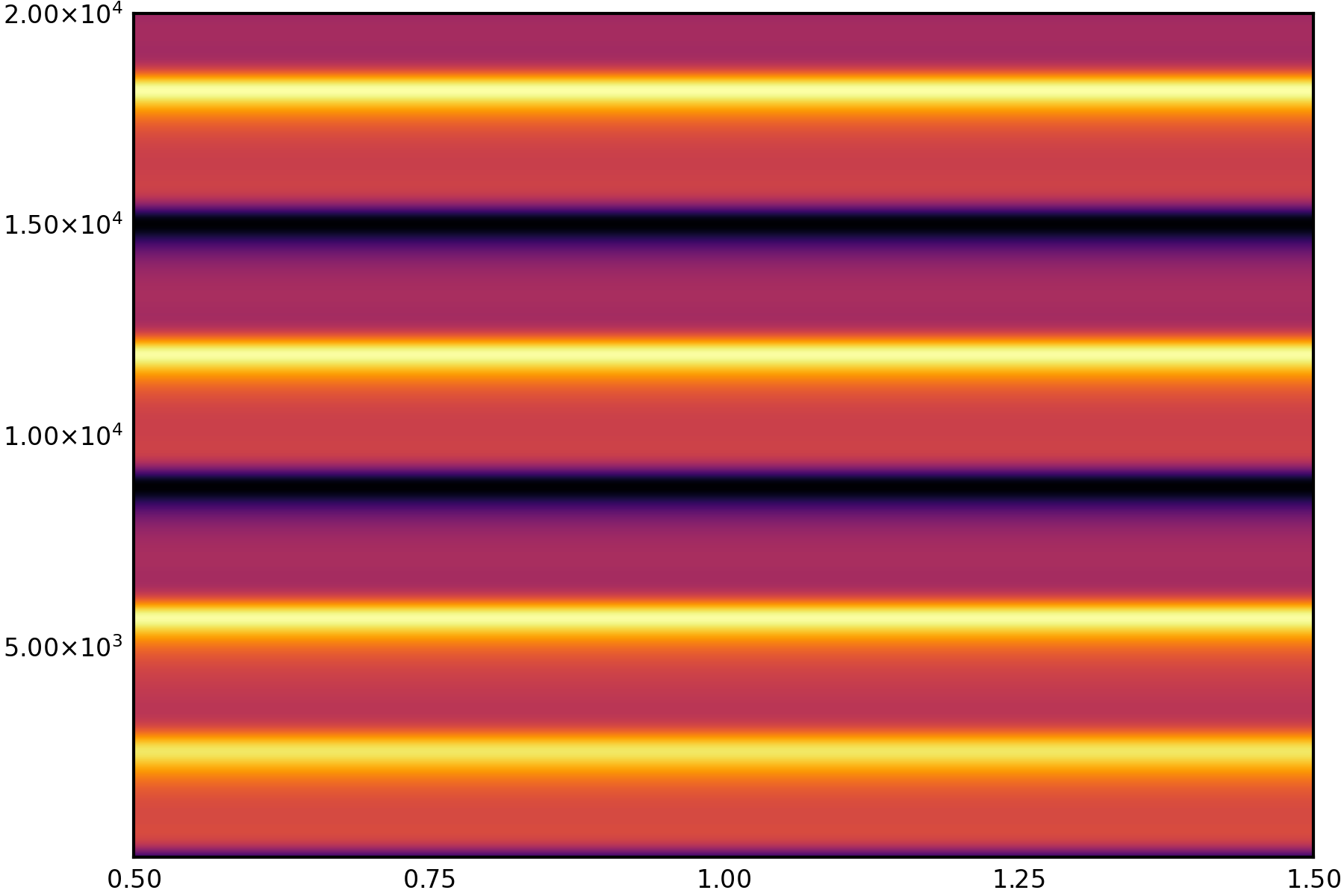}
    \caption{Reconstruction by $\GRBF$}
    \end{subfigure}
    \label{fig:LorenzeDMDExperiment}
    \caption{Spatial reconstruction result for the Lorenz 1963 dynamical system.}
\end{figure}
\begin{figure}[H]
    \centering
    \begin{subfigure}{0.3\textwidth}
        \includegraphics[width=\linewidth]{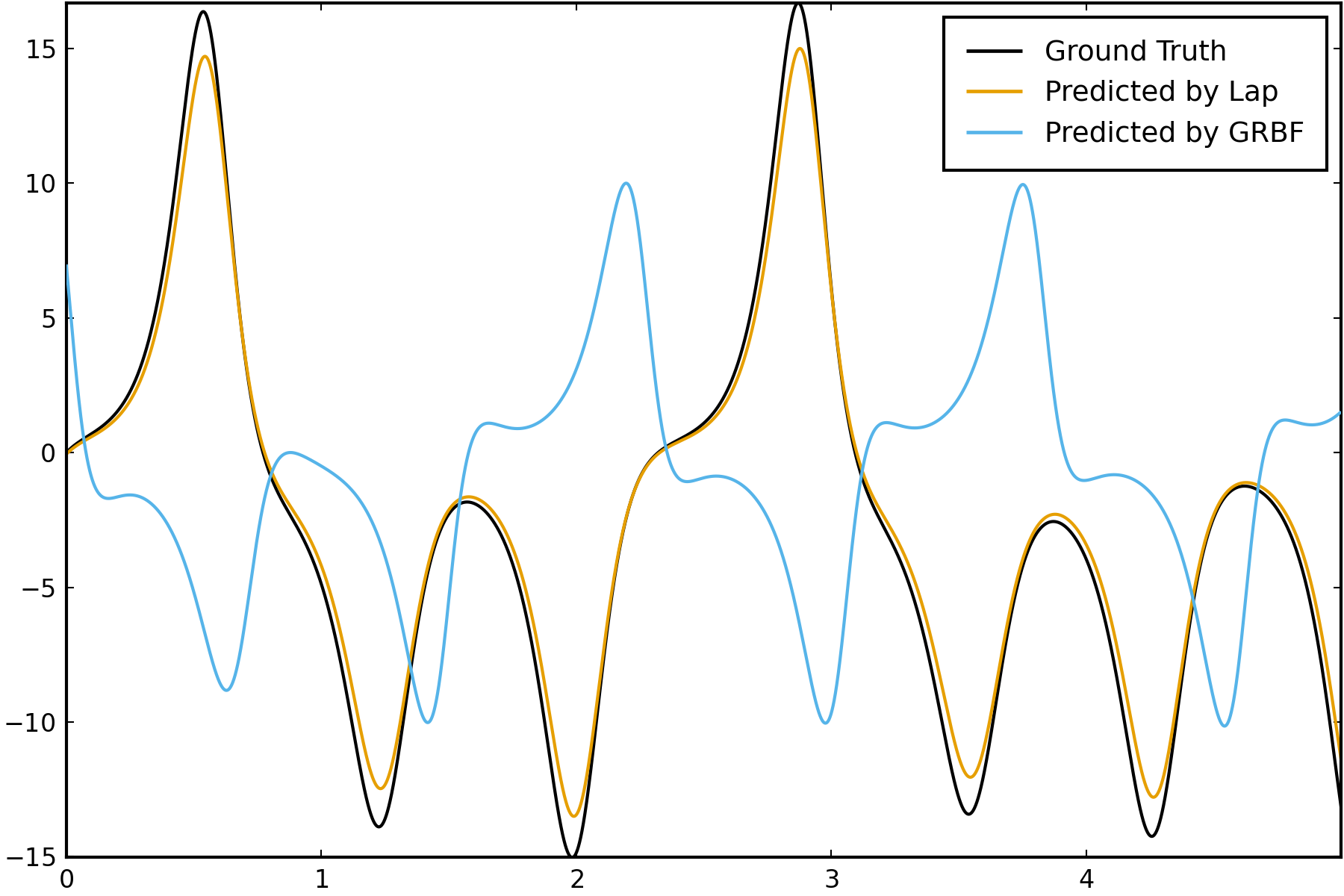}
    \caption{Spatial values}
    \end{subfigure}
    \hfil
    \begin{subfigure}{0.3\textwidth}
        \includegraphics[width=\linewidth]{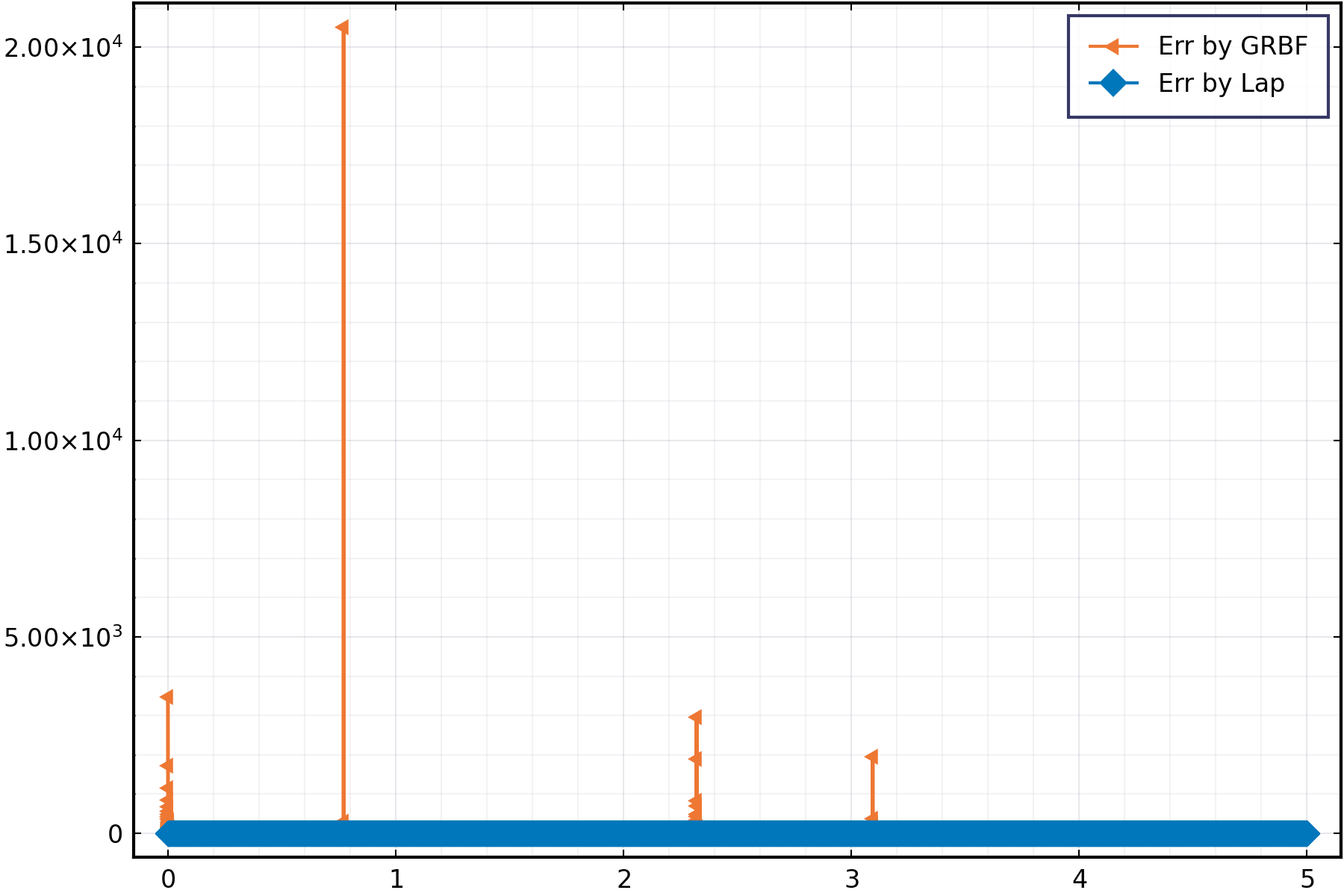}
    \caption{EWE plot}
    \end{subfigure}
    \label{fig:LorenzEDMDComparison}
    \caption{Spatial reconstruction result for the Lorenz 1963 dynamical system.}
\end{figure}
\subsubsection{Experiment 6}
For the following experiment, we will provide first the initial details of the experiment as we did for the Lorenz 1963 dynamical system.
\begin{figure}[H]
    \centering
    \begin{subfigure}{0.3\textwidth}
    \centering
        \includegraphics[width=.7\linewidth]{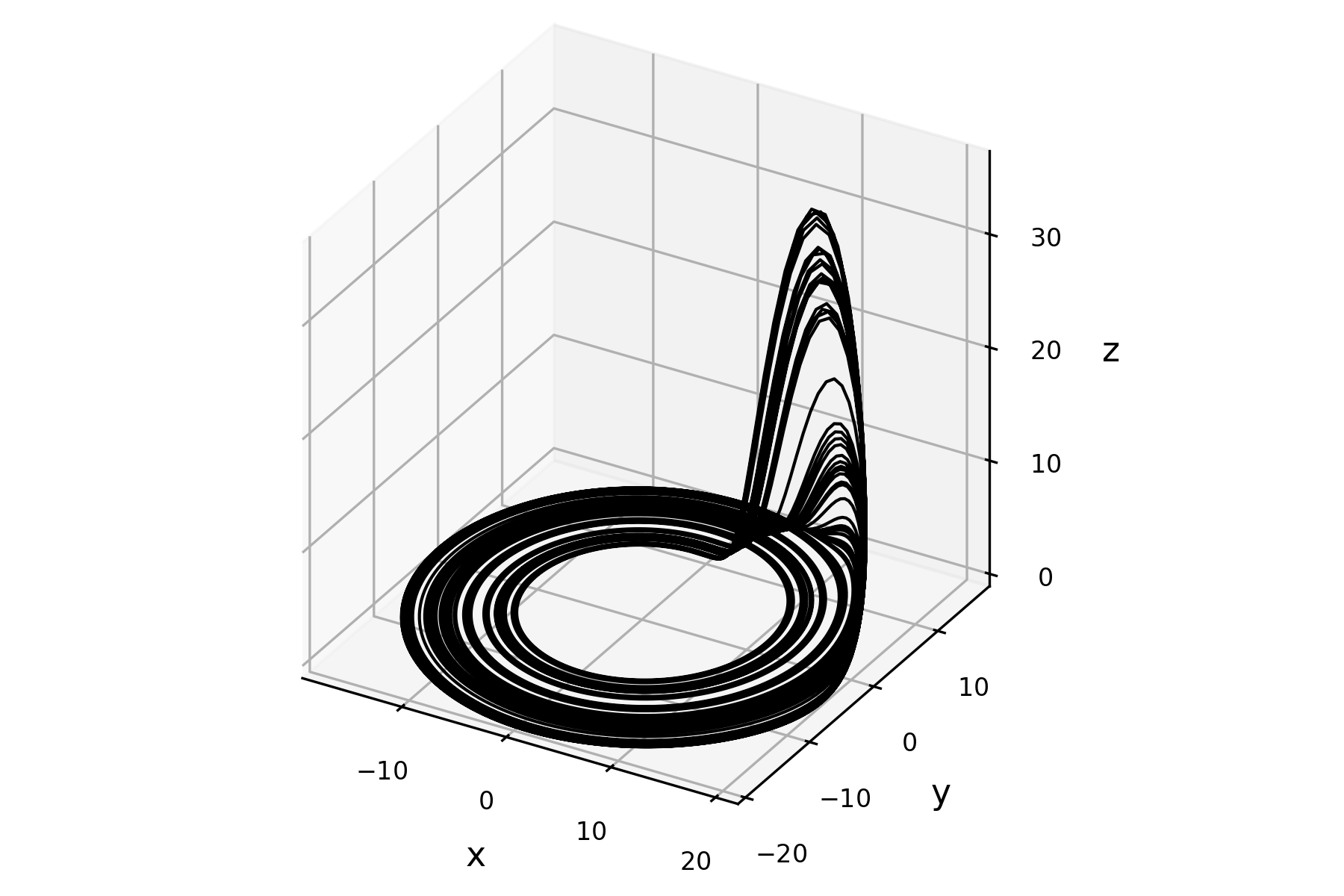}
    \caption{Rössler attractor}
    \end{subfigure}
    \begin{subfigure}{0.3\textwidth}\centering
        \includegraphics[width=.7\linewidth]{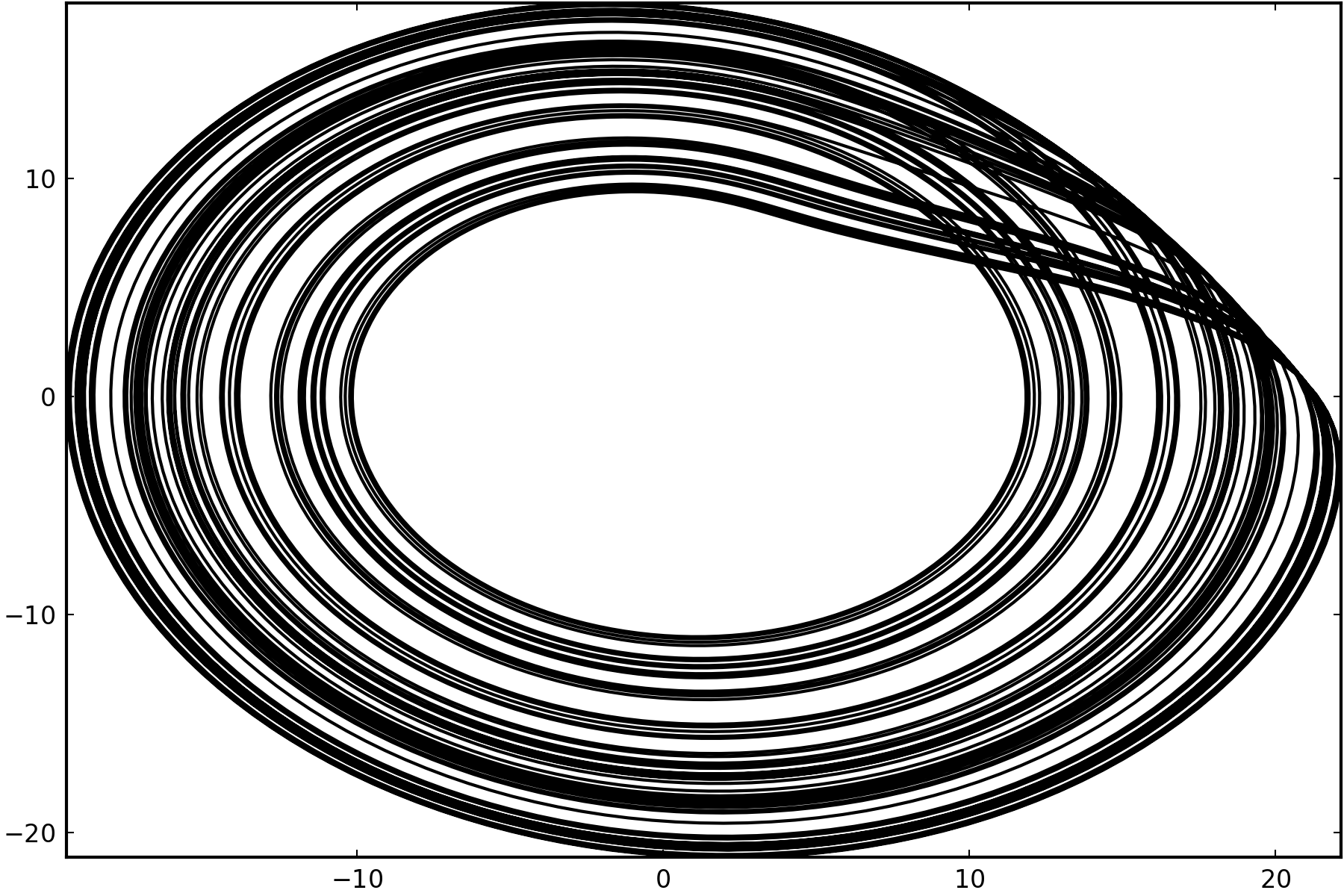}
    \caption{2D $xy$ attractor}
    \end{subfigure}
     \begin{subfigure}{0.3\textwidth}\centering
        \includegraphics[width=.7\linewidth]{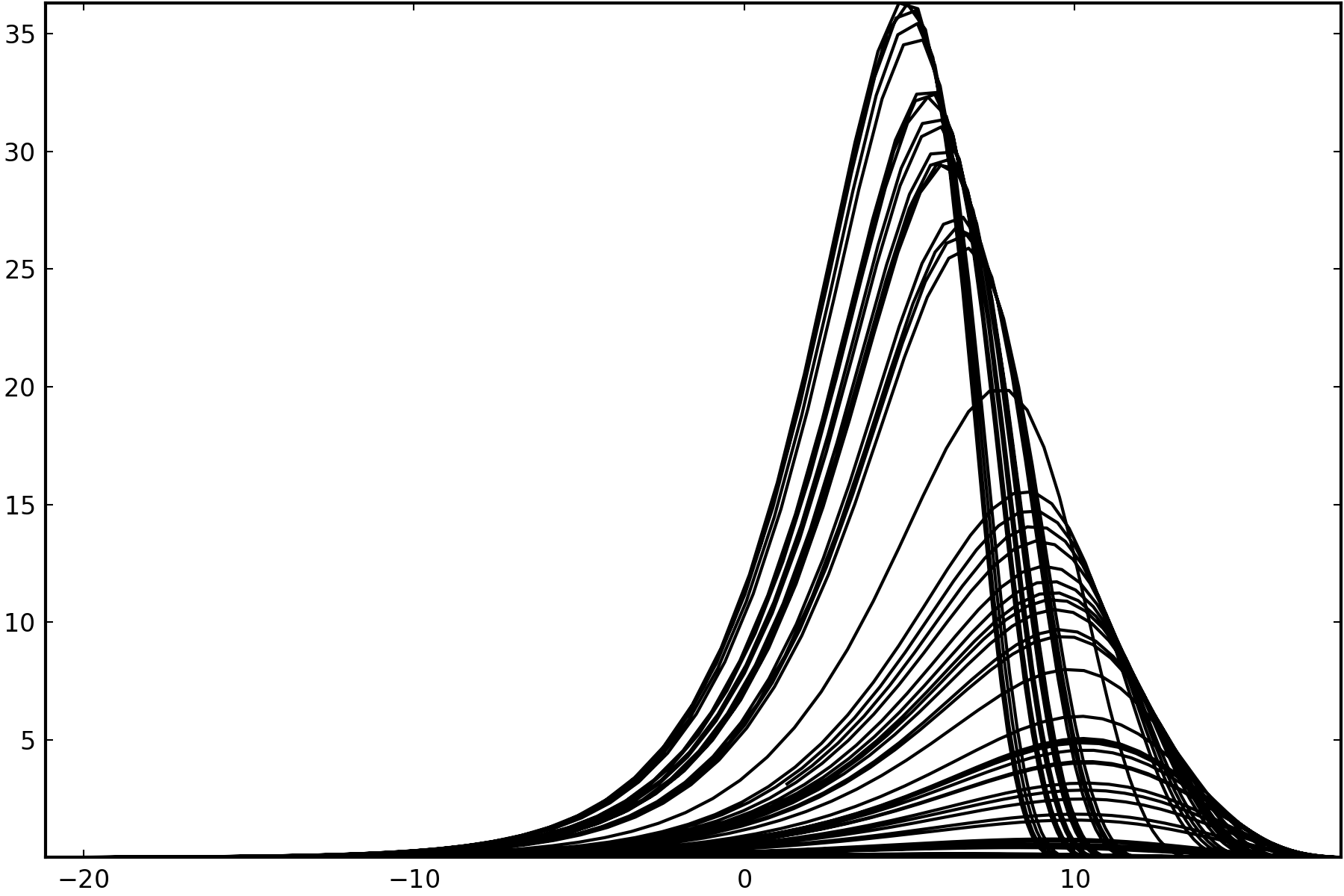}
    \caption{2D $yz$ attractor}
    \end{subfigure}
    \begin{subfigure}{0.3\textwidth}\centering
        \includegraphics[width=.7\linewidth]{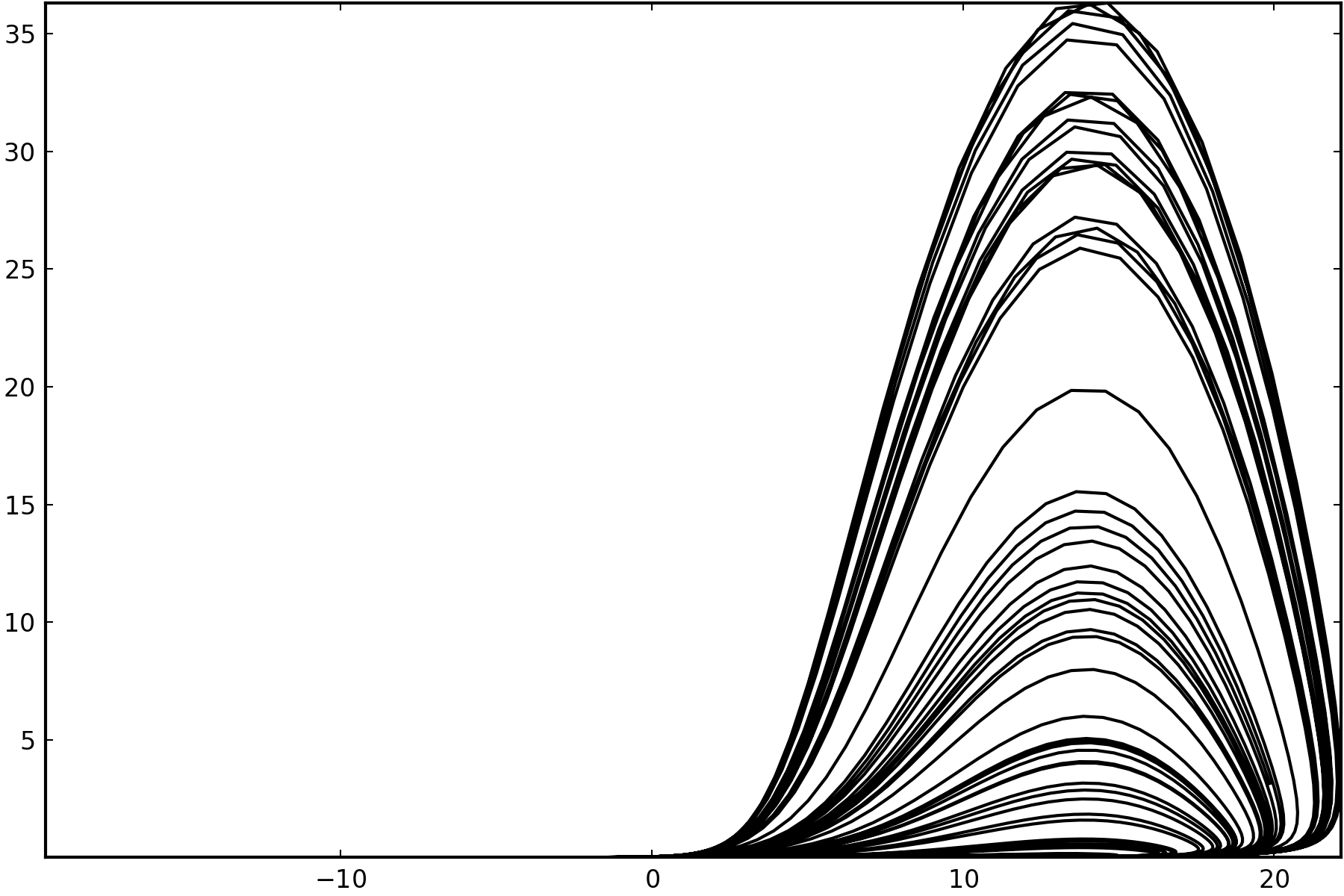}
    \caption{2D $xz$ attractor}
    \end{subfigure}
    \caption{Rössler Attractor followed by its three two dimensional attractor.}
    \label{fig:RosslerActual}
\end{figure}
Now, we provide the time series evolution of the Rössler attractor with respect to time.
\begin{figure}[H]
    \centering
    \begin{subfigure}{0.3\textwidth}\centering
        \includegraphics[width=1\linewidth]{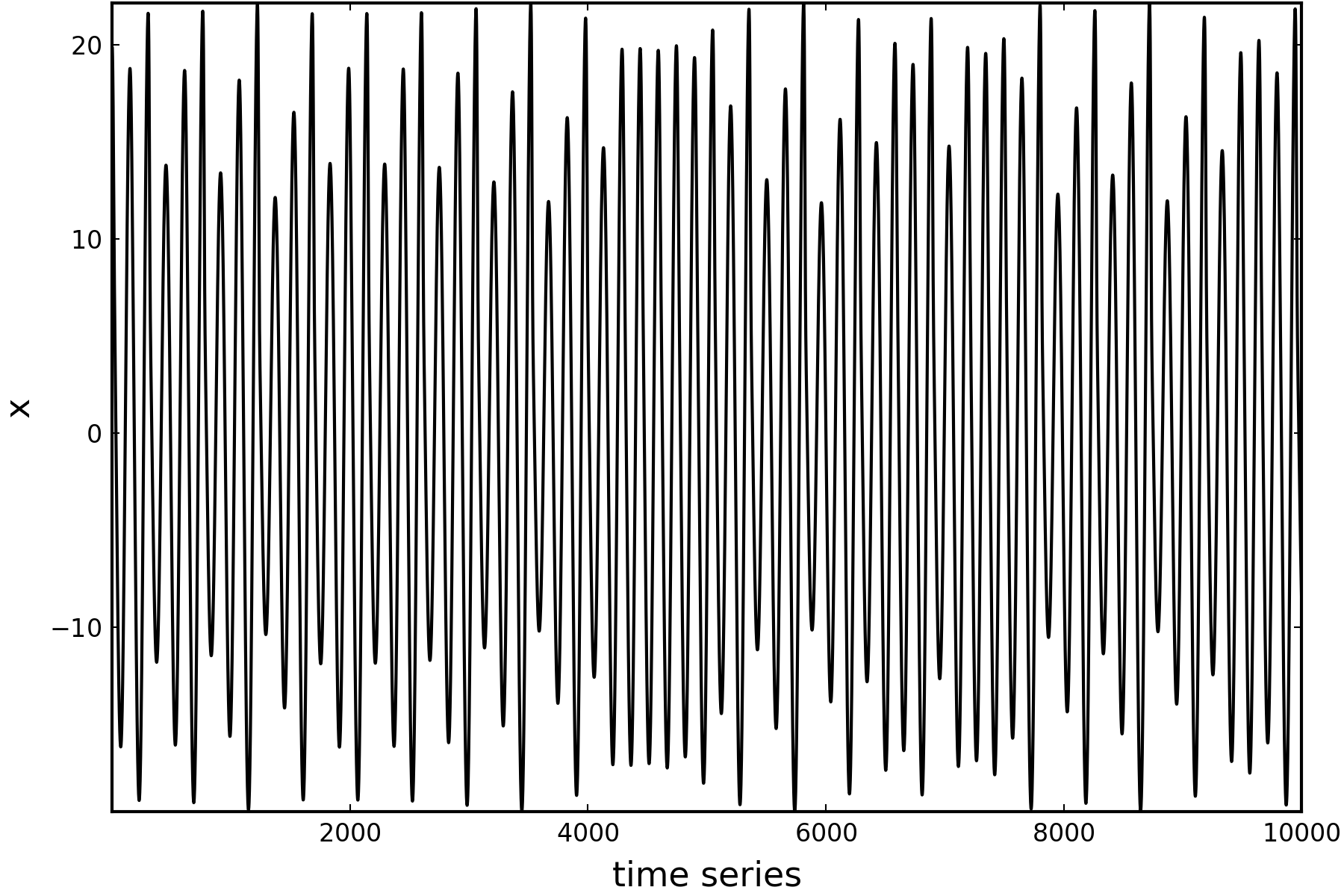}
    \caption{$x$ dimension}
    \end{subfigure}
    \begin{subfigure}{0.3\textwidth}\centering
        \includegraphics[width=1\linewidth]{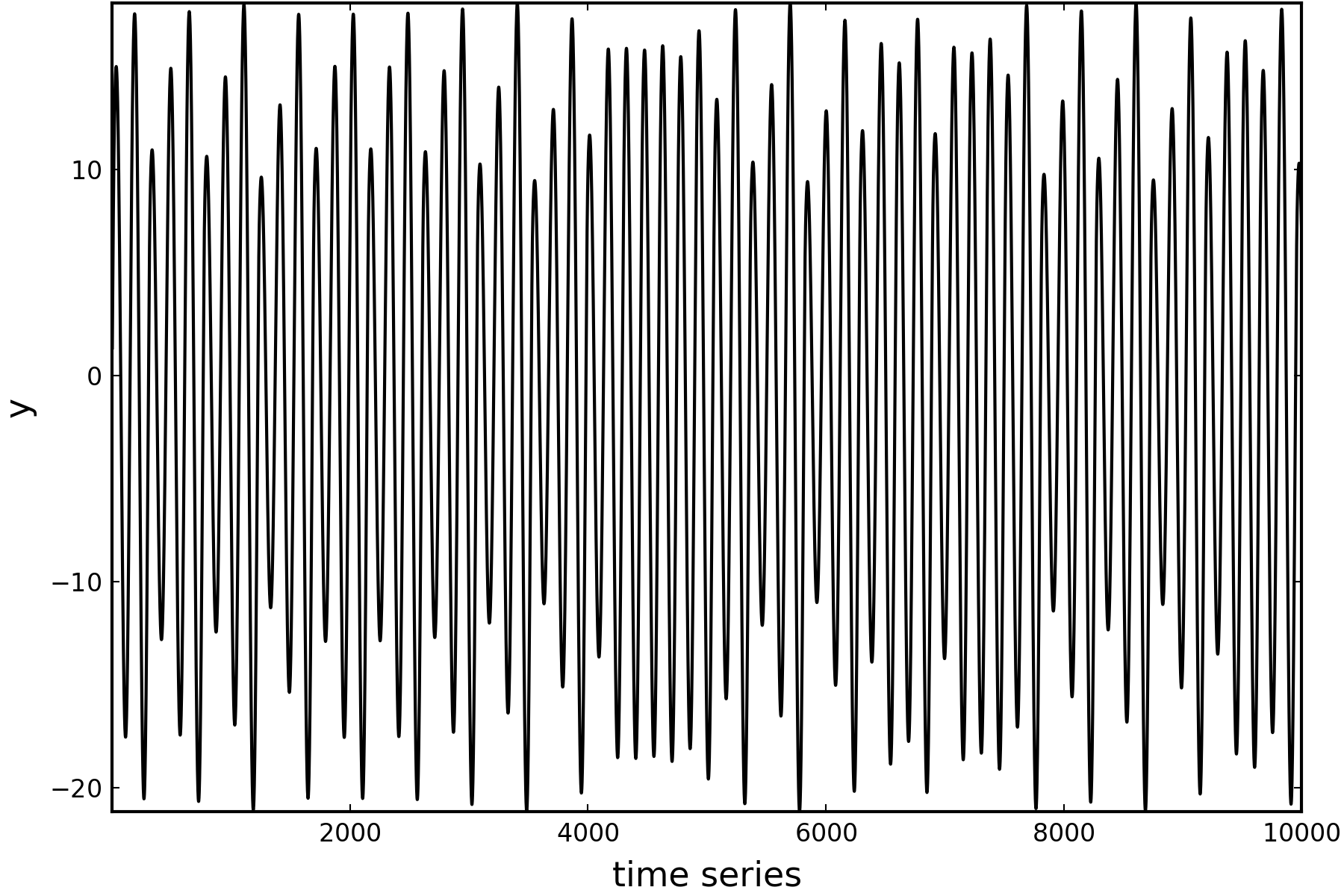}
    \caption{$y$ dimension}
    \end{subfigure}
    \begin{subfigure}{0.3\textwidth}\centering
        \includegraphics[width=1\linewidth]{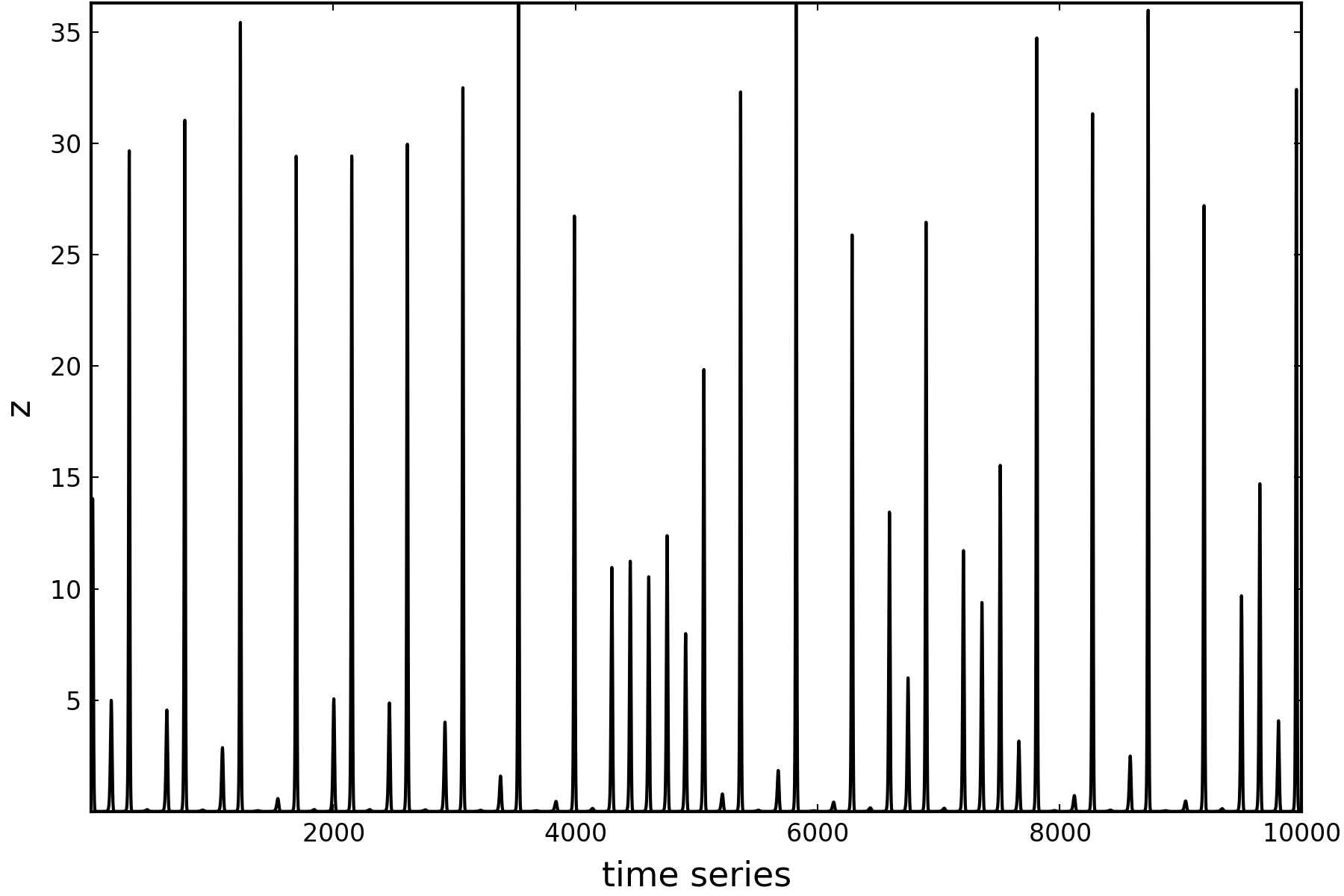}
    \caption{$z$ dimension}
    \end{subfigure}
\caption{Rössler Attractor time series across all three dimensions.}
    \label{fig:RosslerActualdatatimeseries}
\end{figure}
Now, we are ready to provide the kernel-extended DMD result yielded by the choice of both kernels.
\begin{figure}[H]
    \centering
    \begin{subfigure}{0.3\textwidth}\centering
        \includegraphics[width=1\linewidth]{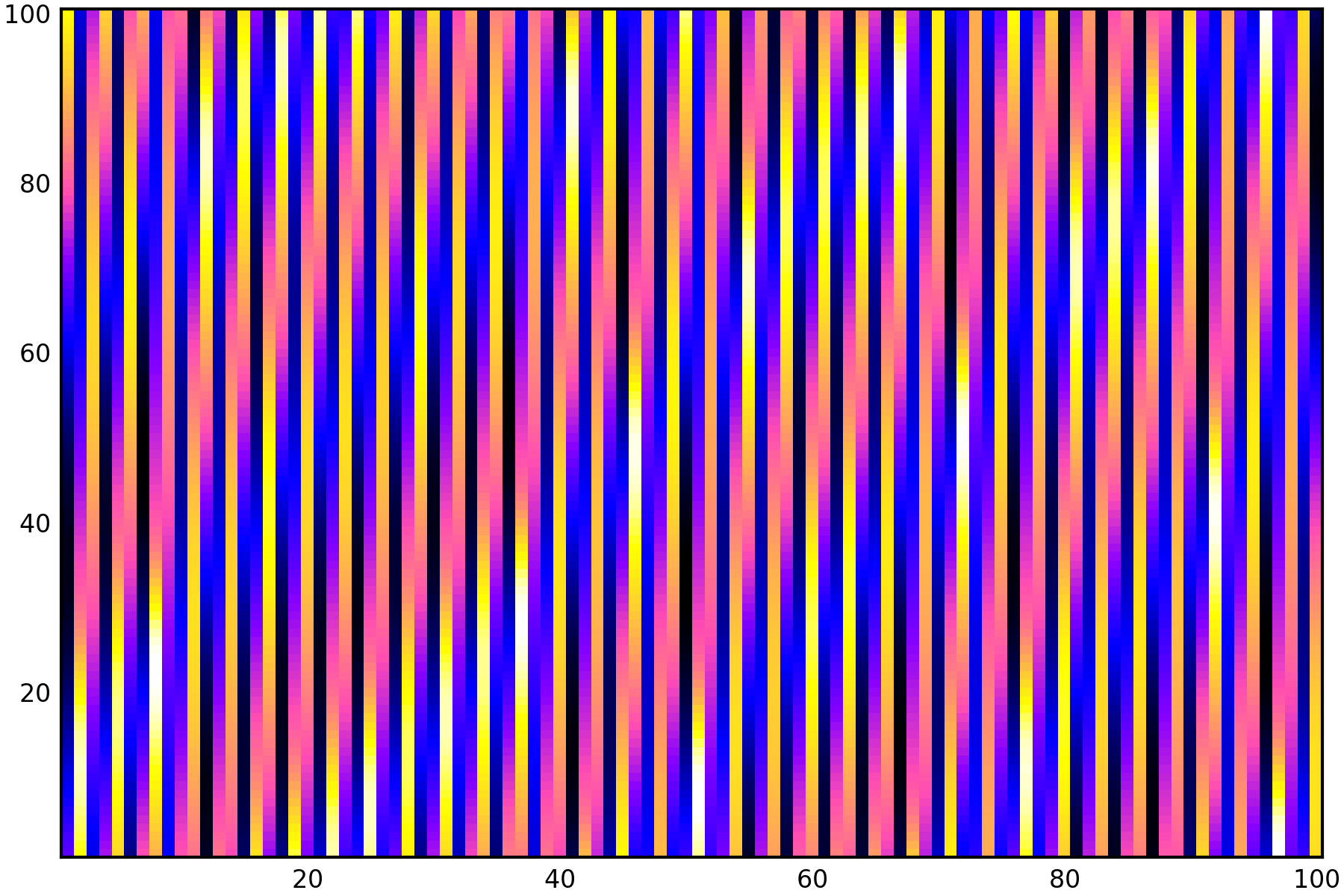}
    \caption{Actual}
    \end{subfigure}
    \begin{subfigure}{0.3\textwidth}\centering
        \includegraphics[width=1\linewidth]{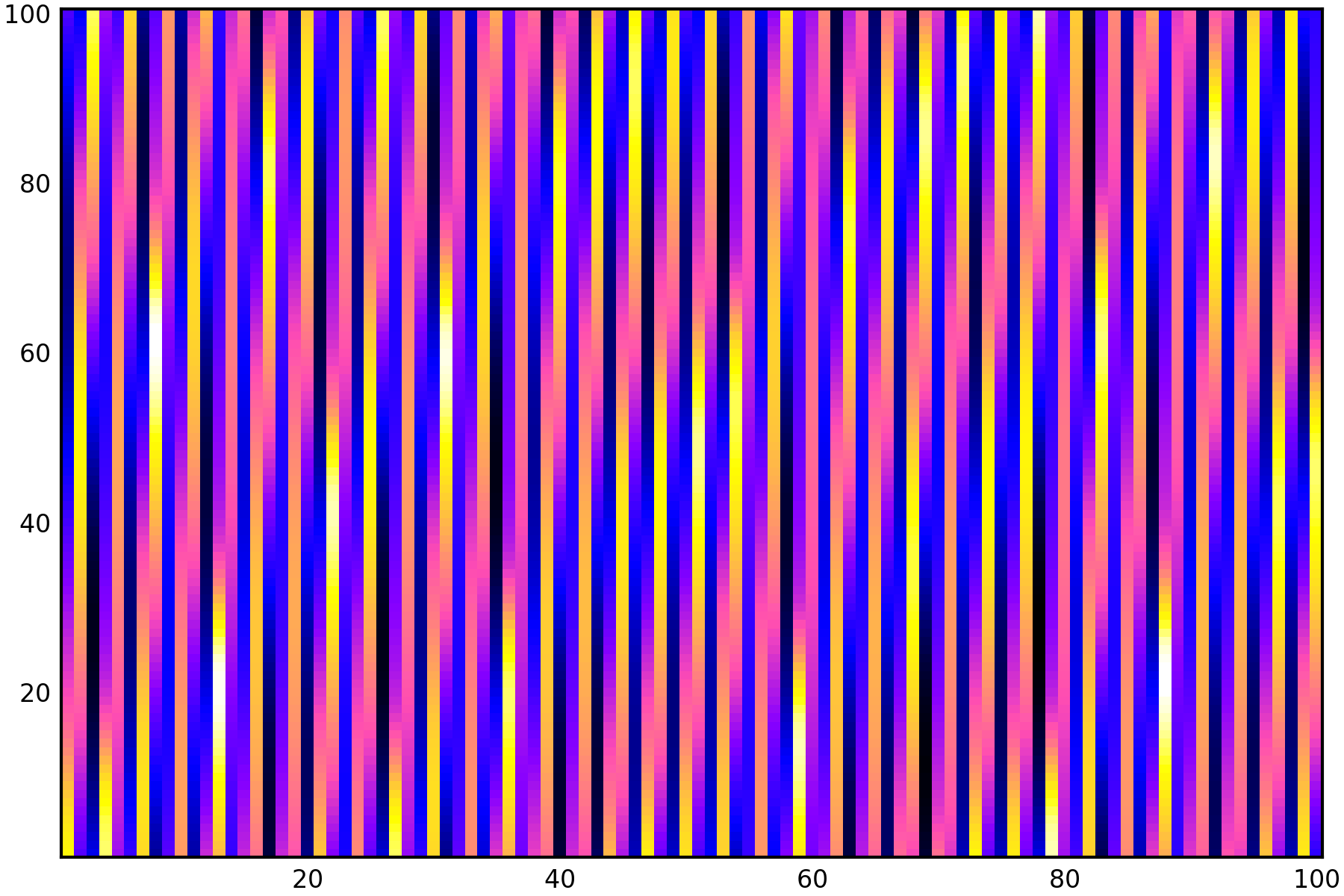}
    \caption{Reconstruction by $\Lap$}
    \end{subfigure}
    \begin{subfigure}{0.3\textwidth}\centering
        \includegraphics[width=1\linewidth]{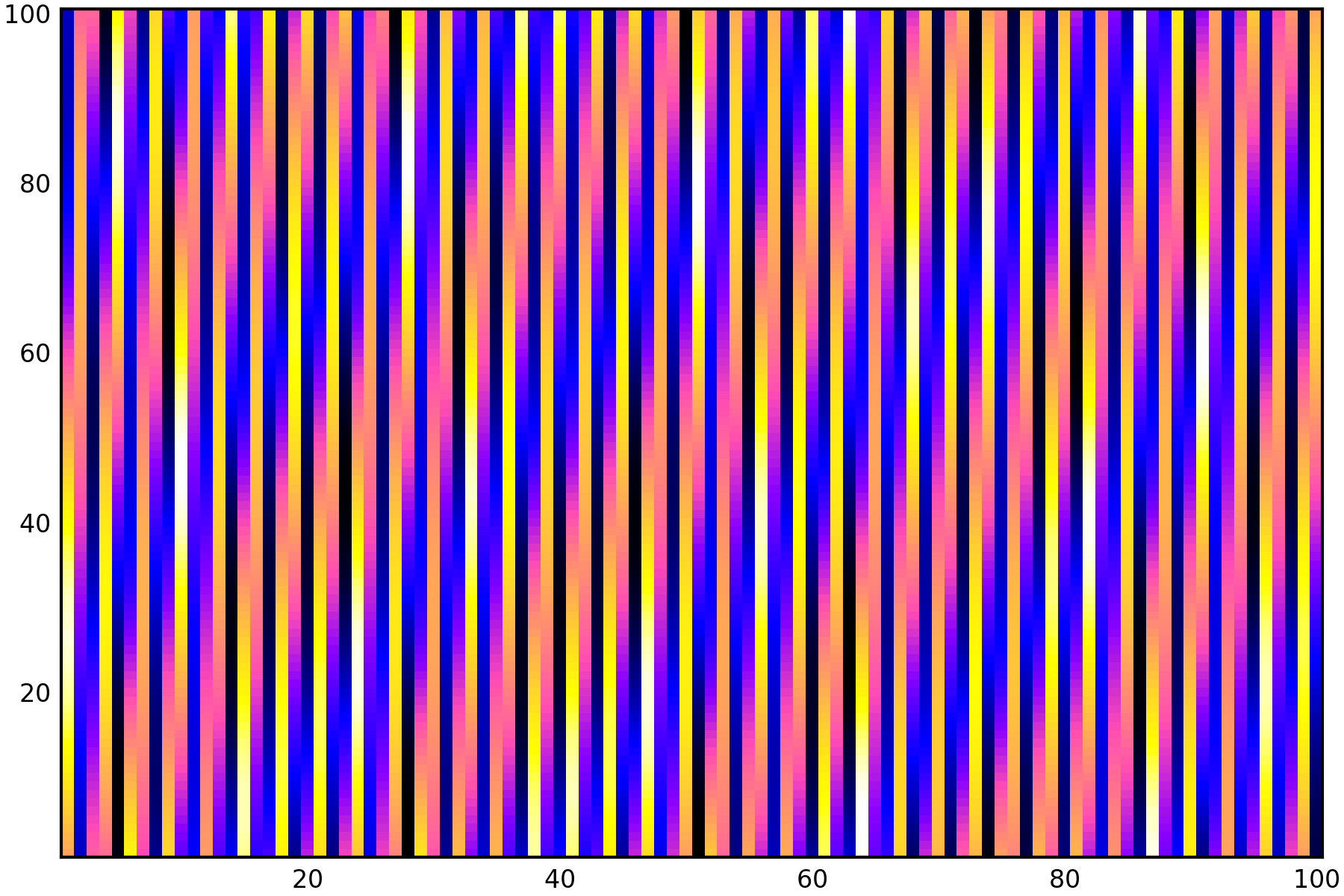}
    \caption{Reconstruction by $\GRBF$}
    \end{subfigure}
    \caption{Spatial values reconstruction for Rössler Attractor.}
    \label{fig:ReconstructionResultRossler}
\end{figure}
The following is now the spatial plots and error plots as well.
\begin{figure}[H]
    \centering
    \begin{subfigure}{0.35\textwidth}\centering
         \includegraphics[width=\linewidth]{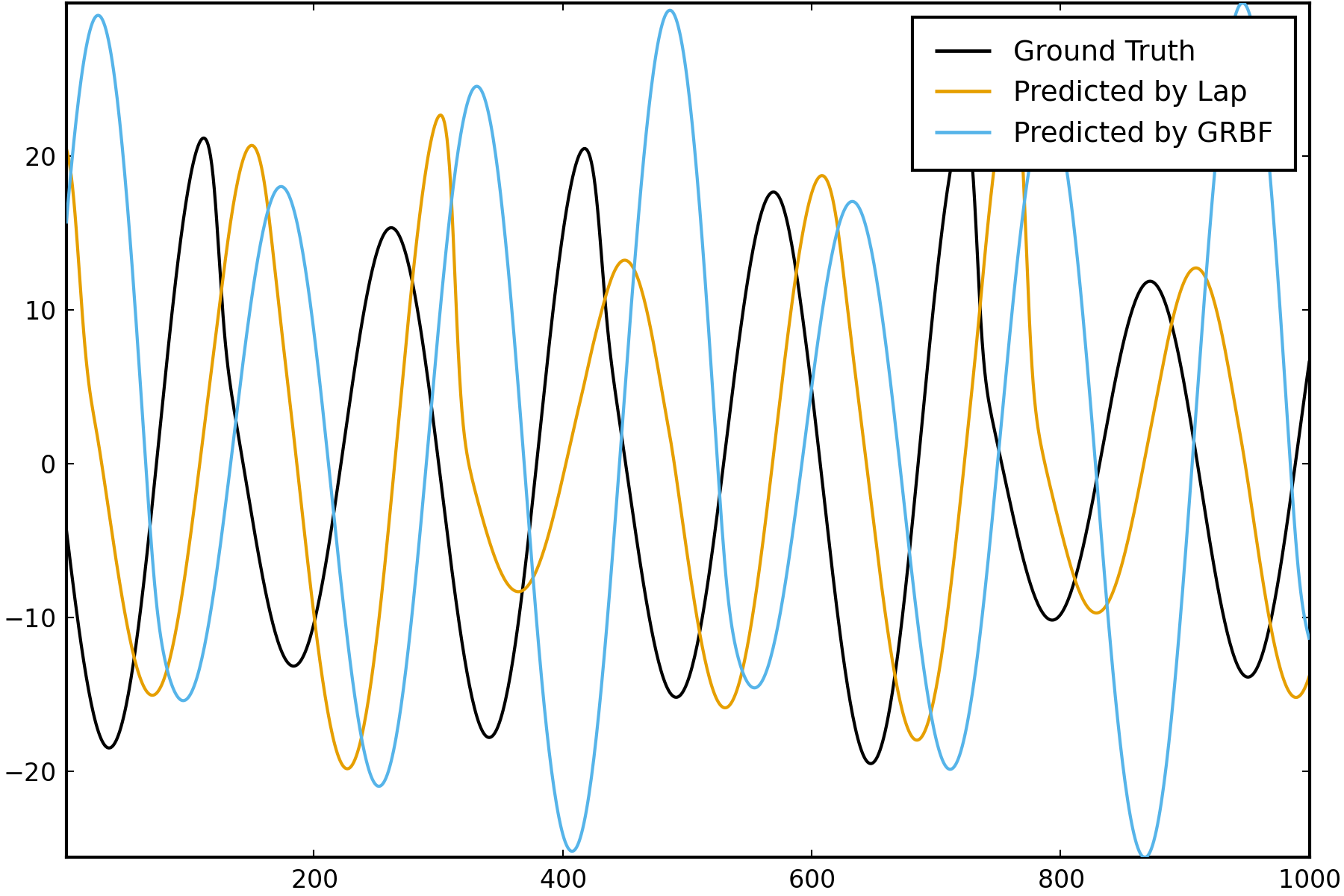}
    \caption{Spatial values}
    \end{subfigure}
    \hfil
   \begin{subfigure}{0.35\textwidth}\centering
         \includegraphics[width=\linewidth]{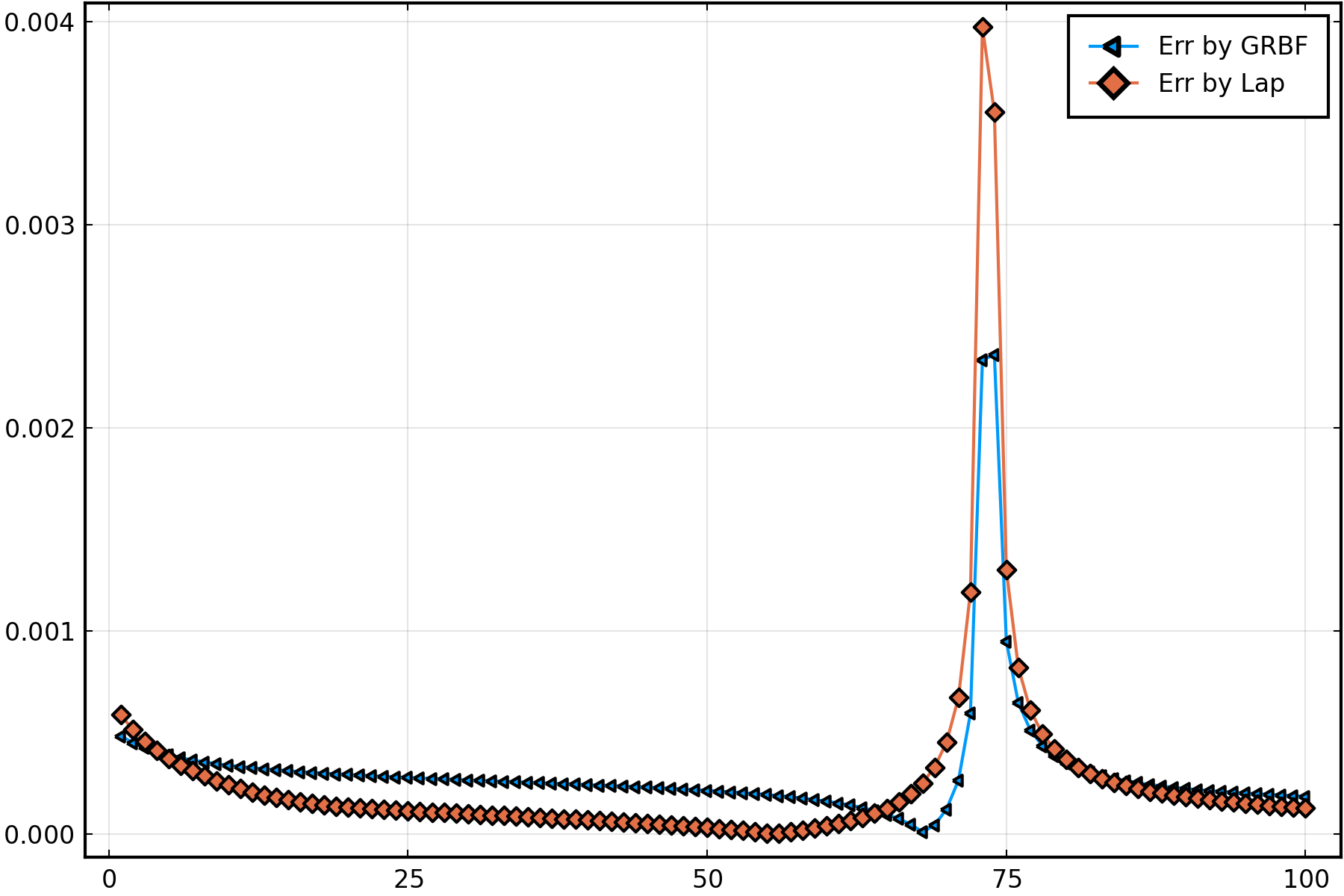}
    \caption{EWE plot}
    \end{subfigure}
    \caption{Spatial values reconstruction for Rössler Attractor followed by error plot.}
    \label{fig:SpatialValuesRossler}
\end{figure}
\subsubsection{Experiment 7}
The following data are from the US Department of Commerce of the National Oceanic and Atmospheric Administration (NOAA). We are interested in the reconstruction of the mean annual temperature anomalies for the year $2023$.  
\begin{figure}[H]
    \centering
       \includegraphics[width=\linewidth]{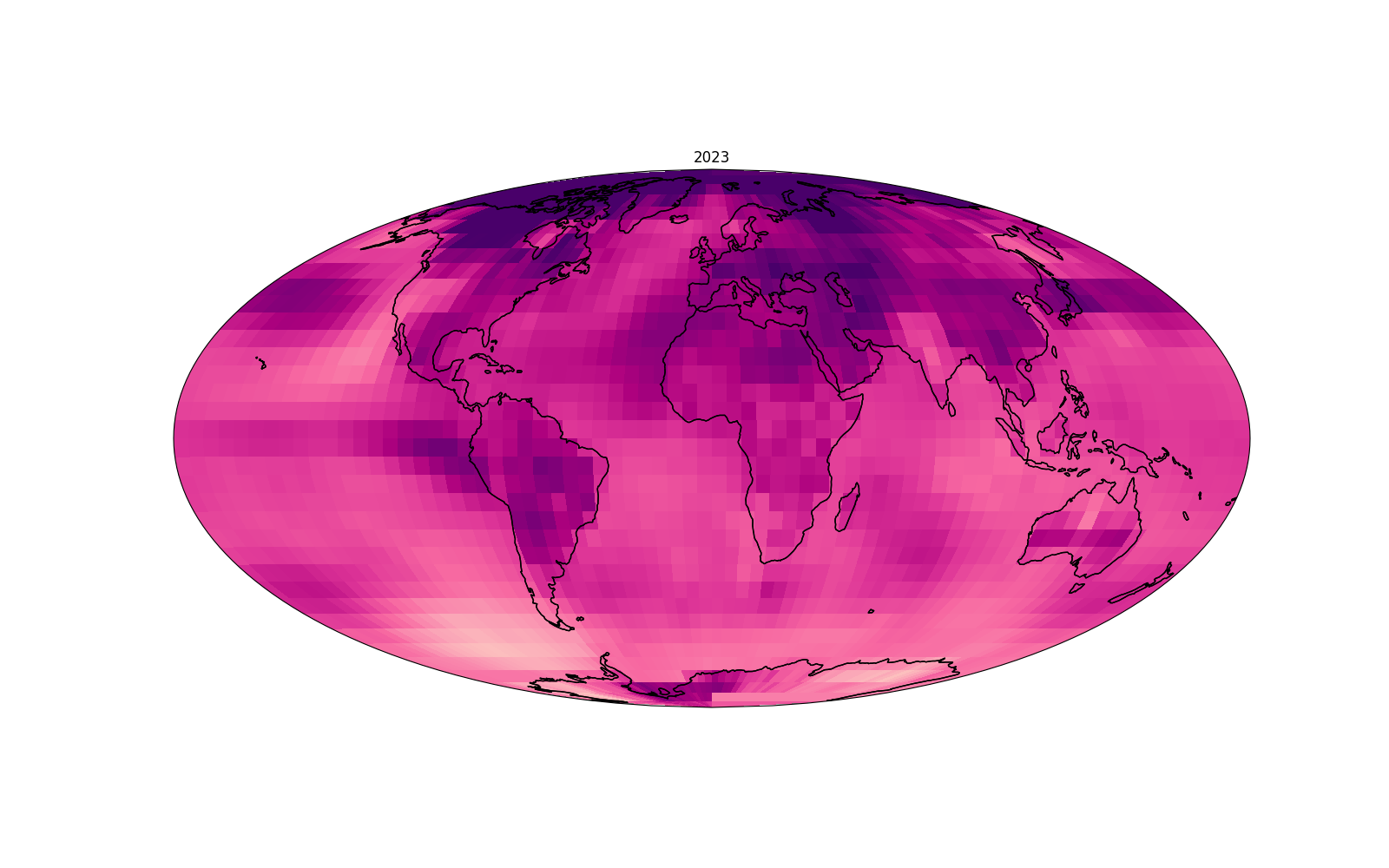} 
       \caption{NOAA Global Sea Surface Temperature Anomaly for year 2023.}
    \label{fig:Weatherdatasnaps1-10}
\end{figure}
The following sets of results present the reconstructed results corresponding to the $9$\textsuperscript{th} snapshot which corresponds to the year $2023$.
\begin{figure}[H]
    \centering
    \begin{subfigure}{0.3\textwidth}\centering
    \includegraphics[width=\linewidth]{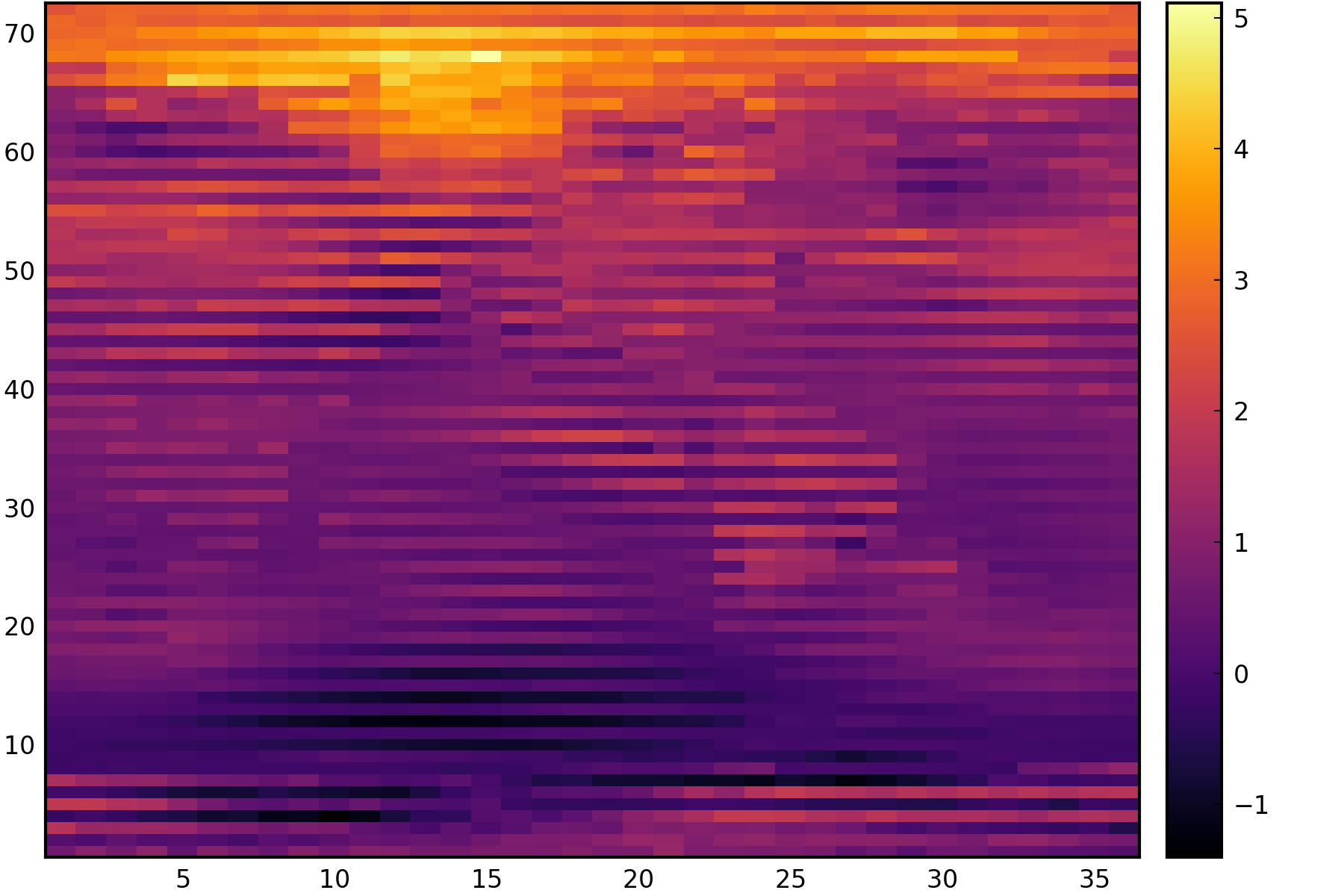}
    \caption{Actual}
    \end{subfigure}
    \begin{subfigure}{0.3\textwidth}\centering
    \includegraphics[width=\linewidth]{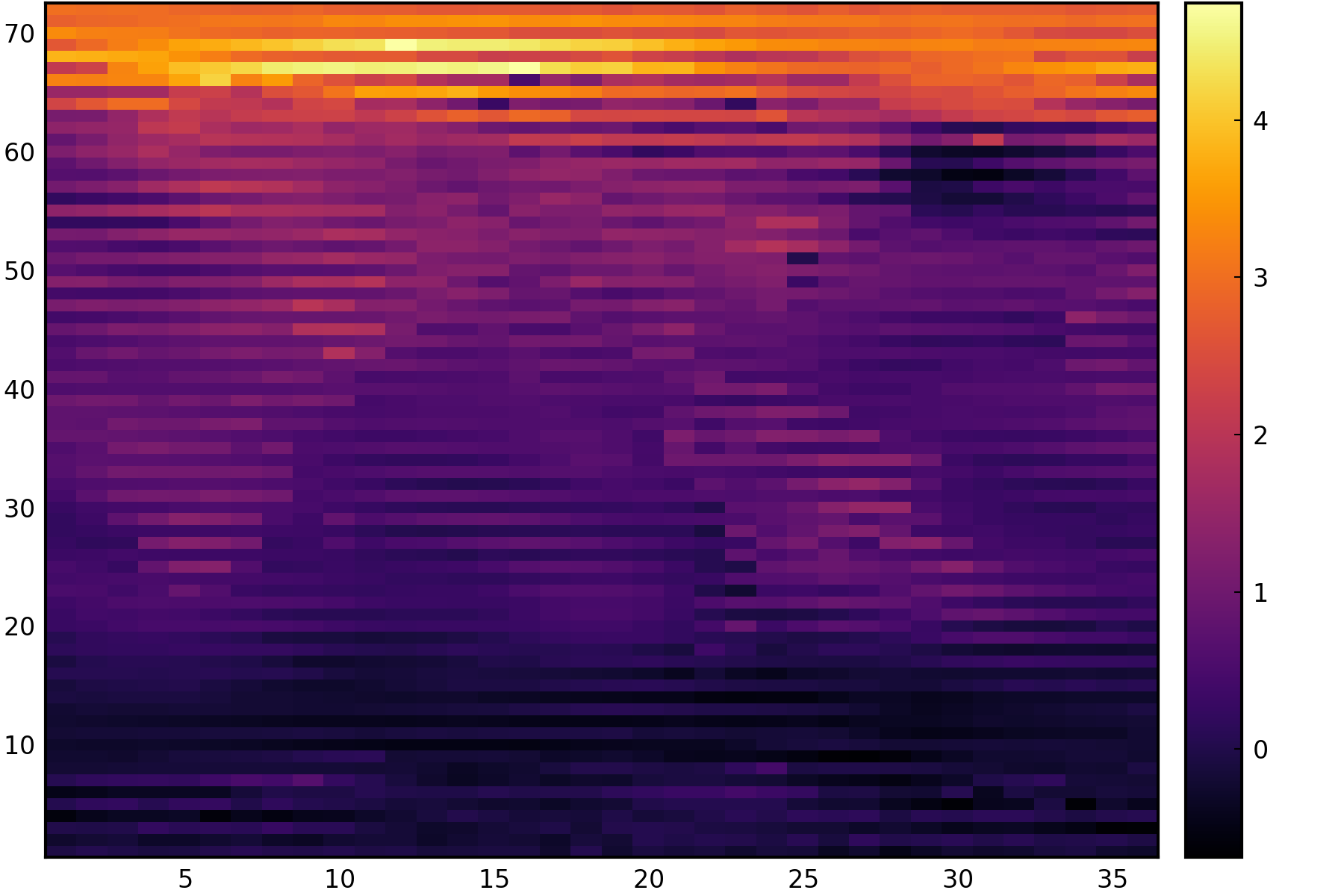}
    \caption{Reconstruction by $\Lap$}
    \end{subfigure}
    \begin{subfigure}{0.3\textwidth}\centering
    \includegraphics[width=\linewidth]{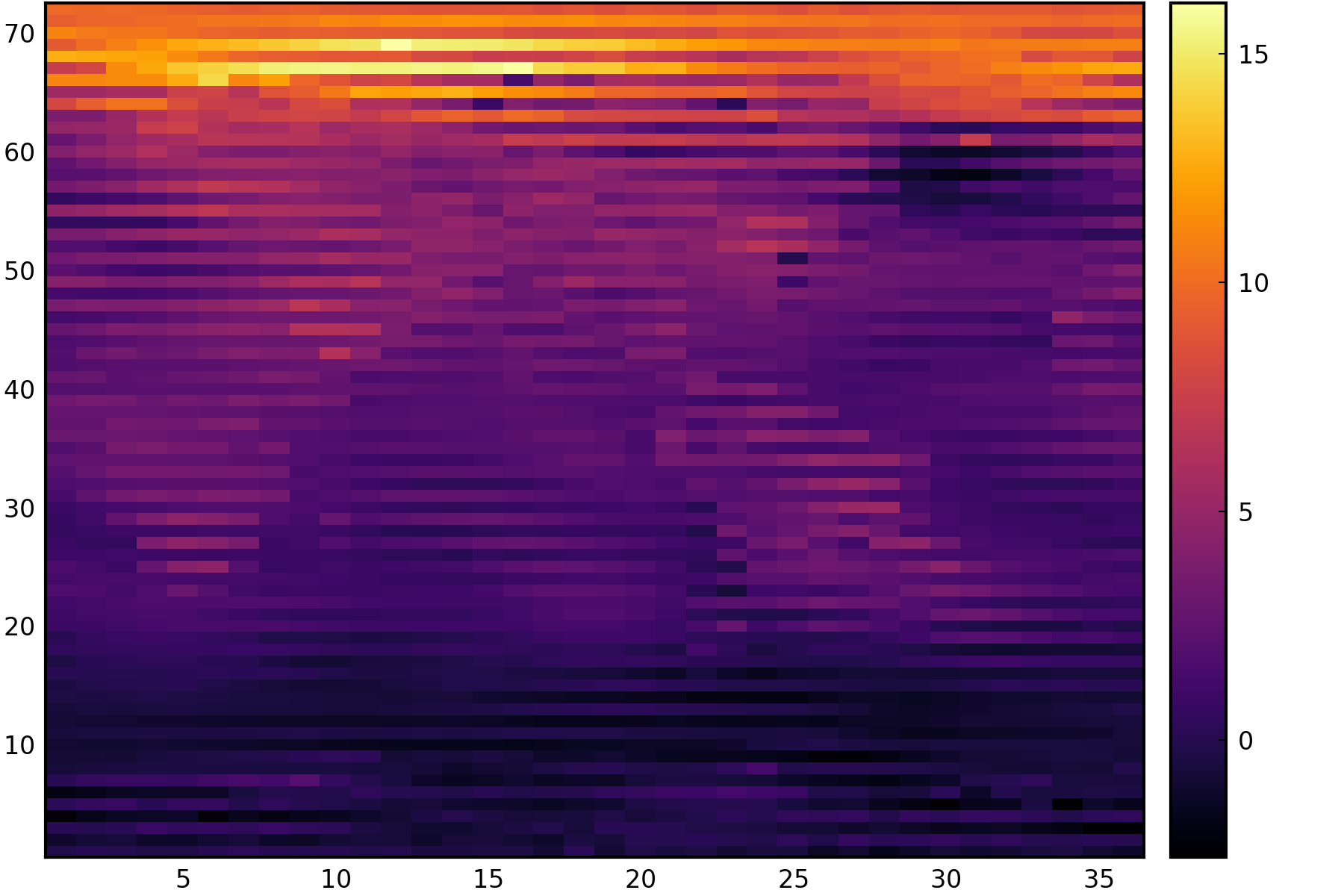}
    \caption{Reconstruction by $\GRBF$}
    \end{subfigure}
    \caption{Reconstruction of NOAA Global Sea Surface Temperature anomalies across latitude and longitude $(36\times72)$ for the year 2023.}
    \label{fig:WeatherSSTReco}
\end{figure}
\begin{figure}[H]
    \centering
    \begin{subfigure}{0.35\textwidth}\centering
        \includegraphics[width=\linewidth]{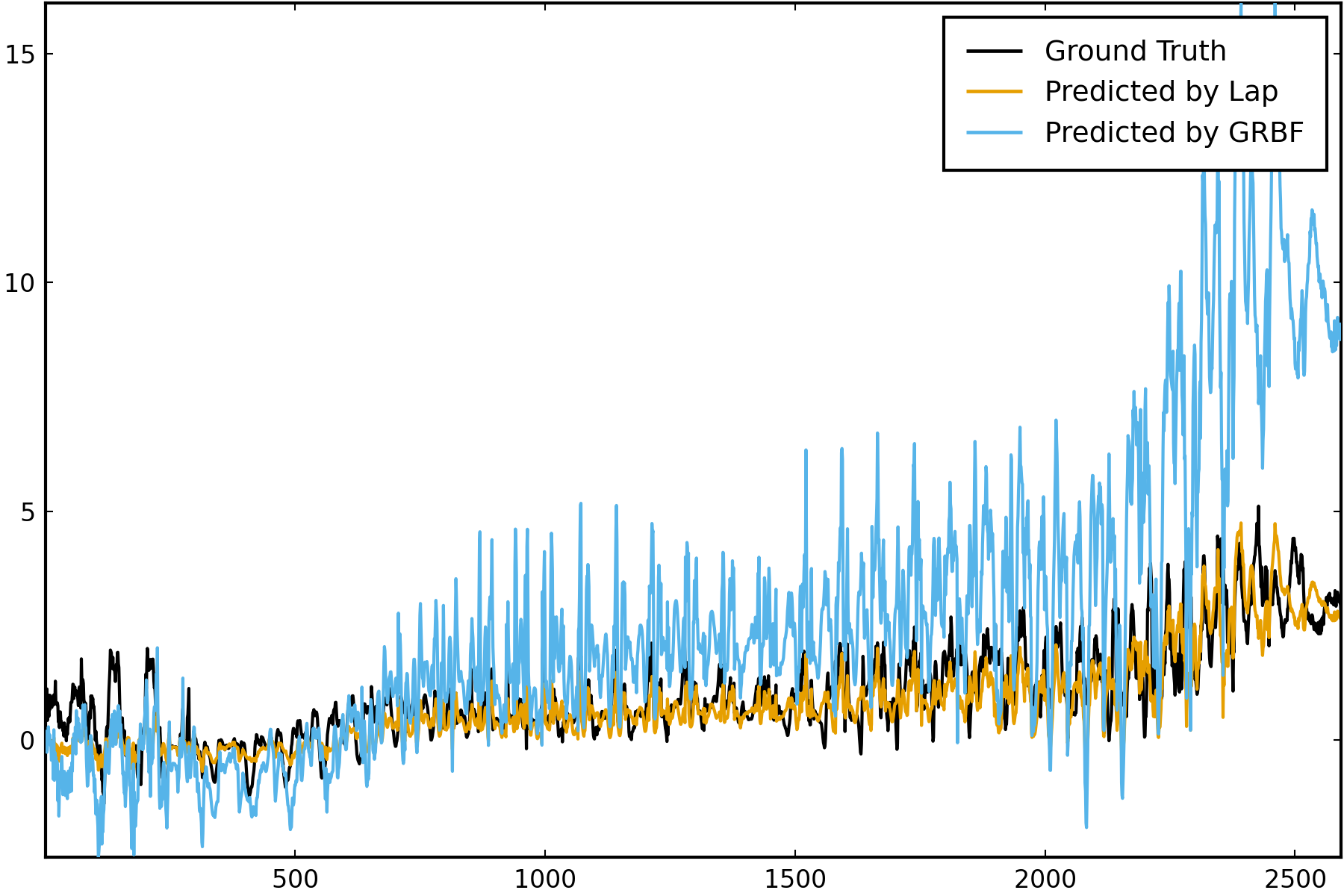}
        \caption{Temp. anomaly reconstruction}
    \end{subfigure}
    \hfil
    \begin{subfigure}{0.35\textwidth}\centering
        \includegraphics[width=\linewidth]{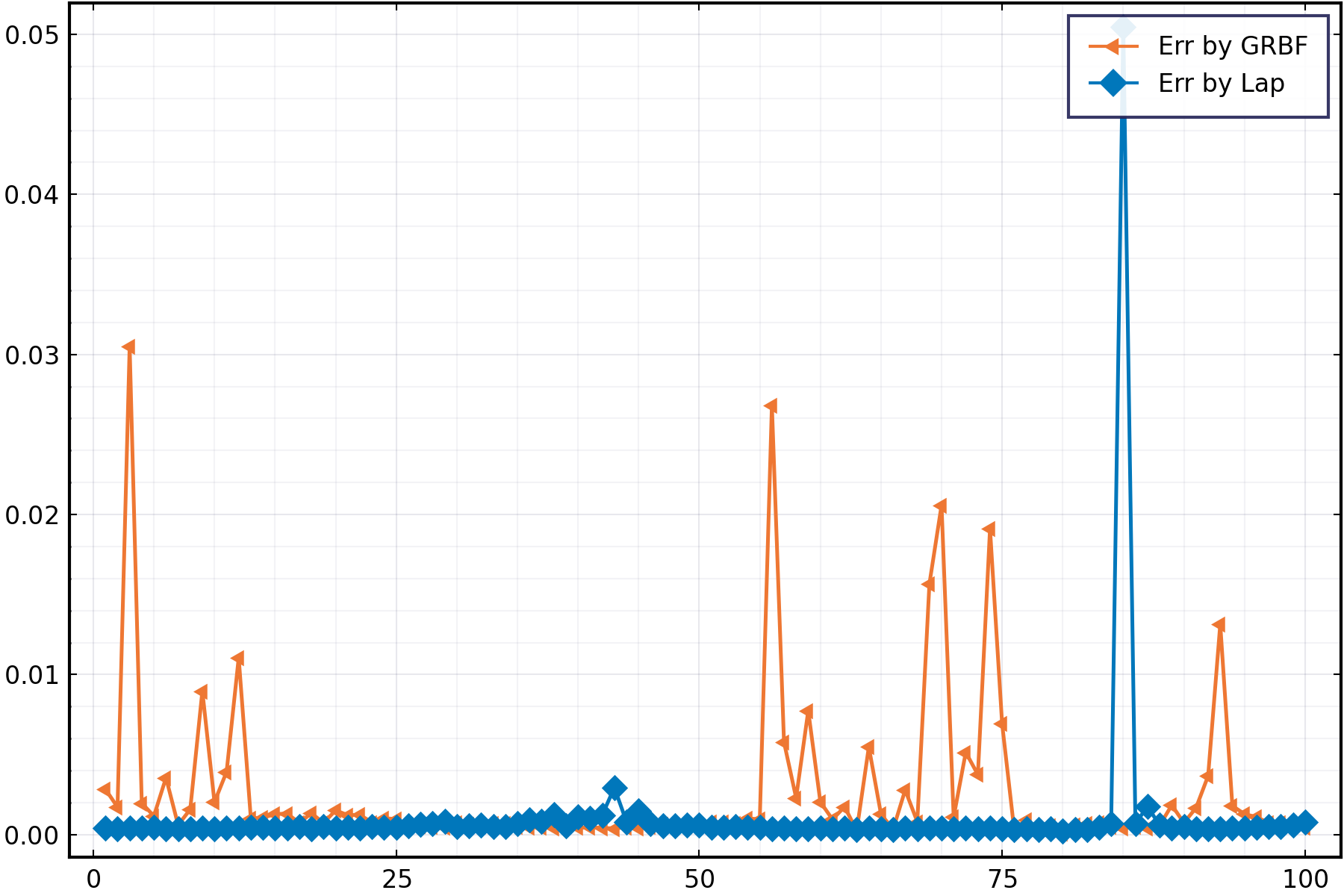}
        \caption{EWE plot}
    \end{subfigure}
    \caption{Reconstruction of values of NOAA Global Sea Surface Temperature anomalies and its comparison for the year $2023$.}
    \label{fig:WeatherSSTSpatialValues}
\end{figure}
\subsection{Experimental Inference}
We presented Lap-KeDMD to discover modal information from a given limited spatial-temporal data set and compared the results in contrast to those of the GRBF kernel using spectral observables of Koopman operators. From the spatial-temporal reconstructions and also with the error plots of every experiment through both $\Lap$ and $\GRBF$, the results obtained by the former are indeed superior to those of the latter. Hence, the spatial temporal reconstruction by the $\GRBF$ under the \Cref{aim} is poor.

This leads to raising the question of the applicability of the Laplacian kernel in modeling of data in contrast to what the GRBF kernel can offer. One possible argument, in light of the appropriate data modeling through the GRBF kernel being \emph{not so appropriate choice} that can be presented here, is due to \cite[Chapter 2, \S Section 2.6]{stein2012interpolation} and \cite[Chapter 4, Page 83]{rasmussen2006gaussian} being in terms of \emph{mean squared derivative} of the GRBF kernel, and hence the kernel function is very smooth. Unfortunately, the same scenario fails to exist for the Laplacian kernel which arose from the special class of kernel function called \emph{Mat\'ern} class of kernel function (see next section).

From the operator-theoretic perspective, the other promising direction, which also remains unexplored is the relationship between the Koopman operators with these kernel functions in the light of their corresponding RKHS. 
The whole situation around the choice of kernel function with respect to the operator-theoretic action of the Koopman operators demands \emph{closability} of it in the respective RKHS \cite{ikeda2022koopman,COLBROOK2024127} or \cite[Page 44]{colbrook2023multiverse}. We aim to now explore this theoretic endeavor in the next section by embedding the Laplacian kernel in the traditional $L^2$-measure sense.
\section{Laplacian kernel as \texorpdfstring{$L^2-$}{}measure}
\subsection{Preliminaries to the Laplacian kernel}
Let $\mathcal{X}\subset\mathbb{R}^D$ be an arbitrary input domain and let $\vx,\vx'\in\mathcal{X}$. Then, any function of $\bm{\tau}=\vx-\vx'$ is called a \emph{stationary} covariance function \cite[Chapter 4, Page 79]{rasmussen2006gaussian}. When the covariance function is only a function of $r=\|\vx-\vx'\|$, it is called \emph{isotropic}. Both the Laplacian kernel and the GRBF kernel are stationary covariance functions. The stationary covariance function can be represented as the Fourier transform of a certain positive finite measure. This representation is given by Bochner's theorem \cite{rasmussen2006gaussian,stein2012interpolation}.
\begin{theorem}[Bochner's Theorem]\label{theorem_Bochner'sTheorem}
    A complex-valued measurable function $k$ defined on $\mathbb{R}^D$ is the covariance function of a weakly stationary mean square continuous complex valued random process on $\mathbb{R}^D$ if and only if it can be represented as 
    \begin{align}
        k\left(\bm{\tau}\right)=\int_{\mathbb{R}^D}e^{2\pi i\mathbf{s}\cdot\bm{\tau}}d\mu(\mathbf{s}),\label{eq_Bochner'stheorem}
    \end{align}
    where $\mu$ is a positive finite measure.
\end{theorem}
If the positive finite measure $\mu$ in \eqref{eq_Bochner'stheorem} has a density $S(\mathbf{s})$ then $S$ is known as \emph{power spectrum} or \emph{spectral density} with respect to the kernel $k$. When the covariance function is isotropic (function of $r=\|\vx-\vx'\|_2$), then $S(\mathbf{s})$ is a function of $s=\|\mathbf{s}\|_2$ only \cite[Theorem 2.5.2]{adler2010geometry}. To gain some intuition about the construction in \cref{theorem_Bochner'sTheorem}, we illustrate the same for both $k=K_{\text{exp}}^{1,\sigma}(\bm{x},\bm{z})$ and $k=K_{\text{exp}}^{2,\sigma}(\bm{x},\bm{z})$, which we already introduced in the beginning.
\begin{enumerate}
    \item GRBF kernel $K_{\text{exp}}^{2,\sigma}(\bm{x},\bm{x'})=k_{\textsc{se}}(r)\circ\left(\|\vx-\vx'\|_2\right)$, where $k_{\textsc{se}}(r)=\exp(-\frac{r^2}{\sigma})$ is the squared exponential (SE) covariance function. The spectral density of the SE covariance function is $S(s)=(\sigma\pi)^{\frac{D}{2}}\exp(-\pi^2\sigma^2)$.
    \item For Laplacian kernel function, we need to first recall the Mat\'ern class of covariance function $k_{\textsc{mat\'ern}}(r)$; mathematically given as follows:
    \begin{align}\label{eq_Matern}
        k_{\textsc{mat\'ern}}(r)=\frac{2^{1-\nu}}{\Gamma(\nu)}\left(\frac{\sqrt{2}\nu r}{\sigma}\right)^\nu K_{\nu}\left(\frac{\sqrt{2}\nu r}{\sigma}\right)
    \end{align}
    where $\nu,\sigma>0$ and $K_\nu$ is a modified Bessel function \cite[Section 9.6]{abramowitz1965handbook}. When $\nu=\frac{1}{2}$ in \eqref{eq_Matern}, we get the exponential covariance function $k_{\textsc{exp}}(r)=\exp(-r/\sigma)$. 
    
    Laplacian kernel $K_{\text{exp}}^{1,\sigma}(\bm{x},\bm{x'})=k_{\textsc{exp}}(r)\circ\left(\|\vx-\vx'\|_2\right)$, where $k_{\textsc{exp}}(r)$ is the exponential covariance function coming from the Mat\'{e}rn class of covariance function. The spectral density of this is given as 
    \begin{align*}
        S(s)=\frac{2^D\pi^{\frac{D}{2}}\Gamma(1/2+D/2)}{\Gamma(1/2)\sigma}\left(\frac{1}{\sigma^2}+4\pi^2s^2\right)^{-(1/2+D/2)}
    \end{align*}
    in $D$ dimensions. For the single-dimension case, $k_{\textsc{exp}}(r)$ is the covariance function of the Ornstein-Uhlenbeck (OU) process \cite{uhlenbeck1930theory}, which is the mathematical model of the velocity of a particle undergoing Brownian motion.
\end{enumerate}
\subsection{Gaussian Measure and GRBF kernel}\label{subsection_GRBF}
The embedding of Laplacian kernel as in the traditional Lebesgue measure $L^2-$ is deeply inspired by the work of embedding function $k_{\textsc{se}}(r)$ coming from the GRBF kernel as the measure $L^2-$, which is referred to as the Gaussian measure $d\mu_{\mathds{G}}(\bm{z})\coloneqq\left(\sigma\pi\right)^{-D}\exp\left(-\|\bm{z}\|_2^2/\sigma\right)dV_D(\bm{z})$ 
\cite{bargmann1961hilbert,segal1963mathematical,baez2014introduction}. 

\begin{figure}[H]
    \centering
    \includegraphics[width=.75\linewidth]{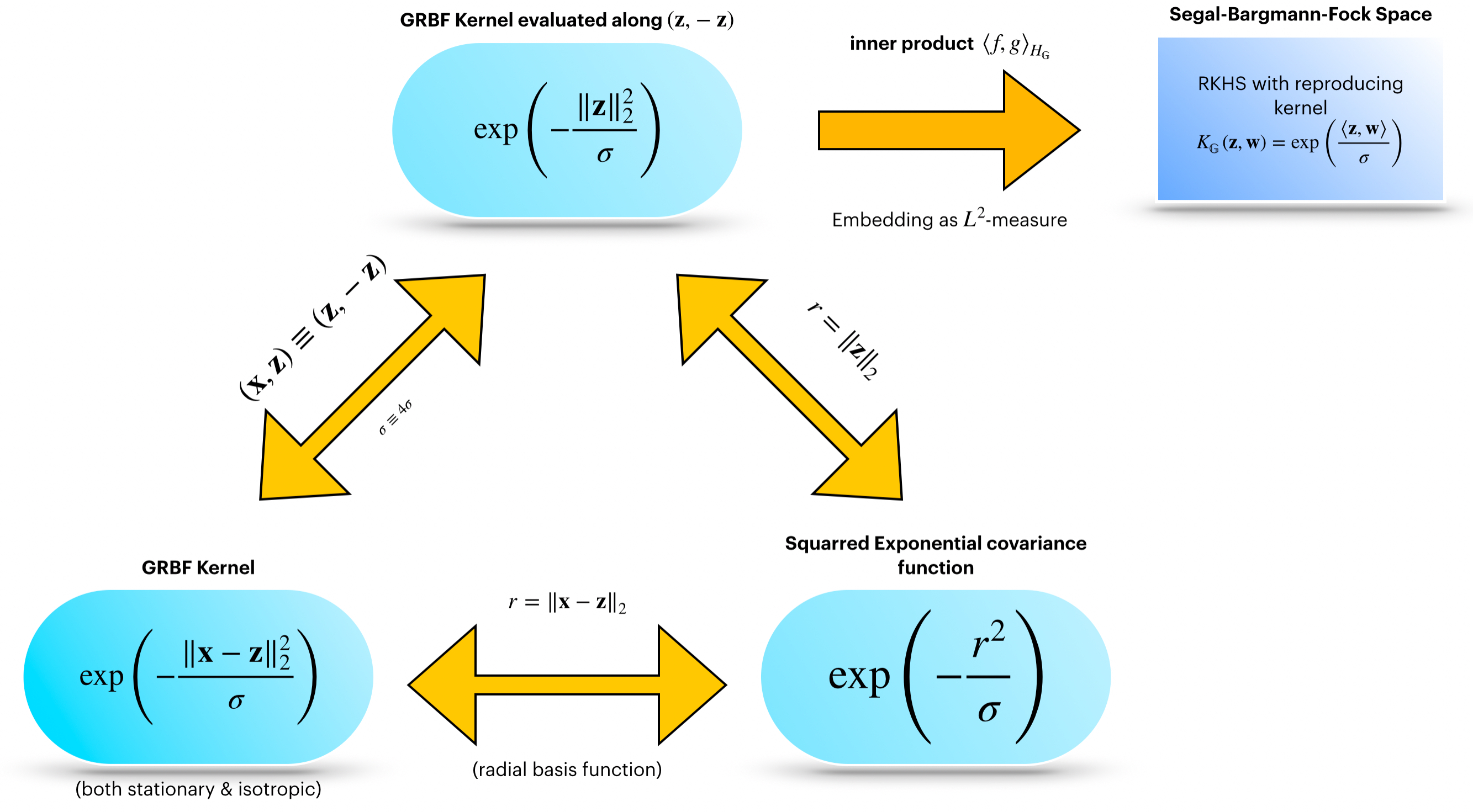}
    \caption{Schematic presentation of important aspects of GRBF kernel in relation with SE covariance function and the Gaussian measure. GRBF kernel evaluated along $(\bmz,-\bmz)$ provides the GRBF measure which leads to RKHS Segal-Bargmann-Fock Space with exponential dot product kernel as the reproducing kernel.}
    \label{fig:GRBF_Kernel}
\end{figure}

Consider an injective linear operator $\mathfrak{I}_{-}$ over $\Cn$ as $\mathfrak{I}_{-}\bmz=-\mI_{D}\bmz$. Then we can immediately observe that when the GRBF kernel is evaluated along the graph of operator $\mathfrak{I}_{-}$, that is, $\left(\bmz,\mathfrak{I}_{-}\bmz\right)=\left(\bmz,-\bmz\right)$ is actually the Gaussian measure $d\mu_{\mathds{G}}(\bm{z})$ up to the scaling of parameter $\sigma$; see \cref{fig:GRBF_Kernel}. 

Let the Hilbert space corresponding to the measure `$d\mu_{\mathds{G}}$' be denoted by $H_{\mathds{G}}$. It was introduced to comprehend the mathematical physics by \emph{Irving Segal} and \emph{Valentine Bargmann} and formally, its inner product is given as follows:
\begin{align}
    \langle f,g\rangle_{H_{\mathds{G}}}\coloneqq
    \int_{\mathbb{C}^D}f(\bm{z})\overline{g(\bm{z})}d\mu_{\mathds{G}}(\bm{z}),
\end{align}
followed by 
\begin{align}
    H_{\mathds{G}}\coloneqq\left\{\text{holomorphic}~f:\mathbb{C}^D\to\mathbb{C}\mid\|f\|_{H_{\mathds{G}}}\coloneqq\sqrt{\langle f,f\rangle_{H_{\mathds{G}}}}<\infty\right\}.
\end{align}
Hilbert space $H_{\mathds{G}}$ is referred as the \emph{Segal-Bargmann-Fock Space} \cite{zhu2012analysis,janson1987hankel} which indeed is a RKHS whose reproducing kernel is exponential dot product kernel given as
\begin{align}
    K_{\mathds{G}}(\bm{z},\bm{w})=\exp\left(\frac{\langle\bm{z},\bm{w}\rangle}{\sigma}\right).
\end{align}
\subsection{Embedding Laplacian kernel as \texorpdfstring{$L^2-$}{} measure}\label{subsec_LapKernelL2}
We plan to treat similarly for the Laplacian kernel and then construct the corresponding RKHS. Once we embed the Laplacian kernel as a $L^2-$measure, we will get the orthonormal basis which, further upon the application of Moore-Aronszajn Theorem (see \cref{theorem_aronsjan}) will provide the desired RKHS. We are now ready to start the analysis.

Let `$d\mu_{\mathds{L}}(\bm{z})$' be the Laplacian measure yielded by the Laplacian kernel $K_{\text{exp}}^{1,\sigma}(\vx,\bmz)$ when evaluated along $\left(\bmz,-\bmz\right)$ up to the scaling of the parameter $\sigma$. Explicitly it is given as follows:
\begin{align}
    d\mu_{\mathds{L}}(\bm{z})\coloneqq\left(2\pi\sigma^2\right)^{-D}\exp\left(-\frac{\|\bm{z}\|_2}{\sigma}\right)dV_D(\bm{z}).
\end{align}
We provide the inner product between two holomorphic functions $f:\Cn\to\mathbb{C}$ and $g:\Cn\to\mathbb{C}$ associated with the Laplacian measure $d\mu_{\mathds{L}}(\bm{z})$ as:
\begin{align}\label{eq_normLaplacemeasure}
    \langle f,g\rangle_{H_{\mathds{L}}}\coloneqq
    \int_{\mathbb{C}^D}f(\bmz)\overline{g(\bmz)}d\mu_{\mathds{L}}(\bm{z}),
\end{align}
followed by 
\begin{align}
    H_{\mathds{L}}\coloneqq\left\{\text{holomorphic}~f:\mathbb{C}^D\to\mathbb{C}\mid\|f\|_{H_{\mathds{L}}}\coloneqq\sqrt{\langle f,f\rangle_{H_{\mathds{L}}}}<\infty\right\}.
\end{align}
\subsubsection{Reproducing Kernel}
The result of the above theorem provides us with the orthonormal basis for the RKHS $\Hl$. Once we have the information for the orthonormal basis of the RKHS, we can employ the Moore-Aronszajn Theorem to construct the reproducing kernel for the associated RKHS. Therefore, in light of determining the reproducing kernel for $\Hl$, the following theorem provides the same for $\Hl$ using the Moore-Aronszajn Theorem.
\begin{theorem}\label{theorem_RKviaOrthonormalbasis}
    The reproducing kernel for RKHS $\Hl$ is 
    \begin{align}\label{eq_21}
        K_{\bm{w}}^\sigma\left(\bm{z}\right)=K^\sigma\left(\bm{z},\bm{w}\right)\coloneqq\left[\sqrt{\frac{\langle\bm{z},\bm{w}\rangle}{\sigma^2}}\right]^{-1}\cdot\sinh{\left(\sqrt{ \frac{\langle\bm{z},\bm{w}\rangle}{\sigma^2}}
        \right)}.
    \end{align}
\end{theorem}
The norm of the reproducing kernel $K_{\bmz}^\sigma$ is $\|K_{\bmz}^\sigma\|=\sqrt{\langle K_{\bmz}^\sigma,K_{\bmz}^\sigma\rangle_{\Hl}}=\sqrt{K^\sigma(\bmz,\bmz)}$. Now, we will provide the definition of \emph{normalized reproducing kernel} for the RKHS $\Hl$.
\begin{definition}
    At $\bmz\in\Cn$ and let $K_{\bmz}^\sigma$ be the reproducing kernel for the RKHS $\Hl$, then the normalized reproducing kernel $\bm{k}_{\bmz}^\sigma$ is given as
    \begin{align}
        \bm{k}_{\bmz}^\sigma\coloneqq\|K_{\bmz}^\sigma\|^{-1}\cdot K_{\bmz}^\sigma.
    \end{align}
    The norm of the normalized reproducing kernel is $\|\bm{k}_{\bmz}^\sigma\|=1$.
\end{definition}
The following result shows that the Laplacian kernel $\Lap(\cdot,\bm{x})$ belongs to the RKHS $\Hl$, $\forall \bm{x}\in\Cn$.
\begin{theorem}\label{thm_LapKernelinsideHl}
   $\forall\bm{x}\in\Cn$, $\left\|\Lapsigma\left(\cdot,\bm{x}\right)\right\|_{\Hl}\leq3^{-D}\exp\left(\frac{\|\bm{x}\|_2}{\sigma}\right)$.
\end{theorem}
\begin{proof}
    We start by evaluating the norm of $\Lapsigma(\bm{z},\bm{x})$ with respect to a fixed $\bm{x}\in\Cn$ inside the RKHS $\Hl$ as follows:
    \begin{align*}
        \left\|\Lapsigma\left(\cdot,\bm{x}\right)\right\|_{\Hl}^2=&
        \int_{\Cn}\left|\exp\left(-\frac{\|\bmz-\bm{x}\|_2}{\sigma}\right)\right|^2d\mu_{\mathds{L}}(\bmz)\\
        =&\left(2\pi\sigma^2\right)^{-D}\int_{\Cn}\left|\exp\left(-\frac{2\|\bmz-\bm{x}\|_2}{\sigma}\right)\right|\exp\left(-\frac{\|\bmz\|_2}{\sigma}\right)dV_D(\bmz)\\
        \leq&\left(2\pi\sigma^2\right)^{-D}\int_{\Cn}\exp\left(-\frac{2}{\sigma}\|\bm{z}\|_2\right)\exp\left(\frac{2\|\bm{x}\|_2}{\sigma}\right)\exp\left(-\frac{\|\bmz\|_2}{\sigma}\right)dV_D(\bm{z}).
    \end{align*}
In the last step, we use the reverse triangle inequality ($\|\bmz-\bm{x}\|_2\geq\|\bm{z}\|_2-\|\bm{x}\|_2$). Converting above into $\bm{z}\in\Cn$, $D$-dimensional polar coordinate system yields:
    \begin{align*}
        \left\|\Lapsigma\left(\cdot,\bm{x}\right)\right\|_{\Hl}^2\leq&\left(2\pi\sigma^2\right)^{-D}\exp\left(\frac{2\|\bm{x}\|_2}{\sigma}\right)\left[\int_{0}^{2\pi}d\theta\int_0^\infty\exp(-3r/\sigma)rdr\right]^D\\=&\left(2\pi\sigma^2\right)^{-D}\exp\left(\frac{2\|\bm{x}\|_2}{\sigma}\right)\left[\frac{2\pi}{(-3/\sigma)^2\Gamma(2)}\right]^D\\
        =&\frac{1}{9^D}\exp\left(\frac{2\|\bm{x}\|_2}{\sigma}\right)<\infty\quad\text{at $\bm{x}\in\Cn$.}
    \end{align*}
    Therefore, taking the square root of the above, provides the desired result.
\end{proof}
\Cref{thm_LapKernelinsideHl} establishes an important observation that the Laplacian kernel $\Lapsigma(\cdot,\cdot)$ exists in the RKHS $\Hl$ and hence can be employed as a choice of observables for the data-driven Koopman operators.
\section{Koopman operator-theoretic quantification over the RKHS \texorpdfstring{$\Hl$}{}}
In this section, we provide only the crucial operator-theoretic quantification of the Koopman operator that will be used in the next immediate section; that is the \emph{compactness} of it over the RKHS $\Hl$. For that, we need to first recall the main operator-theoretic quantities that laid the foundational basis for the operator-theoretic behavior of the Koopman operators over the RKHS $\Hl$.
\begin{definition}
        Let $\mathcal{K}_\varphi$ be the Koopman operator induced by the holomorphic function $\varphi:\Cn\to\Cn$ acting on the RKHS $\Hl$ with the reproducing kernel $K_{\bmz}^{\sigma}$ at $\bmz\in\Cn$. We define
        \begin{align}\label{eq_48}
            \Pi_{\bmz}(\varphi;\sigma)
            \coloneqq\left[\|K_{\bmz}^\sigma\|^{-2}\cdot\|\mathcal{K}_\varphi^*K_{\bmz}^\sigma\|^2\right].
        \end{align}
        We also define the supremum of above over $\bmz\in\Cn$ as follows:
        \begin{align}\label{eq_49r}
            \Pi\left(\varphi;\sigma\right)
            \coloneqq\sup_{\bmz\in\Cn}\Pi_{\bmz}(\varphi;\sigma).
        \end{align}
    \end{definition}
    \begin{note}\label{supplementary}
    The proofs for \emph{boundedness, essential norm estimates} and finally \emph{compactness} of the Koopman operators over the RKHS $\Hl$ are settled in the Supplementary Materials corresponding to the present article for better exposition.
    \end{note}
\begin{theorem}
\label{theorem_compact}
    A bounded $\koop$ achieves its compactification over RKHS $\Hl$ if and only if $\lim_{\|\bmz\|_2\to0}\Pi_{\bmz}(\varphi;\sigma)=0$ where $\varphi(\bmz)=\mA\bmz+\vb$ with $0\not\equiv\mA\in\mathbb{C}^{D\times D}$, with $0<\|\mA\|_{\text{Frob}}<1$ and $\mA$ is invertible.
    \end{theorem}
\section{Closable Koopman Operators}\label{review_closable}
\subsection{Review of closable operators}
We recall when we mean an operator $T$ in a Hilbert space $\mathfrak{H}$ as \emph{closable} or \emph{preclosed} as given in standard references \cite[Chapter 13]{rudin1991functional}, \cite[Chapter X, Page 304]{conway2019course} \cite[Chapter 5, Page 193]{pedersen2012analysis}, or \cite[Chapter VIII, Page 250]{reed2012methods}. 
\begin{definition}[Graph of an operator]
    For an (unbounded) operator $T$ in Hilbert space $\mathfrak{H}$ with its domain $\mathcal{D}(T)$, we define the \emph{graph of $T$} in $\mathfrak{H}$ as follows:
    \begin{align}\label{eq_graphofT}
        \bm{\Gamma}(T)\coloneqq\left\{(x,Tx):x\in\mathcal{D}(T)\right\}.  
    \end{align}
\end{definition}
\begin{definition}
    Let $T_{\boxdot}$ and $T$ be operators in the Hilbert space $\mathfrak{H}$. Let $\bm{\Gamma}(T_{\boxdot})$ and $\bm{\Gamma}(T)$ be the respective graphs of $T_{\boxdot}$ and $T$ as defined in \eqref{eq_graphofT}. If $\bm{\Gamma}(T)\subset\bm{\Gamma}(T_{\boxdot})$, then $T_\boxdot$ is said to be an extension of $T$ and we write $T\subset T_\boxdot$ and equivalently if $T\subset T_\boxdot$ if and only if $\mathcal{D}(T)\subset\mathcal{D}(T_\boxdot)$ and $T_\boxdot\Lambda=T\Lambda$ for all $\Lambda\in\mathcal{D}(T)$.
\end{definition}
\begin{definition}
    An operator is \textbf{closable} if it has a closed extension.
\end{definition}
\begin{lemma}\label{lemma_closable_characterization}
    The operator $T$ in the Hilbert space $\mathfrak{H}$ is closable if and only if for each sequence $\left\{x_n\right\}_n\in\mathcal{D}(T)$ converging to $0$, the only accumulation point of $\left\{Tx_n\right\}_n$ is $0$.
\end{lemma}

\subsection{Closability Results}\label{Lap_Koopman_closable_proof}
\begin{theorem}
    \label{closable_Koopman_theorem}
    Let $\koop:\Hl\to\Hl$ acts compactly over the RKHS $\Hl$ with $\varphi(\bmz)=\mA\bmz+\vb$ with an invertible $\mA\in\mathbb{C}^{D\times D}$ and $0<\|\mA\|_{\text{Frob}}<1$.
    Then $\koop$ can be closable in the RKHS $\Hl$.
\end{theorem}
\begin{proofoutline*}
\Cref{lemma_closable_characterization} (cf. \cite[Chapter 5, Page 193]{pedersen2012analysis}) can be interpreted as: \emph{a linear operator $T:\mathcal{D}(T)\to\mathfrak{H}$ is closable if and only if for any arbitrary sequence $g_m, Tg_m\in\mathfrak{H}$ such that $\lim_{g\to\infty}g_m=0$ and $Tg_m\to h$, then $h=0$.} We will use \Cref{lemma_closable_characterization} to demonstrate that compact Koopman operators $\koop$ induced by the holomorphic function $\varphi:\Cn\to\Cn$ are closable on the RKHS $\Hl$. For this, we will proceed by considering an arbitrary sequence $g_m\in\Hl$ such that $\lim_{m\to\infty}g_m=0$, then we will show that the action of the compact Koopman operator on these $g_m$ is in $\Hl$, and then finally, the limit of this is $0$ as well. However, the same cannot be established for the RKHS of the GRBF kernel.
\end{proofoutline*}
\begin{proof} 
Let $g_m$ be any arbitrary sequence defined on $\Hl$ which lives in $\Hl$, that is, $\|g_m\|_{\Hl}<\infty$ for all $m$. Suppose that for this arbitrary sequence, the accumulation point is $0$, that is, $\lim_{m\to\infty}g_m=0$. Allow the compact Koopman operator to act on $\bmz$ as follows:
\begin{align}
    \koop g_m(\bmz)=g_m(\varphi(\bmz))=g_m(\mA\bmz+\bm{b}).
\end{align}
Then the norm of above with respect to the RKHS $\Hl$ is 
\begin{align*}
    \left\|\koop g_m\right\|_{\Hl}^2=&\left\|g_m(\mA\bmz+\bm{b})\right\|_{\Hl}^2\\=&
    \int_{\Cn}|g_m\left(\mA\bmz+\bm{b}\right)|^2d\mu_{\mathds{L}}(\bmz)\\
    =&(2\pi\sigma^2)^{-D}
    \int_{\Cn}|g_m\left(\mA\bmz+\bm{b}\right)|^2\exp\left(-\frac{\|\bmz\|_2}{\sigma}\right)dV_d(\bmz)\\
    =&(2\pi\sigma^2)^{-D}\int_{\Cn}|g_m(\bm{u})|^2\exp\left(-\frac{\|\mA^{-1}(\bm{u-b})\|_2}{\sigma}\right)dV_D(\bm{u})(|\det(\mA)|)^{-1}\\
    \leq&(2\pi\sigma^2)^{-D}(|\det(\mA)|)^{-1}\exp\left(-\frac{\|\mA\bm{b}\|_2}{\sigma}\right)\int_{\Cn}|g_m(\bm{u})|^2\exp\left(-\frac{\|\mA\bm{u}\|_2}{\sigma}\right)dV_D(\bm{u})\\
    =&
    (|\det(\mA)|)^{-2}\exp\left(-\frac{\|\mA\bm{b}\|_2}{\sigma}\right)\int_{\Cn}|g_m(\bm{u})|^2d\mu_{\mathds{L}}(\bm{u})\\
    =&
    (|\det(\mA)|)^{-2}\exp\left(-\frac{\|\mA\bm{b}\|_2}{\sigma}\right)\|g_m\|_{\Hl}^2\\
    <&\infty.
\end{align*}
This means that the sequence $g_m\in\Hl$ after acting by the compact Koopman operator $\koop$ is still in the RKHS $\Hl$. Furthermore, since $\mA$ is invertible, the chain of inequalities shown above still holds.

\begin{align}
    \lim_{m\to0}\koop g_m(\bmz)=0=\lim_{m\to0}g_m(\mA\bmz+\bm{b})=\lim_{m\to\infty}\langle g_m,K_{\mA\bmz+\bm{b}}^\sigma\rangle_{\Hl}=0.
\end{align}
Therefore, the compact Koopman operator is closable over the RKHS $\Hl$ if and only if the symbol of the Koopman operator is \emph{affine}, that is, $\varphi(\bmz)=\mA\bmz+\bm{b}$ with invertible $\mA$ satisfying $0<\|\mA\|_{\text{Frob}}<1$.
\end{proof}
\subsubsection{Closability failure of Koopman operators over the RKHS of GRBF kernel}\label{GRBF_Koopman_failure_closable_proof}
In order to discuss the closability of Koopman operators on the GRBF kernel, we recall it as follows subsequently followed by its RKHS as well.
\begin{align*}
    K_{{exp}}^{2,\sigma}\triangleq K_{{exp}}^{2,\sigma}(\bm{x},\bm{z})&\coloneqq\exp\left(-\frac{\|\bm{x}-\bm{z}\|_2^2}{\sigma}\right).
\end{align*}
Norm corresponding to the Hilbert function space for the GRBF kernel \cite{steinwart2006explicit,steinwart2008support} is the reproducing kernel is 
\begin{align}\label{eq_6GRBFfunctionspace}
     \left\|f\right\|_{\sigma}^2\coloneqq\frac{2^D\sigma^{2D}}{\pi^D}\int_{\mathbb{C}^D}|f(\bm{z})|^2e^{\sigma^2\sum_{i=1}^D\left|\operatorname{Im}(z_i)\right|^2}dV_D(\bm{z}),\quad\bmz\coloneqq(z_1,z_2,\ldots,z_D)\in\Cn,
\end{align}
where $\operatorname{Im}(z_i)$ denotes the imaginary part of $z_i\in\mathbb{C}$. 
RKHS for kernel $K_{{exp}}^{2,\sigma}(\bm{x},\bm{z})$ is:
\begin{align}\label{eq_7GRBFHilbertspace}
    H_{\sigma}\coloneqq\left\{f:\mathbb{C}^D\to\mathbb{C}:\text{$f$ is holomorphic and $\|f\|_{\sigma}<\infty$}\right\}.
\end{align}
The following theorem uses the bounded Koopman operators in the RKHS $H_\sigma$ defined in \eqref{eq_7GRBFHilbertspace} whose result is borrowed from \cite[Corollary 1, Page 8]{gonzalez2024the}.
\begin{theorem}\label{theorem_GRBF_Koop_closable_failure}
Consider the Koopman operator $\koop:H_{\sigma}\to H_{\sigma}$ with symbol $\varphi(\bmz)=\mA\bmz+\bm{b}$. Then for any sequence $g_m\in H_{\sigma}$ such that $\lim_{m\to\infty}g_m=0$, we have $\|\koop g_m\|_{H_\sigma}\not<\infty$ which implies that $\koop g_m\not\in H_{\sigma}$ and therefore $\koop$ is not closable.
\end{theorem}
\begin{proof}
We have to incorporate a different style of methodology, since at this point we are not aware whether $\mA$ is invertible or not. In order to establish the proof, we start by considering $g_m\in H_{\sigma}$ such that $\lim_{m\to\infty}g_m=0$. Then we determine the norm of $\koop g_m$ for all $m$ with respect to the RKHS $H_{\sigma}$ as follows:
\begin{align*}
    \left\|\koop g_m\right\|_{H_\sigma}^2=&\frac{2^D\sigma^{2D}}{\pi^D}\int_{\mathbb{C}^D}|g_m(\mA\bmz+\bm{b})|^2\exp\left(\sigma^2\sum_{i=1}^D|\operatorname{Im}(z_i)|^2\right)dV_D(\bmz)\\
    =&\frac{2^D\sigma^{2D}}{\pi^D}\int_{\mathbb{C}^D}|\langle g_m,\GRBFsigma(\cdot,\mA\bmz+\bm{b})\rangle|^2\exp\left(\sigma^2\sum_{i=1}^D|\operatorname{Im}(z_i)|^2\right)dV_D(\bmz)\\
    \leq&\frac{2^D\sigma^{2D}}{\pi^D}\|g_m\|^2\int_{\mathbb{C}^D}\|\GRBFsigma(\cdot,\mA\bmz+\bm{b})\|^2\exp\left(\sigma^2\sum_{i=1}^D|\operatorname{Im}(z_i)|^2\right)dV_D(\bmz)\\
    =&\frac{2^D\sigma^{2D}}{\pi^D}\|g_m\|^2\int_{\mathbb{C}^D}\exp\left(\sigma^2\sum_{i=1}^D|\operatorname{Im}(z_i)|^2\right)dV_D(\bmz)\\
    \simeq&\frac{2^D\sigma^{2D}}{\pi^D}\|g_m\|^2\int_{\mathbb{C}^D}\mathcal{O}\left(\sum_{i=1}^D|\operatorname{Im}(z_i)|\right)\\
    \not<&\infty.
\end{align*}
Since the norm of $\koop g_m$ is not finite, this implies that $\koop g_m\not\in H_{\sigma}$. 
Failure of this concludes that the Koopman operator is not closable on the RKHS $H_{\sigma}$.
\end{proof}
\section{Theoretical implication of closability of Koopman operators for the dynamical system}
This section sits under the assumption of choosing the observables which lie in both $\LOmega$ and $\Hl$. Taking into account this assumption allows us to draw useful conclusions that tie the theory of spatial temporal reconstruction through the Koopman mode decomposition and the closability of the Koopman operators over the RKHS $\Hl$. We first provide a set-theoretic structure for \Cref{aim} as follows.
\subsection{On irregular and sparse data snapshots}
For the system \eqref{eq_1}, let $\left\{\mathbf{x}_i\right\}_{i=1}^N\in\Omega$ be a total of $N$ snapshots of data. Let 
\begin{align}
    \mathcal{S}^N\coloneqq\left\{\left\{\mathbf{x}_i\right\}_{i=1}^N\in\Omega:\text{each}~\mathbf{x}_{i}\in\Omega, \forall i=1,\ldots,N~\text{collected uniformly}\right\}.
\end{align}
Then, with respect to the uniform collection of data snapshots, we definte the following set which is the irregular and sparse variant of the same:
\begin{align}
    \mathcal{S}_{\Game}^N\coloneqq\left\{\left\{\mathbf{x}_i\right\}_{i=1}^N\in\Omega:\text{each}~\mathbf{x}_{i}\in\Omega, \forall i=1,\ldots,N~\text{collected irregularly and sparsely}\right\}.
\end{align}
It is immediate now that the following is true:
    $\mathcal{S}^N\cong\mathcal{S}_{\Game}^N,$
means that the entries of both sets are the same. In the light of our \Cref{aim}, we define the subset of $\mathcal{S}_{\Game}^N$ as follows for some $\mathbb{N}\ni N_1<N$:
\begin{align}
    \mathcal{S}_{\Game}^{N_1<N}\coloneqq
    \left\{\left\{\mathbf{x}_i\right\}_{i=1}^{N_1}\in\Omega:\text{each}~\mathbf{x}_{i}\in\Omega, \forall i=1,\ldots,N_1~\text{collected irregularly \& sparsely}\right\}.
\end{align}
Clearly, the following is true:
    $\mathcal{S}_{\Game}^{N_1<N}\subset\mathcal{S}_{\Game}^{N}.$
The above set-theoretic setup will be used in faithful spatial-temporal reconstruction in the upcoming \Cref{theorem_partialsum}.
\subsection{Theorems for theoretical implication of Koopman operators}
In this section, we provide the theoretical implication of the closability of Koopman operators corresponding to a dynamical system \eqref{eq_1}. We recall the spectral measures and the diagonalization of the Koopman operators of \Cref{subsection_KoopmanSpectralMeasure} and our closability result of the compact Koopman operators of \Cref{closable_Koopman_theorem} under a special condition, and we have the following theorem.
\begin{theorem}\label{theorem_spectralapprox}
    Let $g\in L^2(\Omega)\cap\Hl$ and $\varphi:\Cn\to\Cn$ be a holomorphic function. Let $\mathcal{K}_{\varphi_\mA}:\Hl\to\Hl$ be the compact closable Koopman operator in the RKHS $\Hl$ with $\varphi_{\mA}(\bmz)\coloneqq\mA\bmz$ where $\mA\in\mathbb{C}^{D\times D}$ and $0<\|\mA\|_{\text{Frob}}<1$ and $\mathcal{K}_{\mathbf{F}_t}:L^2(\Omega)\to L^2(\Omega)$ be the Koopman operator corresponding to the dynamics in \eqref{eq_1}.
    Then
    \begin{align}
        \|g\||\mathds{V}(\rho\left(\mathcal{K}_{\mathbf{F}_t}\right))|\sqrt{D+1-2\operatorname{Tr}(\mA)}<_{(1)}\|\mathcal{K}_{\mathbf{F}_t}g-\mathcal{K}_{\varphi_\mA}g\|<_{(2)}\|g\|\sqrt{D+1},
    \end{align}
    where $\operatorname{Tr}(\mA)\coloneqq\sum_{i=1}^Da_{ii}$ is the trace of $\mA$. Define $\epsilon(D,\mA)\coloneqq\frac{D+1}{D+1-2\operatorname{Tr}(\mA)}$ for $D\gg1$. Choosing a positive arbitrary $\epsilon$ such that $\epsilon>\epsilon(D,\mA)$, then
    \begin{align}
        \|\mathcal{K}_{\mathbf{F}_t}-\mathcal{K}_{\varphi_\mA}|_{\Hl}\|<\epsilon.\label{eq_7.3}
    \end{align}
\end{theorem}
\begin{proof} 
We will first establish inequality (2) as follows. Observe that if $g\in L^2(\Omega)$ with $\mathcal{K}_{\mathbf{F}_t}:L^2(\Omega)\to L^2(\Omega)$, then
    \begin{align}
        \mathcal{K}_{\mathbf{F}_t}g(\mathbf{x})=\left[\left(\int_{\mathbb{D}}\lambda d\mathds{V}(\lambda)\right)g\right](\mathbf{x}).
    \end{align}
    Now, if $g\in\Hl$ with closable $\koop:\Hl\to\Hl$ where $\varphi(\bmz)=\mA\bmz$ and $0<\|\mA\|_{\text{Frob}}<1$, due to \Cref{closable_Koopman_theorem} and then
    \begin{align}
        \mathcal{K}_{\varphi_\mA}g(\mathbf{x})=&\left[\left(\int_{\mathbb{D}}\lambda d\mathds{V}(\lambda)\right)g\right]\left(\mA\mathbf{x}\right)\label{eq_7.5}.
    \end{align}
    Then,
    \begin{align*}
        \left\|\mathcal{K}_{\mathbf{F}_t}g(\mathbf{x})-\mathcal{K}_{\varphi_\mA}g(\mathbf{x})\right\|^2=&\left\|\left[\left(\int_{\mathbb{D}}\lambda d\mathds{V}(\lambda)\right)g\right](\mathbf{x})-\left[\left(\int_{\mathbb{D}}\lambda d\mathds{V}(\lambda)\right)g\right](\mA\mathbf{x})\right\|^2\\
        =&\left\|\left[\left(\int_{\mathbb{D}}\lambda d\mathds{V}(\lambda)\right)g\right]\left(\mathbf{x}-\mA \mathbf{x}\right)\right\|^2\numberthis\label{eq_7.7}\\
        \leq&\left\|\left[\left(\int_{\mathbb{D}}\lambda d\mathds{V}(\lambda)\right)g\right]\right\|\left\|\mathbf{x}-\mA \mathbf{x}\right\|^2\\
        \leq&\|g\|^2\left\|\int_{\mathbb{D}}\lambda d\mathds{V}(\lambda)\right\|^2\|\mathbf{x}\|^2\|\mathbf{I}_D-\mA\|_{\text{Frob}}^2\\
        \leq&\|g\|^2\left\|\int_{\mathbb{D}}|\lambda| |d\mathds{V}(\lambda)|\right\|^2\|\mathbf{x}\|^2\|\mathbf{I}_D-\mA\|_{\text{Frob}}^2\\
        =&\|g\|^2|\mathds{V}(\mathbb{D})|^2\|\mathbf{x}\|^2\|\mathbf{I}_D-\mA\|^2_{\text{Frob}}.
    \end{align*}
    Since $\|\mathbf{I}_D-\mA\|^2_{\text{Frob}}=\sum_{i=1}^D|1-a_{ii}|^2=\sum_{i=1}^D(1+|a_{ii}^2|-2a_{ii})=D+\|\mA\|^2_{\text{Frob}}-2\operatorname{Tr}(\mA)$, and also $|\mathds{V}(\mathbb{D})|=1$. Thus
    \begin{align}
        \left\|\mathcal{K}_{\mathbf{F}_t}g(\mathbf{x})-\mathcal{K}_{\varphi_\mA}g(\mathbf{x})\right\|^2\leq&\|g\|^2\|\mathbf{x}\|^2\left(D+\|\mA\|^2_{\text{Frob}}-2\operatorname{Tr}(\mA)\right)
        <\|g\|^2(D+1).
    \end{align}
    Choosing $\epsilon>\frac{D+1}{D+1-2\operatorname{Tr}(\mA)}$, we have the desired result. Now, we will establish inequality (1) as follows. For that recall \eqref{eq_7.7} and also $\rho(\mathcal{K}_{\mathbf{F}_t})\subset\mathbb{D}$, therefore $\mathds{V}\left(\rho(\mathcal{K}_{\mathbf{F}_t})\right)<\mathds{V}\left(\mathbb{D}\right)=1$, then
\begin{align*}
    \left\|\mathcal{K}_{\mathbf{F}_t}g(\mathbf{x})-\mathcal{K}_{\varphi_\mA}g(\mathbf{x})\right\|\geq&\left\|\left[\left(\int_{\rho(\mathcal{K}_{\mathbf{F}_t})}\lambda d\mathds{V}(\lambda)\right)g\right]\left(\mathbf{x}-\mA\mathbf{x}\right)\right\|^2\\
    =&\|g\|^2|\mathds{V}(\rho(\mathcal{K}_{\mathbf{F}_t}))|^2\|\mathbf{x}\|^2\left(D+1-2\operatorname{Tr}(\mA)\right).
\end{align*}
It should be noted that we assume that $\mathbf{x}\neq0$ in another way, which means that the result of this theorem will be redundantly true.
\end{proof}
\Cref{theorem_spectralapprox} shows the spectral approximation of the Koopman operators that arose from the dynamical system by the one that acts closable on the RKHS $\Hl$. This result is drawn under the assumption that we chose an observable that lives in both Hilbert space, that is, $\LOmega$ and $\Hl$. It is worth mentioning that one can choose $\epsilon$ arbitrarily large because we are assuming that the dimension of the state space $D$ is large in \eqref{eq_7.3} and therefore we can achieve such an approximation. 
\begin{figure}
    \centering
    \includegraphics[scale=1]{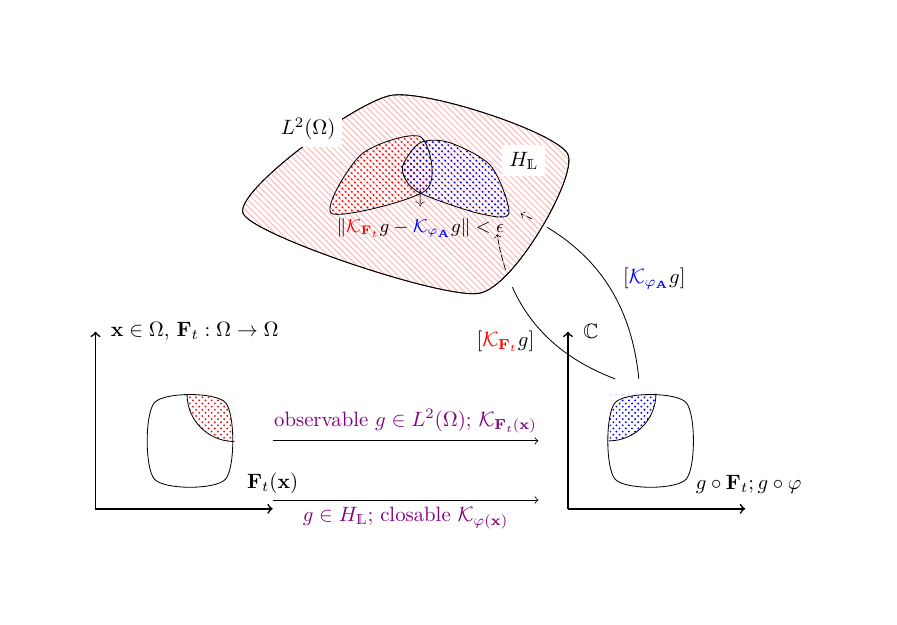}
    \caption{An illustration of the spectral measure approximations of the data-driven Koopman operators over $\LOmega$ due to the closability of the function theoretic Koopman operators over the RKHS $\Hl$. Note that, this is under an assumption that the operator-theoretic Koopman operator has symbol \emph{linear} in $\Cn$.}
    \label{fig:spectral_Closable}
\end{figure}

So far, we gave a theoretical justification in the form of \Cref{theorem_spectralapprox} that concerns only the spectral measures of the Koopman operators. One can draw the following result in terms of Koopman mode decomposition as well with a similar setup provided in the preceding theorem, which has more in-depth practical meaning for the spatial modes reconstruction.
\begin{theorem}\label{theorem_partialsum}
    Let $\mathbf{x}_m$ be the snapshot of the $m$\textsuperscript{th} data corresponding to the system \eqref{eq_1} drawn from $\mathcal{S}_{\Game}^{N_1<N}$. Let $\mathcal{K}_{\mathbf{F}_t}:\LOmega\to\LOmega$ be the corresponding Koopman operator, and let $\zeta_n$ be its $n$\textsuperscript{th}-eigenfunction, and let $\mathbf{c}_n$ be the Koopman modes of \eqref{eq_KoopmanModes}. Assume $\zeta_n\in\LOmega\cap\Hl$. Let the Koopman operator $\koop:\Hl\to\Hl$ induced by $\varphi:\Cn\to\Cn$ be closable over RKHS $\Hl$. 
    
    Then the faithful Koopman mode decomposition with respect to the RKHS $\Hl$ is given as 
    \begin{align}\label{eq_PartialSum}
        \digamma\Delta_{\mathcal{\Hl}}^{(0,m)}
        =\sum_{n=1}^\mathcal{\infty}\mathbf{c}_n\left[\mathcal{K}_{\mathbf{F}_t}^m\zeta_n(\mathbf{x}_0)-\zeta_n\left(\mA^m\mathbf{x}_0+\bm{b}\sum_{i=0}^{m-1}\mA^i\right)\right]
    \end{align} 
    where $\mA\in\mathbb{C}^{D\times D}$ and $0<\|\mA\|_{\text{Frob}}<1$. Lastly, 
    \begin{align}
        \lim_{\mathcal{N}\to\infty}\mathcal{P}_\mathcal{N}\digamma\Delta_{\mathcal{\Hl}}^{(0,m)}=\digamma\Delta_{\mathcal{\Hl}}^{(0,m)}\equiv0.
    \end{align}
\end{theorem}
\begin{proof}
    Let $\zeta_n\in\LOmega\cap\Hl$ for all $n\in\mathbb{N}$. Then recall the Koopman mode decomposition equation in \eqref{eq_KoopmanModes} for the tray of observables $\mathbf{g}$ as follows:
    \begin{align*}
        \mathbf{g}(\mathbf{x}_m)=\sum_{n\in\mathbb{N}}\mathbf{c}_n\mathcal{K}_{\mathbf{F}_t}^m\zeta_n(\mathbf{x}_0).
    \end{align*}
    Allowing the Koopman operator $\mathcal{K}_{\varphi}:\Hl\to\Hl$ to act on the RKHS $\Hl$ closable for the observable $\zeta_n\in\Hl$, then after recalling \Cref{closable_Koopman_theorem}, we have
    \begin{align*}
        \mathbf{g}(\mathbf{x}_m)=\sum_{n\in\mathbb{N}}\mathbf{c}_n\mathcal{K}_{\varphi}^m\zeta_n(\mathbf{x}_0)=\sum_{n\in\mathbb{N}}\mathbf{c}_n\zeta_n\left(\mA^m\mathbf{x}_0+\bm{b}\sum_{i=0}^{m-1}\mA^i\right).
    \end{align*}
    Equating both equations with each other yields the following:
    \begin{align*}
        \sum_{n\in\mathbb{N}}\mathbf{c}_n\left[\mathcal{K}_{\mathbf{F}_t}^m\zeta_n(\mathbf{x}_0)-\zeta_n\left(\mA^m\mathbf{x}_0+\bm{b}\sum_{i=0}^{m-1}\mA^i\right)\right]=0.
    \end{align*}
    Since, the above infinite series sums up to $0$ therefore its partial sum as already noted in \eqref{eq_PartialSum} is also $0$.
\end{proof}
\Cref{theorem_partialsum} implies that when we have accumulated a countable number of these faithful partial sums, the result leads to $0$. Therefore, we achieve the faithful reconstruction of the spatial modes through those observables coming from the RKHS $\Hl$. 
\section{Conclusion}
We present an impressive application of the Laplacian kernel in KeDMD to reconstruct spatial temporal reconstruction under \Cref{aim} for a variety of data sets in \Cref{tab:dataset}. We learn through these chain of experiments that, under identical experimental conditions, the results delivered by the Laplacian kernel outperform that of the GRBF kernel. 
In addition to our empirical study, we also explored the operator-theoretic paradigms of Koopman operators over the RKHS $\Hl$\footnote{follow \Cref{supplementary}.} generated by the Laplacian kernel as a $L^2-$measure. These operator-theoretic paradigms eventually led to the application of the Koopman operators in the RKHS $\Hl$ providing a concrete theoretical justification as well. Ultimately, from such an operator-theoretical investigation, we learn that RKHS $\Hl$ can provide the closability of Koopman operators and therefore we settle this challenge in this paper.

We also devised a new definition for the Koopman mode difference \eqref{eq_KoopmanModes} in \Cref{def_0msnapdiff}. We used the aforementioned definitions to establish a concrete link between closable Koopman operators over the RKHS $\Hl$ and data-driven Koopman operators over the Hilbert space $\LOmega$. In view of these results, we are now able to understand why the Laplacian kernel is capable of performing the spatial-temporal reconstruction in light of \Cref{aim} in contrast to the GRBF kernel. In conclusion, we establish empirical suitability of the Laplacian kernel and provide theoretical justification for the same. We also envision in light of \Cref{aim} that spatial temporal reconstruction is directly linked to the determination and existence of RKHS where the Koopman operators can be closable (cf. \Cref{def_0thmthdifference_faithful}, \Cref{closable_Koopman_theorem}, \Cref{theorem_spectralapprox} and \Cref{theorem_partialsum}). Since the RKHS that we devised in the present articles enjoys such a mathematical endeavor, the Laplacian kernel is an appropriate and preferable choice for extended DMD. 

For future directions, the practical side of \eqref{eq_PartialSum} with respect to RKHS $\Hl$ is yet to be seen. Additionally, so far, we have not encountered any literature in which the utility of the Laplacian kernel is explored and discussed to the extent that we have shown in our paper. Therefore, we envision the future of this in various other platforms ranging from data-driven science to modern cutting-edge ML regimes as well.

\end{document}